\documentclass[12pt,bezier]{article}
\usepackage{times}
\usepackage{booktabs}
\usepackage{enumitem}
\usepackage{pifont}
\usepackage{floatrow}
\usepackage{caption}
\usepackage{mathrsfs}
\usepackage[fleqn]{amsmath}
\usepackage{amsfonts,amsthm,amssymb,mathrsfs}
\usepackage{txfonts}
\usepackage{tikz}
\usepackage{graphics,multicol}
\usepackage{graphicx}
\usepackage{geometry}
\usepackage{xcolor}
\usepackage[colorlinks=true,anchorcolor=black,filecolor=black,linkcolor=black,urlcolor=black,citecolor=black]{hyperref}
\usepackage{caption}
\usepackage{tikz}
\usepackage{cite}

\usepackage{latexsym,bm}
\usepackage{indentfirst}
\usepackage{mathtools}
\usepackage{array,subcaption}
\usepackage{amsmath}
\evensidemargin=\oddsidemargin\topmargin=-1.5cm

\newtheorem{theorem}{Theorem}[section]
\newtheorem{fact}{Fact}
\newtheorem{lemma}[theorem]{Lemma}

\newtheorem{corollary}{Corollary}[section]
\theoremstyle{definition}

\addtocounter{section}{0}

\usetikzlibrary{decorations.pathreplacing}
\newcommand{\npos}{n_{+}}
\newcommand{\nzero}{n_{0}}
\newcommand{\one}{\mathbf{1}}

\newcommand{\mi}{\operatorname{mult}_{1}}
\newcommand{\diag}{\operatorname{diag}}

\usetikzlibrary{positioning,fit,calc,decorations.pathreplacing,shapes.geometric}
\hypersetup{
pdftitle={Distance-hereditary split graphs whose second largest eigenvalue does not exceed 1},
pdfauthor={Zihao Geng, Qian-Qian Chen, Ji-Ming Guo}
}

\usepackage{pdfpages}
\begin{document}
\title{Distance-hereditary split graphs whose second largest eigenvalue does not exceed $1$\footnote{This work was supported by the NSFC (Nos.~12171154; 12601690) and the Yancheng Applied Basic Research Program (No.~YCBK2026135).}}
\author{{\Large Zihao Geng$^{a}$,\,\, Qian-Qian Chen$^{b}$,\,\, Ji-Ming Guo$^{a}$\footnote{Corresponding author.}
\setcounter{footnote}{-1}
\footnote{\emph{E-mail addresses:} \mbox{gzh635516986@163.com} (Z.H. Geng), \mbox{qianqian\_chen@yeah.net} (Q.Q. Chen), \mbox{jimingguo@hotmail.com} (J.-M. Guo).}}\\[2mm]
\parbox[t]{0.96\textwidth}{%
\centering\fontsize{10}{12}\selectfont\itshape
\setlength{\parindent}{0pt}\setlength{\parskip}{4pt}
\textsuperscript{a}\,Department of Mathematics, East China University of Science and Technology,\\
Shanghai 200237, China\endgraf
\textsuperscript{b}\,School of Mathematics and Statistics, Yancheng Teachers University,\\
Yancheng 224002, China}}

\date{}
\maketitle

\vspace{\dimexpr\topsep+\partopsep\relax}
\noindent{\large\bfseries Abstract}\\
The gem is obtained by joining a vertex to every vertex of a path on four vertices.
A distance-hereditary split graph is a split graph containing no induced gem.
In this paper, we characterize all connected distance-hereditary split graphs whose second largest  eigenvalue does not exceed $1$ and  determine the multiplicity of the eigenvalue $1$.
Our result extends previous results on threshold graphs whose second largest eigenvalue does not exceed $1$.\\
\vspace{0.1cm}

\par\vspace{1.5\baselineskip}


\section{Introduction}

In this paper, all graphs are finite, simple and undirected.
Let $G$ be a graph with vertex set $V(G)$ and edge set $E(G)$.
For $U\subseteq V(G)$, let $G[U]$ denote the induced subgraph of $G$ induced by $U$.
For a graph $G$, let $\lambda_{1}(G)\geq\lambda_{2}(G)\geq\cdots\geq\lambda_{n}(G)$ be the eigenvalues of its adjacency matrix $A(G)$.
As usual, $P_{n}$, $C_{n}$ and $K_{n}$ denote the path, cycle and complete graph on $n$ vertices, respectively, and $K_{s,t}$ denotes the complete bipartite graph with parts of sizes $s$ and $t$.

The problem of characterizing connected graphs with $\lambda_{2}\leq1$ was posed by Hoffman and first studied by Cvetkovi\'c~\cite{Cvetkovic1982}.
For a graph $F$, a graph is $F$-free if it contains no induced subgraph isomorphic to $F$.
Let $G$ be a graph and $G-e$ denote the graph obtained by deleting an edge $e$ from $G$.
A set of $n$ pairwise adjacent vertices is called an $n$-clique.
This problem has been studied for several classes of graphs.
Trees~\cite{Neumaier1982}, unicyclic graphs~\cite{Xu2004}, bicyclic graphs~\cite{Guo2005}, tricyclic graphs~\cite{LiYang2011}, bipartite graphs~\cite{Petrovic1991}, line graphs~\cite{PetrovicMilekic1998}, generalized line graphs~\cite{PetrovicMilekic2000}, generalized $\theta$-graphs~\cite{GaoHuang2008},  graphs with three distinct eigenvalues~\cite{ChengGavrilyukGreavesKoolen2016,ChengGreavesKoolen2018}, $r$-regular graphs with $r\leq10$ or $n\geq2r$~\cite{KoledinStanic2013,Stanic2010,WangZhu2026}, $(K_{1,3},K_{5}-e)$-free graphs~\cite{LiuChenStanic2021}, $(K_{4}-e)$-free graphs~\cite{LiuChenStanic2024}, claw-free graphs~\cite{ChenShanLiuStanic2026}, non-bipartite graphs with girth at least $4$~\cite{XuShao2006},   $K_{4}$-minor-free graphs~\cite{Cheng2026} and planar graphs~\cite{GengWangLiLiChenGuo2026,LiSun2023,GengGuoWang2024,ChengLi2025} with $\lambda_{2}\leq1$ have been characterized.
In addition, Zhang et al.~\cite{ZhangLiangKoolen2022} classified $2$-partially distance-regular graphs whose local graphs  have second largest eigenvalue at most $1$.

A graph $G=(C,Z)$ is said to be split if its vertex set can be partitioned into a clique $C$ and an independent set $Z$.
Such a partition $(C,Z)$ is called a split partition of $G$.
Without loss of generality, we suppose that $C$ is a maximum clique of $G$ throughout this paper.
The gem is obtained by joining a vertex to every vertex of $P_{4}$.
A split graph $G=(C,Z)$ is distance-hereditary if and only if it is gem-free~\cite{Brandstadt2005}.
Let $N(x)$ denote the set of neighbors of $x$.
For two sets $X,Y$, we write $X\subseteq Y$ if every element of $X$ belongs to $Y$.
We write $X\subset Y$ if $X\subseteq Y$ and $X\neq Y$.
A split graph is called a nested split graph if $N(x_{1})\subseteq\cdots\subseteq N(x_{t})$, where $Z=\{x_{1},\ldots,x_{t}\}$.
Nested split graphs, also called threshold graphs, are precisely the graphs with no induced $2K_{2}$, $P_{4}$ or $C_{4}$~\cite{MilatovicStanic2012,Stanic2009}.
Every threshold graph is distance hereditary~\cite{ChangHungKloksPeng2011}.

In 2009, Stani\'c~\cite{Stanic2009} characterized all nested split graphs $G$ with $\lambda_{2}(G)<1$.
In 2012, Milatovi\'c and Stani\'c~\cite{MilatovicStanic2012} characterized all nested split graphs $G$ with $\lambda_{2}(G)=1$.
These two results give a complete characterization of threshold graphs with $\lambda_{2}(G)\leq1$.
Since threshold graphs form a proper subclass of distance-hereditary split graphs, the same problem for distance-hereditary split graphs is a natural extension.
In this paper, we characterize all connected distance-hereditary split graphs whose second largest eigenvalue does not exceed $1$ and determine the multiplicity of the eigenvalue $1$.

%Let $V(G)$ and $E(G)$ be the vertex set and edge set of $G$, respectively.
%If $G$ is a connected distance-hereditary split graph of order at least three and $V(G)=C\cup Z$, where $C$ is a maximum clique and $Z$ is an independent set, then for any $x,y\in Z$, $N(x)\cap N(y)=\varnothing$, $N(x)\subseteq N(y)$ or $N(y)\subseteq N(x)$ by~\cite{Cameron2016}.
%Since $G$ is connected, $N(x)\neq\varnothing$ for every $x\in Z$.
%Since $C$ is a maximum clique, $N(x)\subset C$ for every $x\in Z$.
%This family of neighborhoods is called laminar.

% Begin embedded drawing: figures/chain_families.tex
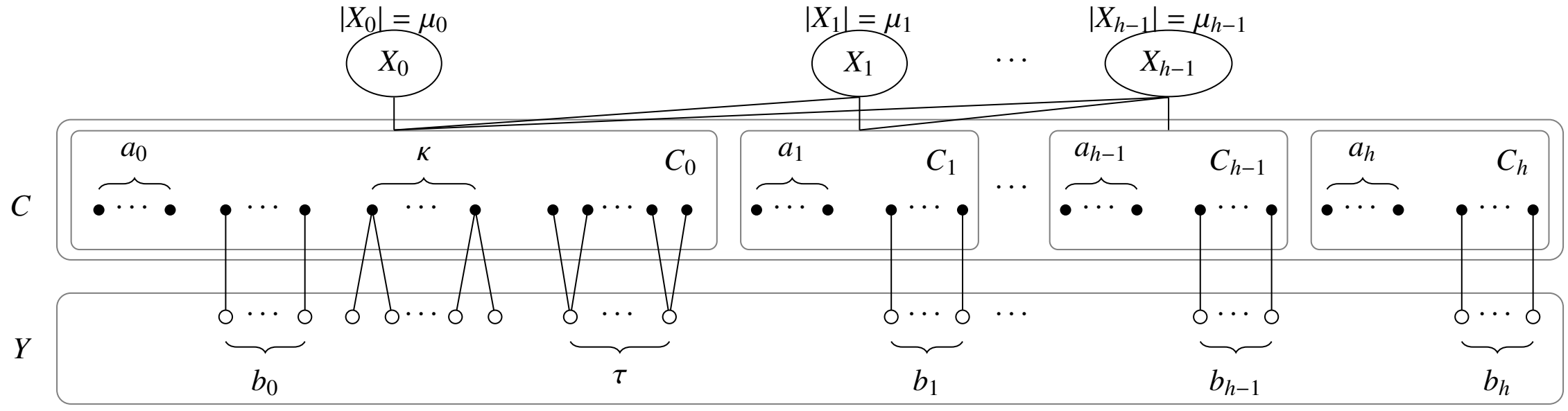
\begin{figure}[!htbp]
\centering
\resizebox{\textwidth}{!}{%
\begin{tikzpicture}[x=1cm,y=1cm,line width=0.5pt,every node/.style={font=\small}]
% All filled vertices belong to the same clique C; clique edges are omitted.
\draw[rounded corners=5pt,gray] (-0.18,-0.53) rectangle (18.83,1.24);
\node[anchor=east] at (-0.35,0.15) {$C$};
\draw[rounded corners=5pt,gray] (-0.18,-0.95) rectangle (18.83,-2.35);
\node[anchor=east] at (-0.35,-1.65) {$Y$};
\node[draw,rounded corners=3pt,gray,minimum width=81.5mm,minimum height=15mm,inner sep=0pt] (C0) at (4.075,0.35) {};
\node[anchor=north east,xshift=-3pt,yshift=-3pt] at (C0.north east) {$C_{0}$};
\node[draw,rounded corners=3pt,gray,minimum width=30mm,minimum height=15mm,inner sep=0pt] (C1) at (9.95,0.35) {};
\node[anchor=north east,xshift=-3pt,yshift=-3pt] at (C1.north east) {$C_{1}$};
\node[draw,rounded corners=3pt,gray,minimum width=30mm,minimum height=15mm,inner sep=0pt] (Chm) at (13.85,0.35) {};
\node[anchor=north east,xshift=-3pt,yshift=-3pt] at (Chm.north east) {$C_{h-1}$};
\node[draw,rounded corners=3pt,gray,minimum width=30mm,minimum height=15mm,inner sep=0pt] (Ch) at (17.15,0.35) {};
\node[anchor=north east,xshift=-3pt,yshift=-3pt] at (Ch.north east) {$C_{h}$};
\node at (11.9,0.35) {$\cdots$};
\node at (11.9,-1.25) {$\cdots$};
% Each X_i is completely adjacent to C_0,...,C_i.
\node[draw,ellipse,fill=white,minimum width=12mm,minimum height=7mm] (X0) at (4.075,1.95) {$X_{0}$};
\node at (4.075,2.5) {$|X_{0}|=\mu_{0}$};
\draw (X0.south)--(C0.north);
\node[draw,ellipse,fill=white,minimum width=12mm,minimum height=7mm] (X1) at (9.95,1.95) {$X_{1}$};
\node at (9.95,2.5) {$|X_{1}|=\mu_{1}$};
\draw (X1.south)--(C0.north);
\draw (X1.south)--(C1.north);
\node at (11.9,1.95) {$\cdots$};
\node[draw,ellipse,fill=white,minimum width=16mm,minimum height=7mm] (Xhm) at (13.85,1.95) {$X_{h-1}$};
\node at (13.85,2.5) {$|X_{h-1}|=\mu_{h-1}$};
\draw (Xhm.south)--(Chm.north);
\draw (Xhm.south)--(C0.north);
\draw (Xhm.south)--(C1.north);
% Two representatives of the a_0 type, with intervening ellipses.
\fill (0.35,0.1) circle (2pt);
\fill (1.25,0.1) circle (2pt);
\node at (0.8,0.1) {$\cdots$};
\draw[decorate,decoration={brace,amplitude=4pt}] (0.35,0.36)--(1.25,0.36) node[midway,above=6pt] {$a_{0}$};
% Two representatives of the b_0 type, with intervening ellipses.
\draw (1.95,0.1)--(1.95,-1.25);
\fill (1.95,0.1) circle (2pt);
\draw[fill=white] (1.95,-1.25) circle (2.4pt);
\draw (2.95,0.1)--(2.95,-1.25);
\fill (2.95,0.1) circle (2pt);
\draw[fill=white] (2.95,-1.25) circle (2.4pt);
\node at (2.45,0.1) {$\cdots$};
\node at (2.45,-1.25) {$\cdots$};
\draw[decorate,decoration={brace,mirror,amplitude=4pt}] (1.95,-1.56)--(2.95,-1.56) node[midway,below=6pt] {$b_{0}$};
% Two representatives of the kappa type, with intervening ellipses.
\draw (3.8,0.1)--(3.55,-1.25);
\draw[fill=white] (3.55,-1.25) circle (2.4pt);
\draw (3.8,0.1)--(4.05,-1.25);
\draw[fill=white] (4.05,-1.25) circle (2.4pt);
\fill (3.8,0.1) circle (2pt);
\draw (5.1,0.1)--(4.85,-1.25);
\draw[fill=white] (4.85,-1.25) circle (2.4pt);
\draw (5.1,0.1)--(5.35,-1.25);
\draw[fill=white] (5.35,-1.25) circle (2.4pt);
\fill (5.1,0.1) circle (2pt);
\node at (4.45,0.1) {$\cdots$};
\node at (4.45,-1.25) {$\cdots$};
\draw[decorate,decoration={brace,amplitude=4pt}] (3.8,0.36)--(5.1,0.36) node[midway,above=6pt] {$\kappa$};
% Two representatives of the tau type, with intervening ellipses.
\draw (6.08,0.1)--(6.3,-1.25)--(6.52,0.1);
\fill (6.08,0.1) circle (2pt);
\fill (6.52,0.1) circle (2pt);
\draw[fill=white] (6.3,-1.25) circle (2.4pt);
\draw (7.33,0.1)--(7.55,-1.25)--(7.77,0.1);
\fill (7.33,0.1) circle (2pt);
\fill (7.77,0.1) circle (2pt);
\draw[fill=white] (7.55,-1.25) circle (2.4pt);
\node at (6.925,0.1) {$\cdots$};
\node at (6.925,-1.25) {$\cdots$};
\draw[decorate,decoration={brace,mirror,amplitude=4pt}] (6.3,-1.56)--(7.55,-1.56) node[midway,below=6pt] {$\tau$};
% Two representatives of the a_1 type, with intervening ellipses.
\fill (8.65,0.1) circle (2pt);
\fill (9.55,0.1) circle (2pt);
\node at (9.1,0.1) {$\cdots$};
\draw[decorate,decoration={brace,amplitude=4pt}] (8.65,0.36)--(9.55,0.36) node[midway,above=6pt] {$a_{1}$};
% Two representatives of the b_1 type, with intervening ellipses.
\draw (10.35,0.1)--(10.35,-1.25);
\fill (10.35,0.1) circle (2pt);
\draw[fill=white] (10.35,-1.25) circle (2.4pt);
\draw (11.25,0.1)--(11.25,-1.25);
\fill (11.25,0.1) circle (2pt);
\draw[fill=white] (11.25,-1.25) circle (2.4pt);
\node at (10.8,0.1) {$\cdots$};
\node at (10.8,-1.25) {$\cdots$};
\draw[decorate,decoration={brace,mirror,amplitude=4pt}] (10.35,-1.56)--(11.25,-1.56) node[midway,below=6pt] {$b_{1}$};
% Representative vertices and pendant neighbors in C_{h-1}.
\fill (12.55,0.1) circle (2pt);
\fill (13.45,0.1) circle (2pt);
\node at (13,0.1) {$\cdots$};
\draw[decorate,decoration={brace,amplitude=4pt}] (12.55,0.36)--(13.45,0.36) node[midway,above=6pt] {$a_{h-1}$};
\draw (14.25,0.1)--(14.25,-1.25);
\fill (14.25,0.1) circle (2pt);
\draw[fill=white] (14.25,-1.25) circle (2.4pt);
\draw (15.15,0.1)--(15.15,-1.25);
\fill (15.15,0.1) circle (2pt);
\draw[fill=white] (15.15,-1.25) circle (2.4pt);
\node at (14.7,0.1) {$\cdots$};
\node at (14.7,-1.25) {$\cdots$};
\draw[decorate,decoration={brace,mirror,amplitude=4pt}] (14.25,-1.56)--(15.15,-1.56) node[midway,below=6pt] {$b_{h-1}$};
% Representative vertices and pendant neighbors in C_{h}.
\fill (15.85,0.1) circle (2pt);
\fill (16.75,0.1) circle (2pt);
\node at (16.3,0.1) {$\cdots$};
\draw[decorate,decoration={brace,amplitude=4pt}] (15.85,0.36)--(16.75,0.36) node[midway,above=6pt] {$a_{h}$};
\draw (17.55,0.1)--(17.55,-1.25);
\fill (17.55,0.1) circle (2pt);
\draw[fill=white] (17.55,-1.25) circle (2.4pt);
\draw (18.45,0.1)--(18.45,-1.25);
\fill (18.45,0.1) circle (2pt);
\draw[fill=white] (18.45,-1.25) circle (2.4pt);
\node at (18,0.1) {$\cdots$};
\node at (18,-1.25) {$\cdots$};
\draw[decorate,decoration={brace,mirror,amplitude=4pt}] (17.55,-1.56)--(18.45,-1.56) node[midway,below=6pt] {$b_{h}$};
\end{tikzpicture}
}
\caption{Graph $D(\boldsymbol a,\boldsymbol b,\boldsymbol\mu;\kappa,\tau)$.}
\label{fig:chain-families}
\end{figure}
% End embedded drawing: figures/chain_families.tex

Let $\boldsymbol a=(a_{0},a_{1},\ldots,a_{h})$, $\boldsymbol b=(b_{0},b_{1},\ldots,b_{h})$  and $\boldsymbol\mu=(\mu_{0},\mu_{1},\ldots,\mu_{h-1})$, where $h$ and $\mu_{i}$ ($0\leq i<h$) are  positive integers and $a_{j},b_{j}$ ($0\leq j\leq h$) are nonnegative integers,  and  $\kappa,\tau$ be nonnegative integers.
Denote by  $D(\boldsymbol a,\boldsymbol b,\boldsymbol\mu;\kappa,\tau)$ the split graph with split partition $(C,Z)$ as follows.
Let $C_{0},\ldots,C_{h}$ be nonempty pairwise disjoint cliques in $C$ such that $C=\bigcup\nolimits_{i=0}^{h} C_{i}$, $|C_{0}|=a_{0}+b_{0}+\kappa+2\tau$ and $|C_{j}|=a_{j}+b_{j}$ for $1\leq j\leq h$.
Let $X_{0},\ldots,X_{h-1},Y$ be pairwise disjoint   sets of $Z$  such that $Z=\bigcup\nolimits_{i=0}^{h-1} X_{i}\cup Y$ and $|X_{i}|=\mu_{i}$.
For every $x\in X_{i}$, $N(x)=\bigcup\nolimits_{j=0}^{i} C_{j}\triangleq Q_{i}$.
For $0\leq j\leq h$, $C_{j}$ has $a_{j}$ vertices with no neighbors in $Y$ and $b_{j}$ vertices each with exactly one pendant neighbor in $Y$.
$C_{0}$ has $\kappa$ vertices each with exactly two pendant neighbors in $Y$,
and the remaining $2\tau$ vertices form $\tau$ disjoint pairs, each with exactly one common neighbor of degree $2$ in $Y$ (see Figure~\ref{fig:chain-families}).
It is easy to see $|Y|=\sum_{i=0}^{h}b_{i}+2\kappa+\tau$.

Let $n$ be a sufficiently large positive integer.
For each $1\leq i\leq355$, let $A_{i}=D(\boldsymbol a,\boldsymbol b,\boldsymbol\mu;0,0)$, where $\boldsymbol a,\boldsymbol b,\boldsymbol\mu$ are given in Table~\ref{tab:s5-explicit-1}.
For each $1\leq j\leq94$, let $\mathcal F_{j}(n)=D(\boldsymbol a,\boldsymbol b,\boldsymbol\mu;\kappa,\tau)$, where $\boldsymbol a,\boldsymbol b,\boldsymbol\mu,\kappa,\tau$ are given in Table~\ref{tab:s5-explicit-F}.
Let $\mathcal H=\{A_{i}:1\leq i\leq355\}\cup\{\mathcal F_{j}(n):1\leq j\leq94,\ n\in\mathbb Z_{>0}\}$.
Our main result is the following theorem.

\begin{theorem}\label{thm:main}
Let $G$ be a graph with $|V(G)|\geq3$.
Then $G$ is a connected distance-hereditary split graph with $\lambda_{2}(G)\leq1$ if and only if $G$ is a connected induced subgraph of some graph $H\in\mathcal H$.
\end{theorem}

The remainder of this paper is organized as follows.
In Section~2, we give some preliminary results on distance-hereditary split graphs and matrix theory.
In Section~3, we establish structural restrictions and spectral conditions for the graph families introduced above.
In Section~4, we prove Theorem~\ref{thm:main} and  determine the multiplicity of the eigenvalue $1$.

\section{Preliminaries}

For graphs $H_{1}$ and $H_{2}$, we write $H_{1}\subset H_{2}$ if $H_{1}$ is isomorphic to an induced subgraph of $H_{2}$.
%For a positive integer $a$, we denote by $aG$ the union of $a$ copies of a graph $G$.
Let $\mi(G)$ denote the multiplicity of the eigenvalue $1$ of $G$.
Let $\diag(M_{1},\ldots,M_{k})$ denote the block diagonal matrix with diagonal blocks $M_{1},\ldots,M_{k}$ and all other blocks zero.
For a matrix $M$, let $(M)_{uv}$ denote the entry in row $u$ and column $v$.
Let $I_{n}$ and $J_{n}$ denote the identity and all-ones matrices of order $n$, respectively.
We write $I,J$ when their orders are clear.
Let $\one$ denote the all-ones  column vector.
Let $\mathbf{0}$ denote a zero column vector of appropriate order and $0_{s\times t}$ the $s\times t$ zero matrix.
For a real symmetric matrix $M$, let $\npos(M)$ and $n_{0}(M)$ denote the numbers of its positive and zero eigenvalues, respectively.
%For $x\in\mathbb{R}$, let $\lfloor x\rfloor=\max\{n\in\mathbb{Z}:n\leq x\}$.
Let $\rho(M)$ denote the spectral radius of $M$.

\begin{lemma}\label{lem:equ}
Let $G=(C,Z)$ be a connected split graph and $B$ be the submatrix of $A(G)$ whose rows correspond to the vertices of $C$ and columns to the vertices of $Z$. Then\\

 $\lambda_{2}(G)\leq1$ if and only if $\npos(J_{|C|}-2I_{|C|}+BB^{\mathsf T})\leq1$; $\mi(G)=\nzero(J_{|C|}-2I_{|C|}+BB^{\mathsf T})$.
\end{lemma}
\begin{proof}
By a suitable ordering to the vertices of $G$, we have $$A(G)-I=\begin{psmallmatrix}J_{|C|}-2I_{|C|}&B\\B^{\mathsf T}&-I_{|Z|}\end{psmallmatrix}.$$
Let $D=\begin{psmallmatrix}I_{|C|}&B\\0_{|Z|\times|C|}&I_{|Z|}\end{psmallmatrix}$. 
Then  $$D(A(G)-I)D^{\mathsf T}=\diag\bigl(J_{|C|}-2I_{|C|}+BB^{\mathsf T},-I_{|Z|}\bigr).$$
This completes the proof.
\end{proof}

\begin{lemma}[{\cite[Corollary~2.5.2]{BrouwerHaemers2012}}]\label{lem:interlacing}
Let $A$ be a real symmetric matrix of order $n$ and $B$ be a principal submatrix of $A$ of order $k$.
Let $\lambda_{1}(A)\geq\cdots\geq\lambda_{n}(A)$ be the eigenvalues of $A$ and $\lambda_{1}(B)\geq\cdots\geq\lambda_{k}(B)$ be the eigenvalues of $B$.
Then
$$
\lambda_{i}(A)\geq\lambda_{i}(B)\geq\lambda_{n-k+i}(A),\qquad i=1,2,\ldots,k.
$$
\end{lemma}

By Lemma~\ref{lem:interlacing}, we have the following corollary.

\begin{corollary}\label{cor:interlacing}
Let $R$ be an induced subgraph of $G$.
If $\lambda_{2}(R)>1$, then $\lambda_{2}(G)>1$.
\end{corollary}

Next we give a property of distance-hereditary split graphs.

\begin{lemma}[\cite{Cameron2016,Pardal2019}]\label{lem:laminar}
A split graph with split partition $(C,Z)$ is gem-free if and only if $N(x)\cap N(y)=\varnothing$  or $N(x)\subseteq N(y)$ or $N(y)\subseteq N(x)$ for all $x,y\in Z$.
\end{lemma}
%\begin{proof}
%Suppose that $G$ is gem-free.
%Then $N(x)\cap N(y)=\varnothing$, $N(x)\subseteq N(y)$ or $N(y)\subseteq N(x)$ for all $x,y\in Z$.
%Otherwise, let $a\in N(x)\cap N(y)$, $b\in N(x)\setminus N(y)$ and $c\in N(y)\setminus N(x)$.
%Then $G[\{x,y,a,b,c\}]\cong\gem$, a contradiction.

%Suppose that the stated neighborhood condition holds and $G$ contains an induced gem.
%In every split partition of a gem, the two ends of its $P_4$ belong to the independent set.
%The two internal vertices and the universal vertex belong to the clique.
%Then the neighborhoods of the two ends intersect and neither contains the other, a contradiction.
%\end{proof}

\begin{lemma}[{Courant--Fischer, \cite[Theorem~2.4.1]{BrouwerHaemers2012}}]\label{lem:courant-fischer}

Let $A$ be a real symmetric matrix of order $n$ with eigenvalues $\lambda_{1}(A)\geq\cdots\geq\lambda_{n}(A)$.
Let $W$ be a $k$-dimensional subspace of $\mathbb R^n$.
Let $W^\perp=\{y\in\mathbb R^n:x^{\mathsf T}y=0\text{ for every }x\in W\}$ denote the orthogonal complement of $W$.
Then
$$
\lambda_{k}(A)\geq\min_{\substack{x\in W\\x\neq\mathbf{0}}}\frac{x^{\mathsf T}Ax}{x^{\mathsf T}x}\quad\textup{and}\quad
\lambda_{k+1}(A)\leq\max_{\substack{x\in W^\perp\\x\neq\mathbf{0}}}\frac{x^{\mathsf T}Ax}{x^{\mathsf T}x}\quad(1\leq k<n).
$$

\end{lemma}

For $x\in Z$, let $m_{x}$ be the number of vertices $y\in Z$ with $N(y)=N(x)$.
For $x\in Z$ with $N(x)\neq\varnothing$ and $S\subseteq N(x)$, let $\one_{S}$ be the column vector indexed by $N(x)$ with entries $1$ on $S$ and $0$ on $N(x)\setminus S$.
Let $L_{x}=\sum_{\substack{y\in Z\\N(y)\subseteq N(x)}}\one_{N(y)}\one_{N(y)}^{\mathsf T}$.
For $x\in Z$, the neighborhood $N(x)$ is said to be heavy, critical or subcritical if $\rho(L_{x})>2$, $\rho(L_{x})=2$ or $\rho(L_{x})<2$, respectively.

\begin{samepage}
\begin{lemma}\label{lem:light}

Let $G=(C,Z)$ be a connected split graph and $x\in Z$.
If $N(x)$ is critical or subcritical, then one of the following holds:
\begin{enumerate}[label=\textup{(\roman*)}]
\item $|N(x)|=1$ and $m_{x}\in \{1,2\}$;
\item $|N(x)|=2$, $m_{x}=1$ and there is no $y\in Z$ with $N(y)\subset N(x)$.
\end{enumerate}

\end{lemma}
\end{samepage}

\begin{proof}

Let $\one$ be the all-ones vector indexed by $N(x)$.
Since every entry of $L_{x}$ is at least $m_{x}$, $m_{x}|N(x)|^2\leq\one^{\mathsf T}L_{x}\one$.
Then $m_{x}|N(x)|\leq\frac{\one^{\mathsf T}L_{x}\one}{\one^{\mathsf T}\one}\leq\rho(L_{x})\leq2$ by Lemma~\ref{lem:courant-fischer}.
Thus, either $|N(x)|=1$ and $m_{x}\in \{1,2\}$, or $|N(x)|=2$ and $m_{x}=1$.

Suppose that $|N(x)|=2$ and $m_{x}=1$.
If $y\in Z$ and $N(y)\subset N(x)$, then $|N(y)|=1$.
Then $\rho(L_{x})\geq\frac{\one^{\mathsf T}L_{x}\one}{\one^{\mathsf T}\one}\geq   \frac{\one^{\mathsf T}(J_{2}+m_{y}\one_{N(y)}\one_{N(y)}^{\mathsf T})\one}{\one^{\mathsf T}\one}\geq\frac{4+m_{y}}{2}>2$ by Lemma~\ref{lem:courant-fischer}, a contradiction.
This completes the proof.
\end{proof}

Let  $M$ be a real symmetric matrix and $\ker M=\{x:Mx=\mathbf{0}\}$ denote the kernel of $M$.
For a subspace $W$ and a vector column  $e$, we write $W\perp e$ if $x^{\mathsf T}e=0$ for every $x\in W$.
We write $W\not\perp e$ if $x^{\mathsf T}e\neq0$ for some $x\in W$.
For a symmetric matrix $M$ and a vector $e$ with $\ker M\perp e$, let $\eta(M,e)=e^{\mathsf T}z$, where $Mz=e$.
If $e=\one$, then we write $\eta(M)$.
Since $M$ is real symmetric and $\ker M\perp e$, the equation $Mz=e$ has a solution.
If $Mz_{1}=Mz_{2}=e$, then $z_{1}-z_{2}\in\ker M$, and hence $e^{\mathsf T}z_{1}=e^{\mathsf T}z_{2}$.
Thus, the value of $\eta(M,e)$ is uniquely determined.

\begin{samepage}
\begin{lemma}\label{lem:update}
Suppose that $M$ is a real symmetric matrix, $\npos(M)=1$, $\ker M\perp e$ and $\omega>0$.
Let $\Delta=1+\omega \eta(M,e)$ and $\widehat M=M+\omega ee^{\mathsf T}$.
Then the following hold:
\begin{enumerate}[label=\textup{(\roman*)}]
\item If $\Delta>0$, then $\npos(\widehat M)=1$, $\nzero(\widehat M)=\nzero(M)$ and $\eta(\widehat M,e)=\eta(M,e)/\Delta$.
\item If $\Delta=0$, then $\npos(\widehat M)=1$, $\nzero(\widehat M)=\nzero(M)+1$ and $\ker\widehat M\not\perp e$.
\item If $\Delta<0$, then $\npos(\widehat M)=2$.
\end{enumerate}
\end{lemma}
\end{samepage}

\begin{proof}
Let $q$ satisfy $Mq=e$ and $e^{\mathsf T}q=\eta(M,e)$.
Let $K=\begin{pmatrix}M&e\\e^{\mathsf T}&-1/\omega\end{pmatrix}$, $E_{1}=\begin{pmatrix}I&\mathbf{0}\\-q^{\mathsf T}&1\end{pmatrix}$ and $E_{2}=\begin{pmatrix}I&\omega e\\\mathbf{0}^{\mathsf T}&1\end{pmatrix}$.
Then $E_{1}KE_{1}^{\mathsf T}=\diag(M,-\Delta/\omega)=H_{1}$ and $E_{2}KE_{2}^{\mathsf T}=\diag(\widehat M,-1/\omega)=H_{2}$.
Thus, $\npos(H_{1})=\npos(H_{2})$ and $\nzero(H_{1})=\nzero(H_{2})$.
Then $\widehat Mq=Mq+\omega e(e^{\mathsf T}q)=\Delta e$.

Suppose that $\Delta>0$. 
Then $-\Delta/\omega<0$.
%Thus, $H_1$ has one more negative eigenvalue than $M$ and $H_2$ has one more negative eigenvalue than $\widehat M$.
Thus, $\npos(\widehat M)=\npos(M)=1$ and $\nzero(\widehat M)=\nzero(M)$.
Then $\widehat Mx=Mx+\omega e(e^{\mathsf T}x)=\mathbf{0}$ for every $x\in\ker M$.
Thus, $\ker M\subseteq\ker\widehat M$.
Since their dimensions are equal, $\ker\widehat M=\ker M$ and $\ker\widehat M\perp e$.
Since $\widehat M(q/\Delta)=e$, $\eta(\widehat M,e)=e^{\mathsf T}q/\Delta=\eta(M,e)/\Delta$.
Thus, (i) holds.

Suppose that $\Delta=0$. 
Then $\diag(M,-\Delta/\omega)=\diag(M,0)$.
%This matrix has one positive eigenvalue and $\nzero(M)+1$ zero eigenvalues.
Since $-1/\omega<0$, $\npos(\widehat M)=1$ and $\nzero(\widehat M)=\nzero(M)+1$.
Since $\widehat Mq=\Delta e=\mathbf{0}$, $q\in\ker\widehat M$.
Since $\Delta=1+\omega\eta(M,e)=0$, $e^{\mathsf T}q=\eta(M,e)=-1/\omega\neq0$.
Then $\ker\widehat M\not\perp e$.
Thus, (ii) holds.

Suppose that $\Delta<0$.
Then $-\Delta/\omega>0$.
Since $\npos(M)=1$, $\diag(M,-\Delta/\omega)$ has two positive eigenvalues.
Then $\diag(\widehat M,-1/\omega)$ also has two positive eigenvalues.
Since $-1/\omega<0$, $\npos(\widehat M)=2$.
Thus, (iii) holds.
\end{proof}

\begin{lemma}[{see~\cite[Theorem~2.5.1(i)]{BrouwerHaemers2012}}]\label{lem:positive-eigenvalues}
Let $A$ be a real symmetric matrix of order $n$ and $Q$ be an $n\times k$ matrix with $Q^{\mathsf T}Q=I_{k}$.
If $Q^{\mathsf T}AQ$ is positive definite, then $A$ has at least $k$ positive eigenvalues.
\end{lemma}

\begin{lemma}\label{lem:mainzero}
{
Let $M$ be a real symmetric matrix of order $n$ and $\npos(M)\geq 1$.
Let $M_{1}$ be a real symmetric negative definite matrix of order $r$, $e\in\mathbb R^n$ and $f\in\mathbb R^r$.
Let $M_{2}=\diag(M,M_{1})$ and $e_{1}=(e^{\mathsf T},f^{\mathsf T})^{\mathsf T}$.
If there is a vector $y\in\ker M$ with $e^{\mathsf T}y\neq0$, then $M+\omega ee^{\mathsf T}$ and $M_{2}+\omega e_{1}e_{1}^{\mathsf T}$ have at least two positive eigenvalues for every $\omega>0$.
}
\end{lemma}

\begin{proof}
{
Let $y_{0}=y/\sqrt{y^{\mathsf T}y}$.
Then $My_{0}=\mathbf{0}$ and $e^{\mathsf T}y_{0}\neq0$.
Let $x$ be a unit eigenvector of $M$ with eigenvalue $\lambda>0$.
Then $x^{\mathsf T}y_{0}=0$.
Let $Q=(x,y_{0})$, $a=e^{\mathsf T}x$ and $b=e^{\mathsf T}y_{0}$.
Then $Q^{\mathsf T}Q=I_{2}$.
Let
$$
H=Q^{\mathsf T}(M+\omega ee^{\mathsf T})Q
=\begin{pmatrix}
\lambda+\omega a^2&\omega ab\\
\omega ab&\omega b^2
\end{pmatrix}.
$$
It is easy to see that both eigenvalues of $H$ are positive.
Then $M+\omega ee^{\mathsf T}$ has at least two positive eigenvalues by Lemma~\ref{lem:positive-eigenvalues}.

Let $\mathbf{0}_{r}$ be the zero column vector of order $r$.
Let $x_{1}=(x^{\mathsf T},\mathbf{0}_{r}^{\mathsf T})^{\mathsf T}$ and $y_{1}=(y_{0}^{\mathsf T},\mathbf{0}_{r}^{\mathsf T})^{\mathsf T}$.
Then $M_{2}x_{1}=\lambda x_{1}$, $M_{2}y_{1}=\mathbf{0}$, $e_{1}^{\mathsf T}x_{1}=a$ and $e_{1}^{\mathsf T}y_{1}=b$.
Let $P=(x_{1},y_{1})$.
Then $P^{\mathsf T}P=I_{2}$ and $P^{\mathsf T}(M_{2}+\omega e_{1}e_{1}^{\mathsf T})P=H$.
Thus, $M_{2}+\omega e_{1}e_{1}^{\mathsf T}$ has at least two positive eigenvalues by Lemma~\ref{lem:positive-eigenvalues}.
This completes the proof.
}
\end{proof}

\section{Structure and spectral conditions}

% Begin embedded drawing: figures/D1_family.tex
\begin{figure}[htbp]
\centering
\resizebox{12.767cm}{!}{%
\begin{tikzpicture}[x=1cm,y=1cm,line width=0.5pt,every node/.style={font=\small}]
% All filled vertices belong to the same clique C.
\draw[rounded corners=4pt,gray] (-0.25,1.40) rectangle (13.95,1.90);
\node[anchor=east] at (-0.42,1.65) {$Z$};
\draw[rounded corners=5pt,gray] (-0.25,-0.50) rectangle (13.95,0.50);
\node[anchor=east] at (-0.42,0) {$C$};

% Four complete induced subgraphs within C.
\draw[rounded corners=3pt,gray] (0,-0.35) rectangle (2.20,0.30);
\draw[rounded corners=3pt,gray] (2.85,-0.35) rectangle (5.15,0.30);
\draw[rounded corners=3pt,gray] (5.80,-0.35) rectangle (8.90,0.30);
\draw[rounded corners=3pt,gray] (9.55,-0.35) rectangle (13.70,0.30);
\foreach \pos/\j in {1.10/0,4.00/1,7.35/2,11.625/3}{
  \node at (\pos,-0.18) {$C_{\j}$};
}

% C_0: two representative vertices and an ellipsis.
\foreach \pos in {0.45,1.75}{\fill (\pos,0.12) circle (2pt);}
\node at (1.10,0.12) {$\cdots$};

% C_1: two pendant-neighbor patterns, repeated beta times.
\foreach \pos in {3.35,4.65}{
  \draw (\pos,0.12)--(\pos,1.65);
  \fill (\pos,0.12) circle (2pt);
  \draw[fill=white] (\pos,1.65) circle (2.4pt);
}
\node at (4.00,0.12) {$\cdots$};
\node at (4.00,1.65) {$\cdots$};

% C_2: gamma clique vertices, each with two pendant neighbors.
\foreach \pos in {6.40,8.30}{
  \draw (\pos,0.12)--({\pos-0.35},1.65)
        (\pos,0.12)--({\pos+0.35},1.65);
  \fill (\pos,0.12) circle (2pt);
  \draw[fill=white] ({\pos-0.35},1.65) circle (2.4pt);
  \draw[fill=white] ({\pos+0.35},1.65) circle (2.4pt);
}
\node at (7.35,0.12) {$\cdots$};
\node at (7.35,1.65) {$\cdots$};

% C_3: delta disjoint pairs, each with one degree-two common neighbor.
\foreach \pos in {10.40,12.85}{
  \draw ({\pos-0.35},0.12)--(\pos,1.65)--({\pos+0.35},0.12);
  \fill ({\pos-0.35},0.12) circle (2pt);
  \fill ({\pos+0.35},0.12) circle (2pt);
  \draw[fill=white] (\pos,1.65) circle (2.4pt);
}
\node at (11.625,0.12) {$\cdots$};
\node at (11.625,1.65) {$\cdots$};

% Braces count C_0/C_2 vertices, and the one-to-one Z patterns for C_1/C_3.
\draw[decorate,decoration={brace,mirror,amplitude=4pt}]
  (0.45,-0.72)--(1.75,-0.72) node[midway,below=6pt] {$\alpha$};
\draw[decorate,decoration={brace,amplitude=4pt}]
  (3.35,2.13)--(4.65,2.13) node[midway,above=6pt] {$\beta$};
\draw[decorate,decoration={brace,mirror,amplitude=4pt}]
  (6.40,-0.72)--(8.30,-0.72) node[midway,below=6pt] {$\gamma$};
\draw[decorate,decoration={brace,amplitude=4pt}]
  (10.40,2.13)--(12.85,2.13) node[midway,above=6pt] {$\delta$};
\end{tikzpicture}
}
\caption{Graph $D(\alpha,\beta,\gamma,\delta)$.}
\label{fig:D1}
\end{figure}
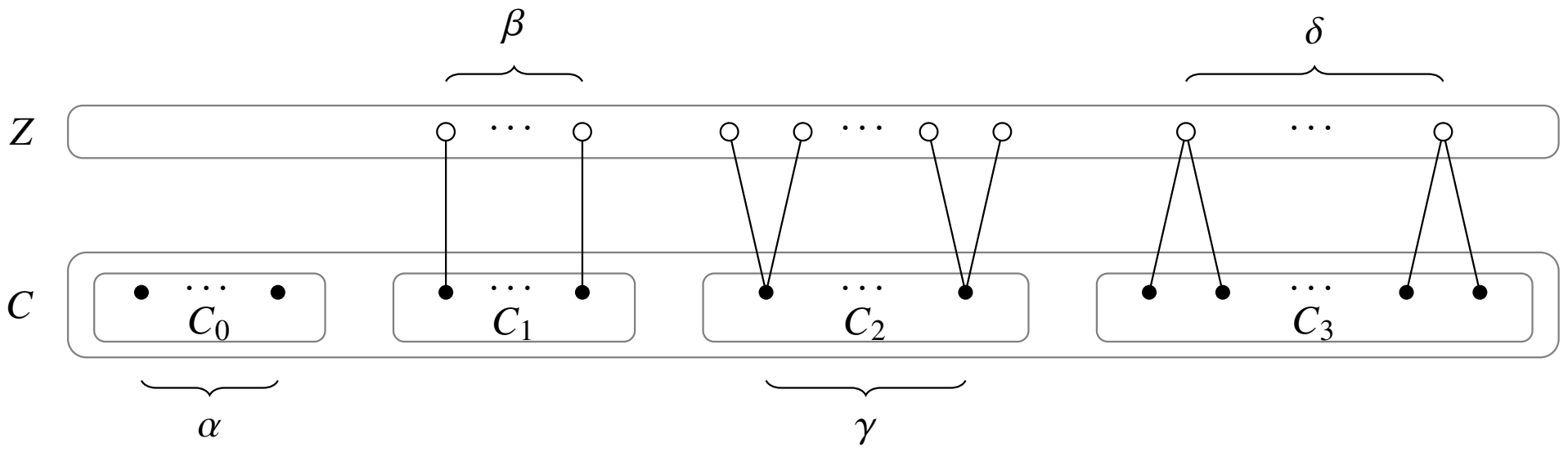
% End embedded drawing: figures/D1_family.tex

Let $\alpha,\beta,\gamma,\delta$ be nonnegative integers.
Denote by $D(\alpha,\beta,\gamma,\delta)$ the split graph with split partition $(C,Z)$ as follows.
Let $C_{0},C_{1},C_{2},C_{3}$ be pairwise disjoint cliques in $C$ such that $C=\bigcup\nolimits_{i=0}^{3} C_{i}$, $|C_{0}|=\alpha$, $|C_{1}|=\beta$, $|C_{2}|=\gamma$, $|C_{3}|=2\delta$ and $|Z|=\beta+2\gamma+\delta$.
The vertices of $C_{0}$ have no neighbors in $Z$.
Each vertex of $C_{1}$ has exactly one pendant neighbor in $Z$, and each vertex of $C_{2}$ has exactly two pendant neighbors in $Z$.
The vertices of $C_{3}$ form $\delta$ disjoint pairs and each pair has exactly one common neighbor of degree $2$ in $Z$ (see Figure~\ref{fig:D1}).
Then $\alpha+\beta+\gamma+2\delta\geq 2$ and if $\delta>0$, then $\alpha+\beta+\gamma+2\delta\geq 3$, for $C$ is a
maximum clique. 
In particular, $D(2,0,0,0)=K_{2}$.
Let $\mathcal{D}_{1}$ be the set of all $D(\alpha,\beta,\gamma,\delta)$.

In $D(\boldsymbol a,\boldsymbol b,\boldsymbol\mu;\kappa,\tau)$, let $\sigma_{j}=a_{j}+2b_{j}$ for $0\leq j\leq h$ and $\mu_{h}=1$.
Let $(c_{0},p_{0})=(\sigma_{0},\mu_{0}\sigma_{0}-2)$ if $\kappa+\tau=0$ and $(c_{0},p_{0})=(1,\mu_{0})$ if $\kappa+\tau\geq 1$.
For $1\leq j\leq h$, let
\refstepcounter{equation}\label{eq:recurrence}$$
c_{j}=2c_{j-1}-\sigma_{j}p_{j-1},\qquad
p_{j}=2p_{j-1}+\mu_{j}c_{j}.
\eqno(\theequation)$$
%For the families $\mathcal{D}_{2}$ and $\mathcal{D}_{3}$ defined below, the parameters satisfy
%\refstepcounter{equation}\label{eq:admissible}$$
%p_j>0\quad(1\leq j<h),\qquad p_h\geq0.
%\eqno(\theequation)$$

Let $\mathcal{D}_{2}=\{D(\boldsymbol a,\boldsymbol b,\boldsymbol\mu;0,0):
1\leq h\leq4,\ \mu_{0}(a_{0}+2b_{0})>2,\ p_{j}>0\ (1\leq j<h),\ p_{h}\geq0\}$
and $\mathcal{D}_{3}=\{D(\boldsymbol a,\boldsymbol b,\boldsymbol\mu;\kappa,\tau):
\kappa+\tau\geq1,\ 1\leq h\leq3,\ p_{j}>0\ (1\leq j<h),\ p_{h}\geq0\}$.

\begin{lemma}\label{lem:base}
Let $G\cong D(\boldsymbol a,\boldsymbol b,\boldsymbol\mu;\kappa,\tau)$, $x_{0}\in X_{0}$ and $M_{0}=L_{x_{0}}-2I$.
Then the following hold:
\begin{enumerate}[label=\textup{(\arabic*)}]
\item If $\kappa+\tau=0$, then $\npos(M_0)=1$ if and only if $\mu_{0}(a_{0}+2b_{0})>2$;
if $\kappa+\tau=0$ and $\mu_{0}(a_{0}+2b_{0})>2$, then $\nzero(M_{0})=0$ and $\eta(M_{0})=(a_{0}+2b_{0})/(\mu_{0}(a_{0}+2b_{0})-2)$.

\item If $\kappa+\tau>0$, then $\npos(M_0)=1$, $\nzero(M_{0})=\kappa+\tau-1$, $\ker M_{0}\perp\one$ and $\eta(M_{0})=1/\mu_{0}$.
\end{enumerate}
\end{lemma}

\begin{proof}
For $u,v\in C_{0}$, not necessarily distinct, $(L_{x_{0}})_{uv}$ equals the number of vertices $x\in Z$ adjacent to both $u$ and $v$ and satisfying $N(x)\subseteq C_{0}$.
For any vertex $y\in X_{0}$, $N(y)=C_{0}$.
Thus, $(L_{x_{0}})_{uv}$ is at least $\mu_{0}$.
Every vertex $x\in Z\setminus X_{0}$ with $N(x)\subseteq C_{0}$ belongs to $Y$.

Suppose that $\kappa+\tau=0$.
Then $|C_{0}|=a_{0}+b_{0}$ and $(L_{x_{0}})_{uv}=\mu_{0}$ for distinct $u,v\in C_{0}$.
For $u\in C_{0}$, $(L_{x_{0}})_{uu}=\mu_{0}+i$.
Here $i=0$ if $N(u)\cap Y=\varnothing$ and $i=1$ if $N(u)\cap Y\neq\varnothing$.
Let $F_{1}=\diag(-2I_{a_{0}},-I_{b_{0}})$.
Then $F_{1}$ is invertible and $M_{0}=F_{1}+\mu_{0}J$.
Let $q=F_{1}^{-1}\one$.
Then $q=(-1/2,\ldots,-1/2,-1,\ldots,-1)^{\mathsf T}$, where $-1/2$ occurs $a_{0}$ times and $-1$ occurs $b_{0}$ times.
Then $\one^{\mathsf T}q=-a_{0}/2-b_{0}=-\sigma_{0}/2$.
Let $K=\begin{pmatrix}F_{1}&\one\\\one^{\mathsf T}&-1/\mu_{0}\end{pmatrix}$, $E_{1}=\begin{pmatrix}I&\mathbf{0}\\-q^{\mathsf T}&1\end{pmatrix}$ and $E_{2}=\begin{pmatrix}I&\mu_{0}\one\\\mathbf{0}^{\mathsf T}&1\end{pmatrix}$.
Then $E_{1}KE_{1}^{\mathsf T}=\diag(\widehat F_{2},\sigma_{0}/2-1/\mu_{0})$ and $E_{2}KE_{2}^{\mathsf T}=\diag(M_{0},-1/\mu_{0})$.
Then $\nzero(E_{1}KE_{1}^{\mathsf T})=\nzero(E_{2}KE_{2}^{\mathsf T})$ and $\npos(E_{1}KE_{1}^{\mathsf T})=\npos(E_{2}KE_{2}^{\mathsf T})$ for $E_{1},E_{2}$ are nonsingular.
Since all eigenvalues of $F_{1}$ are negative and $-1/\mu_{0}<0$, $M_{0}$ has exactly one positive eigenvalue if and only if $\sigma_{0}/2-1/\mu_{0}>0$, which means that $\mu_{0}\sigma_{0}>2$.
If $\mu_{0}\sigma_{0}>2$, then $\nzero(M_{0})=0$.
Since $F_{1}q=\one$ and $\one^{\mathsf T}q=-\sigma_{0}/2$, $M_{0}q=(1-\mu_{0}\sigma_{0}/2)\one$.
If $\mu_{0}\sigma_{0}>2$, then $M_{0}z=\one$ for $z=q/(1-\mu_{0}\sigma_{0}/2)$.
Then $\eta(M_{0})=\one^{\mathsf T}z=(-\sigma_{0}/2)/(1-\mu_{0}\sigma_{0}/2)=\sigma_{0}/(\mu_{0}\sigma_{0}-2)$.

Suppose that $\kappa+\tau>0$.
If $u\in C_{0}$ has exactly $r$ pendant neighbors in $Y$ and no other neighbors in $Y$, where $r\in\{0,1,2\}$, then $(L_{x_{0}})_{uu}=\mu_{0}+r$ and $(M_{0}-\mu_{0}J)_{uu}=r-2$.
For a pair $\{u,v\}$ with a common neighbor $y\in Y$ satisfying $N(y)=\{u,v\}$, $(L_{x_{0}})_{uu}=(L_{x_{0}})_{uv}=(L_{x_{0}})_{vv}=(L_{x_{0}})_{vu}=\mu_{0}+1$.
Then $(M_{0}-\mu_{0}J)_{uu}=(M_{0}-\mu_{0}J)_{vv}=-1$ and $(M_{0}-\mu_{0}J)_{uv}=(M_{0}-\mu_{0}J)_{vu}=1$.
Let $R=J_{2}-2I_{2}$.
Then $M_{0}=F_{2}+\mu_{0}J$, where $F_{2}=\diag(-2I_{a_{0}},\allowbreak -I_{b_{0}},\allowbreak 0_{\kappa\times\kappa},\allowbreak R,\allowbreak\ldots,\allowbreak R)$ and $R$ occurs $\tau$ times.
Here $0_{\kappa\times\kappa}$ is the zero matrix of order $\kappa$.
The eigenvalues of $R$ are $0$ and $-2$.
Thus, all eigenvalues of $F_{2}$ are nonpositive and $\dim\ker F_{2}=\kappa+\tau$.
Let $v=(v_{1},\ldots,v_{|C_{0}|})^{\mathsf T}$ be a column vector, where
$$
v_{i}=
\begin{cases}
1,&\kappa>0\text{ and }i=a_{0}+b_{0}+1;\\
1,&\kappa=0\text{ and }i\in\{a_{0}+b_{0}+1,a_{0}+b_{0}+2\};\\
0,&\text{otherwise}.
\end{cases}$$
Let $s=\one^{\mathsf T}v$.
Then $F_{2}v=\mathbf{0}$ and $s>0$.
Since $v^{\mathsf T}M_{0}v=\mu_{0}s^2>0$, $M_{0}$ has a positive eigenvalue.
For every vector $x$ with $\one^{\mathsf T}x=0$, $x^{\mathsf T}M_{0}x=x^{\mathsf T}F_{2}x\leq0$.
Suppose that $M_{0}$ has two positive eigenvalues.
Let $u,w$ be orthonormal eigenvectors corresponding to positive eigenvalues of $M_{0}$.
There are real numbers $a,b$, not both zero, such that $x=au+bw$ satisfies $\one^{\mathsf T}x=0$.
Then $x^{\mathsf T}M_{0}x>0$, a contradiction.
Thus, $\npos(M_{0})=1$.
If $M_{0}x=\mathbf{0}$, then $v^{\mathsf T}M_{0}x=\mu_{0}s\one^{\mathsf T}x=0$.
Since $\mu_{0}s>0$, $\one^{\mathsf T}x=0$.
Then $F_{2}x=M_{0}x-\mu_{0}\one(\one^{\mathsf T}x)=\mathbf{0}$.
Conversely, if $F_{2}x=\mathbf{0}$ and $\one^{\mathsf T}x=0$, then $M_{0}x=\mathbf{0}$.
\mbox{Thus, $\ker M_{0}=\{x:F_{2}x=\mathbf{0},\ \one^{\mathsf T}x=0\}$}.
Then $\nzero(M_{0})=\kappa+\tau-1$ and $\ker M_{0}\perp\one$.
Let $z=v/(\mu_{0}s)$.
Then $F_{2}z=\mathbf{0}$ and $\one^{\mathsf T}z=1/\mu_{0}$.
Thus, $M_{0}z=F_{2}z+\mu_{0}\one(\one^{\mathsf T}z)=\one$.
Then $\eta(M_{0})=\one^{\mathsf T}z=1/\mu_{0}$.
This completes the proof.
\end{proof}

% Begin embedded drawing: figures/R_graphs.tex
\begin{figure}[h]
\centering
\begin{subfigure}[b]{0.23\textwidth}
\centering
% Begin embedded drawing: 禁用子图/R1.tikz.tex
\begin{tikzpicture}[x=1cm,y=1cm,scale=1,every node/.style={font=\scriptsize},line width=0.45pt]
\draw[rounded corners=3pt,gray] (-1.43,-0.26) rectangle (1.43,0.22);
\draw[rounded corners=3pt,gray] (-1.43,0.88) rectangle (1.43,1.36);
\node[anchor=east,inner sep=0pt] at (-1.53,-0.02) {$C$};
\node[anchor=east,inner sep=0pt] at (-1.53,1.12) {$Z$};
\coordinate (c0) at (-0.30,0);
\coordinate (c1) at (0.30,0);
\coordinate (s0) at (-1.20,1.12);
\coordinate (s1) at (-0.60,1.12);
\coordinate (s2) at (0.00,1.12);
\coordinate (s3) at (0.60,1.12);
\coordinate (s4) at (1.20,1.12);
\draw (s0) -- (c0);
\draw (s1) -- (c0);
\draw (s2) -- (c0);
\draw (s3) -- (c1);
\draw (s4) -- (c1);
\fill (c0) circle (2pt) node[below=4pt] {$1$};
\fill (c1) circle (2pt) node[below=4pt] {$2$};
\draw[fill=white] (s0) circle (2.3pt);
\draw[fill=white] (s1) circle (2.3pt);
\draw[fill=white] (s2) circle (2.3pt);
\draw[fill=white] (s3) circle (2.3pt);
\draw[fill=white] (s4) circle (2.3pt);
\end{tikzpicture}
% End embedded drawing: 禁用子图/R1.tikz.tex
\caption*{$R_{1}$}
\end{subfigure}
\hfill
\begin{subfigure}[b]{0.23\textwidth}
\centering
% Begin embedded drawing: 禁用子图/R2.tikz.tex
\begin{tikzpicture}[x=1cm,y=1cm,scale=1,every node/.style={font=\scriptsize},line width=0.45pt]
\draw[rounded corners=3pt,gray] (-1.13,-0.26) rectangle (1.13,0.22);
\draw[rounded corners=3pt,gray] (-1.13,0.88) rectangle (1.13,1.36);
\node[anchor=east,inner sep=0pt] at (-1.23,-0.02) {$C$};
\node[anchor=east,inner sep=0pt] at (-1.23,1.12) {$Z$};
\coordinate (c0) at (-0.60,0);
\coordinate (c1) at (0.00,0);
\coordinate (c2) at (0.60,0);
\coordinate (s0) at (-0.90,1.12);
\coordinate (s1) at (-0.30,1.12);
\coordinate (s2) at (0.30,1.12);
\coordinate (s3) at (0.90,1.12);
\draw (s0) -- (c0);
\draw (s1) -- (c0);
\draw (s2) -- (c0);
\draw (s3) -- (c1);
\draw (s3) -- (c2);
\fill (c0) circle (2pt) node[below=4pt] {$1$};
\fill (c1) circle (2pt) node[below=4pt] {$2$};
\fill (c2) circle (2pt) node[below=4pt] {$3$};
\draw[fill=white] (s0) circle (2.3pt);
\draw[fill=white] (s1) circle (2.3pt);
\draw[fill=white] (s2) circle (2.3pt);
\draw[fill=white] (s3) circle (2.3pt);
\end{tikzpicture}
% End embedded drawing: 禁用子图/R2.tikz.tex
\caption*{$R_{2}$}
\end{subfigure}
\hfill
\begin{subfigure}[b]{0.23\textwidth}
\centering
% Begin embedded drawing: 禁用子图/R3.tikz.tex
\begin{tikzpicture}[x=1cm,y=1cm,scale=1,every node/.style={font=\scriptsize},line width=0.45pt]
\draw[rounded corners=3pt,gray] (-1.13,-0.26) rectangle (1.13,0.22);
\draw[rounded corners=3pt,gray] (-1.13,0.88) rectangle (1.13,1.36);
\node[anchor=east,inner sep=0pt] at (-1.23,-0.02) {$C$};
\node[anchor=east,inner sep=0pt] at (-1.23,1.12) {$Z$};
\coordinate (c0) at (-0.90,0);
\coordinate (c1) at (-0.30,0);
\coordinate (c2) at (0.30,0);
\coordinate (c3) at (0.90,0);
\coordinate (s0) at (-0.60,1.12);
\coordinate (s1) at (0.00,1.12);
\coordinate (s2) at (0.60,1.12);
\draw (s0) -- (c0);
\draw (s0) -- (c1);
\draw (s0) -- (c2);
\draw (s1) -- (c3);
\draw (s2) -- (c3);
\fill (c0) circle (2pt) node[below=4pt] {$1$};
\fill (c1) circle (2pt) node[below=4pt] {$2$};
\fill (c2) circle (2pt) node[below=4pt] {$3$};
\fill (c3) circle (2pt) node[below=4pt] {$4$};
\draw[fill=white] (s0) circle (2.3pt);
\draw[fill=white] (s1) circle (2.3pt);
\draw[fill=white] (s2) circle (2.3pt);
\end{tikzpicture}
% End embedded drawing: 禁用子图/R3.tikz.tex
\caption*{$R_{3}$}
\end{subfigure}
\hfill
\begin{subfigure}[b]{0.23\textwidth}
\centering
% Begin embedded drawing: 禁用子图/R4.tikz.tex
\begin{tikzpicture}[x=1cm,y=1cm,scale=1,every node/.style={font=\scriptsize},line width=0.45pt]
\draw[rounded corners=3pt,gray] (-1.43,-0.26) rectangle (1.43,0.22);
\draw[rounded corners=3pt,gray] (-1.43,0.88) rectangle (1.43,1.36);
\node[anchor=east,inner sep=0pt] at (-1.53,-0.02) {$C$};
\node[anchor=east,inner sep=0pt] at (-1.53,1.12) {$Z$};
\coordinate (c0) at (-1.20,0);
\coordinate (c1) at (-0.60,0);
\coordinate (c2) at (0.00,0);
\coordinate (c3) at (0.60,0);
\coordinate (c4) at (1.20,0);
\coordinate (s0) at (-0.30,1.12);
\coordinate (s1) at (0.30,1.12);
\draw (s0) -- (c0);
\draw (s0) -- (c1);
\draw (s0) -- (c2);
\draw (s1) -- (c3);
\draw (s1) -- (c4);
\fill (c0) circle (2pt) node[below=4pt] {$1$};
\fill (c1) circle (2pt) node[below=4pt] {$2$};
\fill (c2) circle (2pt) node[below=4pt] {$3$};
\fill (c3) circle (2pt) node[below=4pt] {$4$};
\fill (c4) circle (2pt) node[below=4pt] {$5$};
\draw[fill=white] (s0) circle (2.3pt);
\draw[fill=white] (s1) circle (2.3pt);
\end{tikzpicture}
% End embedded drawing: 禁用子图/R4.tikz.tex
\caption*{$R_{4}$}
\end{subfigure}
\par\medskip
\begin{subfigure}[b]{0.23\textwidth}
\centering
% Begin embedded drawing: 禁用子图/R5.tikz.tex
\begin{tikzpicture}[x=1cm,y=1cm,scale=1,every node/.style={font=\scriptsize},line width=0.45pt]
\draw[rounded corners=3pt,gray] (-1.13,-0.26) rectangle (1.13,0.22);
\draw[rounded corners=3pt,gray] (-1.13,0.88) rectangle (1.13,1.36);
\node[anchor=east,inner sep=0pt] at (-1.23,-0.02) {$C$};
\node[anchor=east,inner sep=0pt] at (-1.23,1.12) {$Z$};
\coordinate (c0) at (-0.60,0);
\coordinate (c1) at (0.00,0);
\coordinate (c2) at (0.60,0);
\coordinate (s0) at (-0.90,1.12);
\coordinate (s1) at (-0.30,1.12);
\coordinate (s2) at (0.30,1.12);
\coordinate (s3) at (0.90,1.12);
\draw (s0) -- (c0);
\draw (s0) -- (c1);
\draw (s1) -- (c0);
\draw (s1) -- (c1);
\draw (s2) -- (c2);
\draw (s3) -- (c2);
\fill (c0) circle (2pt) node[below=4pt] {$1$};
\fill (c1) circle (2pt) node[below=4pt] {$2$};
\fill (c2) circle (2pt) node[below=4pt] {$3$};
\draw[fill=white] (s0) circle (2.3pt);
\draw[fill=white] (s1) circle (2.3pt);
\draw[fill=white] (s2) circle (2.3pt);
\draw[fill=white] (s3) circle (2.3pt);
\end{tikzpicture}
% End embedded drawing: 禁用子图/R5.tikz.tex
\caption*{$R_{5}$}
\end{subfigure}
\hfill
\begin{subfigure}[b]{0.23\textwidth}
\centering
% Begin embedded drawing: 禁用子图/R6.tikz.tex
\begin{tikzpicture}[x=1cm,y=1cm,scale=1,every node/.style={font=\scriptsize},line width=0.45pt]
\draw[rounded corners=3pt,gray] (-1.13,-0.26) rectangle (1.13,0.22);
\draw[rounded corners=3pt,gray] (-1.13,0.88) rectangle (1.13,1.36);
\node[anchor=east,inner sep=0pt] at (-1.23,-0.02) {$C$};
\node[anchor=east,inner sep=0pt] at (-1.23,1.12) {$Z$};
\coordinate (c0) at (-0.90,0);
\coordinate (c1) at (-0.30,0);
\coordinate (c2) at (0.30,0);
\coordinate (c3) at (0.90,0);
\coordinate (s0) at (-0.60,1.12);
\coordinate (s1) at (0.00,1.12);
\coordinate (s2) at (0.60,1.12);
\draw (s0) -- (c0);
\draw (s0) -- (c1);
\draw (s1) -- (c0);
\draw (s1) -- (c1);
\draw (s2) -- (c2);
\draw (s2) -- (c3);
\fill (c0) circle (2pt) node[below=4pt] {$1$};
\fill (c1) circle (2pt) node[below=4pt] {$2$};
\fill (c2) circle (2pt) node[below=4pt] {$3$};
\fill (c3) circle (2pt) node[below=4pt] {$4$};
\draw[fill=white] (s0) circle (2.3pt);
\draw[fill=white] (s1) circle (2.3pt);
\draw[fill=white] (s2) circle (2.3pt);
\end{tikzpicture}
% End embedded drawing: 禁用子图/R6.tikz.tex
\caption*{$R_{6}$}
\end{subfigure}
\hfill
\begin{subfigure}[b]{0.23\textwidth}
\centering
% Begin embedded drawing: 禁用子图/R7.tikz.tex
\begin{tikzpicture}[x=1cm,y=1cm,scale=1,every node/.style={font=\scriptsize},line width=0.45pt]
\draw[rounded corners=3pt,gray] (-1.13,-0.26) rectangle (1.13,0.22);
\draw[rounded corners=3pt,gray] (-1.13,0.88) rectangle (1.13,1.36);
\node[anchor=east,inner sep=0pt] at (-1.23,-0.02) {$C$};
\node[anchor=east,inner sep=0pt] at (-1.23,1.12) {$Z$};
\coordinate (c0) at (-0.60,0);
\coordinate (c1) at (0.00,0);
\coordinate (c2) at (0.60,0);
\coordinate (s0) at (-0.90,1.12);
\coordinate (s1) at (-0.30,1.12);
\coordinate (s2) at (0.30,1.12);
\coordinate (s3) at (0.90,1.12);
\draw (s0) -- (c0);
\draw (s0) -- (c1);
\draw (s1) -- (c0);
\draw (s2) -- (c2);
\draw (s3) -- (c2);
\fill (c0) circle (2pt) node[below=4pt] {$1$};
\fill (c1) circle (2pt) node[below=4pt] {$2$};
\fill (c2) circle (2pt) node[below=4pt] {$3$};
\draw[fill=white] (s0) circle (2.3pt);
\draw[fill=white] (s1) circle (2.3pt);
\draw[fill=white] (s2) circle (2.3pt);
\draw[fill=white] (s3) circle (2.3pt);
\end{tikzpicture}
% End embedded drawing: 禁用子图/R7.tikz.tex
\caption*{$R_{7}$}
\end{subfigure}
\hfill
\begin{subfigure}[b]{0.23\textwidth}
\centering
% Begin embedded drawing: 禁用子图/R8.tikz.tex
\begin{tikzpicture}[x=1cm,y=1cm,scale=1,every node/.style={font=\scriptsize},line width=0.45pt]
\draw[rounded corners=3pt,gray] (-1.13,-0.26) rectangle (1.13,0.22);
\draw[rounded corners=3pt,gray] (-1.13,0.88) rectangle (1.13,1.36);
\node[anchor=east,inner sep=0pt] at (-1.23,-0.02) {$C$};
\node[anchor=east,inner sep=0pt] at (-1.23,1.12) {$Z$};
\coordinate (c0) at (-0.90,0);
\coordinate (c1) at (-0.30,0);
\coordinate (c2) at (0.30,0);
\coordinate (c3) at (0.90,0);
\coordinate (s0) at (-0.60,1.12);
\coordinate (s1) at (0.00,1.12);
\coordinate (s2) at (0.60,1.12);
\draw (s0) -- (c0);
\draw (s0) -- (c1);
\draw (s1) -- (c0);
\draw (s2) -- (c2);
\draw (s2) -- (c3);
\fill (c0) circle (2pt) node[below=4pt] {$1$};
\fill (c1) circle (2pt) node[below=4pt] {$2$};
\fill (c2) circle (2pt) node[below=4pt] {$3$};
\fill (c3) circle (2pt) node[below=4pt] {$4$};
\draw[fill=white] (s0) circle (2.3pt);
\draw[fill=white] (s1) circle (2.3pt);
\draw[fill=white] (s2) circle (2.3pt);
\end{tikzpicture}
% End embedded drawing: 禁用子图/R8.tikz.tex
\caption*{$R_{8}$}
\end{subfigure}
\caption{The forbidden graphs.}
\label{fig:forbidden}
\end{figure}
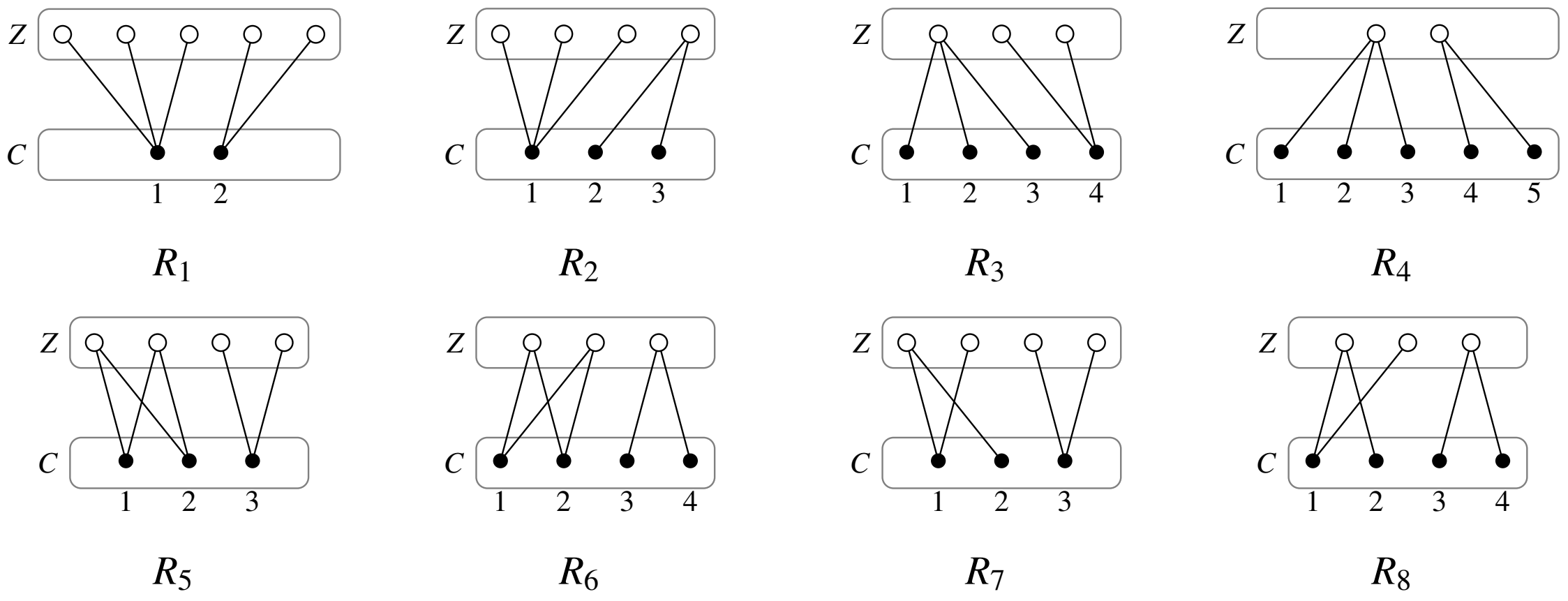
% End embedded drawing: figures/R_graphs.tex

The graphs $R_{1},\ldots,R_{8}$ in Figure~\ref{fig:forbidden} are split graphs with split partition $(C,Z)$.
By direct computation, we have the following lemma.

\begin{lemma}\label{lem:forbidden}
For $1\leq k\leq8$, $\lambda_{2}(R_{k})>1$.
\end{lemma}

If  $\lambda_{2}(G)\leq1$, then $G$ contains none of $R_{1},\ldots,R_{8}$ as an induced subgraph by Corollary~\ref{cor:interlacing} and Lemma~\ref{lem:forbidden}.

%\begin{lemma}\label{lem:chain}
%Let $G=(C,Z)$ be a distance-hereditary split graph with $\lambda_{2}(G)\leq1$.
%Let  $h$ be the number of distinct heavy neighborhoods $N(x)$ with $x\in Z$ and $h\geq1$.
%Then the following hold:
%\begin{enumerate}[label=\textup{(\arabic*)}]
%\item There are vertices $x_{0},\ldots,x_{h-1}\in Z$ such that the distinct heavy neighborhoods are $Q_{j}=N(x_{j})$ for $0\leq j\leq h-1$ and $Q_{0}\subset Q_{1}\subset\cdots\subset Q_{h-1}\subset C$.

%\item  For every $y\in Z$ whose neighborhood $N(y)$ is critical, $N(y)\subset Q_{0}$.
%\end{enumerate}
%\end{lemma}

\begin{lemma}\label{lem:chain}
Let $G=(C,Z)$ be a distance-hereditary split graph with $\lambda_{2}(G)\leq1$.
If $N(x)$   and  $N(y)$  are  heavy for  $x,y\in Z$, then either $N(x)\subseteq N(y)$ or $N(y)\subseteq N(x)$;
if  $N(x)$ is heavy and $N(y)$  is  critical, then $N(y)\subset N(x)$.
\end{lemma}
\begin{proof}
For  $x,y\in Z$, let $\rho(L_{x})>2$ and $\rho(L_{y})\geq2$.
If $|N(y)|=1$, then $L_{y}=[m_{y}]_{1\times1}$ and $m_{y}\geq2$.
If $|N(y)|=1$, let $W_{1}=\{t,z_{1},z_{2}\}$, where  $t\in N(y)$ and $z_{1},z_{2}\in Z$ ($z_{1}\neq z_{2}$) with $N(z_{1})=N(z_{2})=\{t\}$.
If $|N(y)|\geq2$, let $W_{2}=\{t_{1},t_{2},y\}$, where $t_{1},t_{2}$ be distinct vertices in $N(y)$.

Suppose that $N(x)\cap N(y)=\varnothing$.
If $|N(x)|=1$, then $L_{x}=[m_{x}]$ and $m_{x}\geq3$.
If $|N(x)|=1$,  then there are  $u\in N(x)$ and $z_{3},z_{4},z_{5}\in Z$  with $N(z_{3})=N(z_{4})=N(z_{5})=\{u\}$.
Thus, if $|N(x)|=1$, then either $G[W_{1}\cup\{u,z_{3},z_{4},z_{5}\}]\cong R_{1}$ (if $|N(y)|=1$) or $G[W_{2}\cup\{u,z_{3},z_{4},z_{5}\}]\cong R_{2}$ (if $|N(y)|\geq2$), a contradiction.
If $|N(x)|\geq3$ and  $u_{1},u_{2},u_{3}\in N(x)$,  then either $G[W_{1}\cup\{u_{1},u_{2},u_{3},x\}]\cong R_{3}$ (if $|N(y)|=1$) or $G[W_{2}\cup\{u_{1},u_{2},u_{3},x\}]\cong R_{4}$ (if $|N(y)|\geq2$), a contradiction.
If $|N(x)|=2$ and $m_{x}\geq2$, then either $G[W_{1}\cup\{u,v,z_{3},z_{4}\}]\cong R_{5}$ ($u,v\in N(x)$, $z_{3},z_{4}\in Z$ and $N(z_{3})=N(z_{4})=N(x)$) or $G[W_{2}\cup\{u,v,z_{3},z_{4}\}]\cong R_{6}$, a contradiction.
If $|N(x)|=2$ and $m_{x}=1$, then there are $u\in N(x)$ and $z\in Z$ with $N(z)=\{u\}$; otherwise $L_{x}=J_{2}$ and $\rho(L_{x})=2$, a contradiction.
Thus, if $|N(x)|=2$ and $m_{x}=1$, then either $G[W_{1}\cup\{u,v,x,z\}]\cong R_{7}$ ($u,v\in N(x)$) or $G[W_{2}\cup\{u,v,x,z\}]\cong R_{8}$, a contradiction.
Thus, $N(x)\cap N(y)\neq\varnothing$.
Thus, $N(x)\subseteq N(y)$ or $N(y)\subseteq N(x)$ by Lemma~\ref{lem:laminar}.

Let  $N(y_{1})$ be a critical neighborhood with $y_{1}\in Z$.
Then $N(y_{1})\cap N(x)\neq\varnothing$.
Since $\rho(L_{y_{1}})=2$ and $\rho(L_{x})>2$, $N(y_{1})\neq N(x)$.
If $N(x)\subset N(y_{1})$, then $|N(y_{1})|=2$ and $N(x)$ is a nonempty proper subset of $N(y_{1})$, a contradiction by Lemma~\ref{lem:light}.
Then $N(y_{1})\subset N(x)$ by Lemma~\ref{lem:laminar}.
\end{proof}

\begin{corollary}\label{cor:shape}
Let $G=(C,Z)$ be a connected distance-hereditary split graph with $|V(G)|\geq3$ and $\lambda_{2}(G)\leq1$.
Then $G\in \mathcal{D}_{1}$ or $G\cong D(\boldsymbol a,\boldsymbol b,\boldsymbol\mu;\kappa,\tau)$ ($\kappa+\tau>0$) or $G\cong D(\boldsymbol a,\boldsymbol b,\boldsymbol\mu;0,0)$ ($\mu_{0}(a_0+2b_0)>2$).
\end{corollary}

\begin{proof}
%Then $G$ contains none of $R_{1},\ldots,R_{8}$ as an induced subgraph by Corollary~\ref{cor:interlacing} and Lemma~\ref{lem:forbidden}.
Suppose that no heavy neighborhood exists.
By Lemma~\ref{lem:light}, every vertex of $Z$ has degree $1$ or $2$ in $G$.
Let $x,y\in Z$ satisfy $|N(x)|=|N(y)|=2$ and $N(x)\neq N(y)$.
If $N(x)\cap N(y)\neq\varnothing$, then $N(x)\subseteq N(y)$ or $N(y)\subseteq N(x)$ by Lemma~\ref{lem:laminar}.
Thus, if  $N(x)\cap N(y)\neq\varnothing$, then $N(x)=N(y)$, a contradiction.
Thus, the distinct neighborhoods of the vertices of $Z$ of degree $2$ in $G$ are pairwise disjoint.
Thus, $G\in \mathcal{D}_{1}$.

Suppose that there are $h$ distinct heavy neighborhoods ($h\geq1$).
Thus, by Lemma~\ref{lem:chain}, there are vertices $x_{0},\ldots,x_{h-1}\in Z$ such that $Q_{j}=N(x_{j})$ for $0\leq j<h$ and $Q_{0}\subset\cdots\subset Q_{h-1}$. Since $C$ is a maximum clique, $Q_{h-1}\subset C$.
Let~$C_{0}=Q_{0}$ and $C_{h}=C\setminus Q_{h-1}$.
For $1\leq j<h$, let $C_{j}=Q_{j}\setminus Q_{j-1}$.
Let $\mu_{j}=m_{x_{j}}$ for $0\leq j<h$.
If $y\in Z$ has a critical neighborhood, then $N(y)\subset Q_{0}$ by Lemma~\ref{lem:chain}. If $u\in C\setminus Q_{0}$ and $y\in Z$ with $N(y)=\{u\}$, then $N(y)$ is subcritical by Lemma~\ref{lem:chain}. Thus, $m_{y}=1$ for $\rho(L_{y})=m_{y}<2$. Thus, every vertex of $C\setminus Q_{0}$ has at most one pendant neighbor in $Z$.
If $y\in Z$ and $N(y)\subset Q_{0}$, then $N(y)$ is not heavy, for $Q_{0}$ is the smallest heavy neighborhood. Thus, if $y\in Z$ and $N(y)\subset Q_{0}$, then $N(y)$ is critical or subcritical.
If $y\in Z$ has a critical neighborhood $N(y)$ with $|N(y)|=1$, then $\rho(L_{y})=m_{y}=2$ by Lemma~\ref{lem:light}.
Thus, every critical neighborhood $N(y)$ with $y\in Z$ either satisfies $|N(y)|=1$ and $m_{y}=2$ or has the form in Lemma~\ref{lem:light}\textup{(ii)}.
Let $\kappa$ and $\tau$ be the numbers of distinct critical neighborhoods of sizes $1$ and $2$, respectively.
Thus, if there is a vertex $y\in Z$ whose neighborhood $N(y)$ is critical, then $\kappa+\tau>0$.
And  if no vertex of $Z$ has a critical neighborhood, then $\kappa=\tau=0$.
Then $G\cong D(\boldsymbol a,\boldsymbol b,\boldsymbol\mu;\kappa,\tau)$ by Lemmas~\ref{lem:light} and~\ref{lem:chain}.
If $\kappa+\tau=0$, then $\mu_{0}\sigma_{0}>2$ by Lemma~\ref{lem:base} and $N(x_{0})$ is heavy.
\end{proof}

\begin{lemma}\label{prop:criterion}
Let $G\cong D(\boldsymbol a,\boldsymbol b,\boldsymbol\mu;\kappa,\tau)$, where $\mu_{0}(a_{0}+2b_{0})>2$ if $\kappa+\tau=0$.
Then\\

 $\lambda_{2}(G)\leq1$ if and only if $p_{j}>0$ for $1\leq j<h$ and $p_{h}\geq0$ (where $p_{j}$ is defined in~Eq.\eqref{eq:recurrence}); 
If $\lambda_{2}(G)\leq1$ and $\kappa+\tau=0$, then $\mi(G)=\nobreak\left\{\begin{smallmatrix*}[l]1,&p_{h}=0;\\0,&p_{h}\neq0.\end{smallmatrix*}\right.;$
If $\lambda_{2}(G)\leq1$ and $\kappa+\tau>0$, then $\mi(G)=\nobreak\left\{\begin{smallmatrix*}[l]\kappa+\tau,&p_{h}=0;\\\kappa+\tau-1,&p_{h}\neq0.\end{smallmatrix*}\right.$.
\end{lemma}

\begin{proof}
Let $B$ be the submatrix of $A(G)$ in Lemma~\ref{lem:equ}. We consider the matrix $J_{|C|}-2I_{|C|}+BB^{\mathsf T}$.
Let $x_{0}\in X_{0}$ and $M_{0}=L_{x_{0}}-2I$.
For $1\leq j\leq h$, let $F_{j}=\diag(-2I_{a_{j}},-I_{b_{j}})$, $N_{j}=\diag(M_{j-1},F_{j})$ and $M_{j}=N_{j}+\mu_{j}J$.
Then $(F_{j})_{uu}=|N(u)\cap Y|-2$  and   $(F_{j})_{uv}=0$  for distinct vertices $u,v\in C_{j}$.
If $u,v\in Q_{j-1}$, then $(N_{j})_{uv}=(M_{j-1})_{uv}$.
If $u,v\in C_{j}$, then $(N_{j})_{uv}=(F_{j})_{uv}$.
If $u\in Q_{j-1}$ and $v\in C_{j}$, then $(N_{j})_{uv}=(N_{j})_{vu}=0$.
For $1\leq j<h$, $x\in X_{j}$ and $u\in C$, $(B)_{ux}=1$ if $u\in Q_{j}$, and $(B)_{ux}=0$ if $u\in C\setminus Q_{j}$.
Thus, for $1\leq j<h$ and $u,v\in Q_{j}$, $\sum_{x\in X_{j}}(B)_{ux}(B)_{vx}=\mu_{j}$.
For $u\in C$, $(N_{h})_{uu}=|N(u)\cap Z|-2=(BB^{\mathsf T})_{uu}-2$.
For distinct $u,v\in C$, $(N_{h})_{uv}=|N(u)\cap N(v)\cap Z|=(BB^{\mathsf T})_{uv}$.
Thus, $N_{h}=BB^{\mathsf T}-2I_{|C|}$.
Since $\mu_{h}=1$, $M_{h}=J_{|C|}-2I_{|C|}+BB^{\mathsf T}$.

Since all eigenvalues of $F_{j}$ are negative, $\npos(N_{j})=\npos(M_{j-1})$.
For every vector $x$, $x^{\mathsf T}M_{j}x=x^{\mathsf T}N_{j}x+\mu_{j}(\one^{\mathsf T}x)^2\geq x^{\mathsf T}N_{j}x$.
Then $\npos(M_{j})\geq\npos(M_{j-1})$ by Lemma~\ref{lem:courant-fischer}.

We next prove by induction on $j$ that if $p_{1},\ldots,p_{j}>0$, then $\npos(M_{j})=1$, $\ker M_{j}\perp\one$ and $\eta(M_{j})=c_{j}/p_{j}$.
The case $j=0$ holds by Lemma~\ref{lem:equ} and ~\ref{lem:base}.
For $1\leq j\leq h$ and $1\leq i<j$, suppose that $p_{i}>0$, $\npos(M_{j-1})=1$, $\ker M_{j-1}\perp\one$ and $\eta(M_{j-1})=c_{j-1}/p_{j-1}$.
Then $\npos(N_{j})=1$ and $\nzero(N_{j})=\nzero(M_{j-1})$.
Then $y=(x^{\mathsf T},\mathbf{0}^{\mathsf T})^{\mathsf T}$, where $y\in\ker N_{j}$ and $x\in\ker M_{j-1}$.
Thus, $\ker N_{j}\perp\one$.
Let $q$ be a  vector satisfying $M_{j-1}q=\one$.
Let $z=F_{j}^{-1}\one$.
Then $z=(-1/2,\ldots,-1/2,-1,\ldots,-1)^{\mathsf T}$, where $-1/2$ occurs $a_{j}$ times and $-1$ occurs $b_{j}$ times.
Then $N_{j}(q^{\mathsf T},z^{\mathsf T})^{\mathsf T}=\one$.
Thus, $\eta(N_{j})=\one^{\mathsf T}q+\one^{\mathsf T}z=\frac{c_{j-1}}{p_{j-1}}-\frac{a_{j}}{2}-b_{j}=\frac{c_{j-1}}{p_{j-1}}-\frac{\sigma_{j}}{2}=\frac{c_{j}}{2p_{j-1}}$ by~Eq.\eqref{eq:recurrence}.
Let $\Delta_{j}=1+\mu_{j}\eta(N_{j})$.
Then $\Delta_{j}=p_{j}/(2p_{j-1})$ by~Eq.\eqref{eq:recurrence}.
If $p_{j}>0$, then $\npos(M_{j})=1$ and $\nzero(M_{j})=\nzero(M_{j-1})$ by Lemma~\ref{lem:update}.
Thus, if $p_{j}>0$, then $\ker M_{j}\perp\one$ and $\eta(M_{j})=c_{j}/p_{j}$.
If at least one of $p_{1},\ldots,p_{h}$ is nonpositive, let $k=\min\{i:1\leq i\leq h,\ p_{i}\leq0\}$.
If $p_{k}<0$, then $\Delta_{k}<0$ and $\npos(M_{k})=2$ by Lemma~\ref{lem:update}, which implies that $\npos(M_{h})\geq2$.
If $p_{k}=0$ and $k<h$, then $\Delta_{k}=0$, $\npos(M_{k})=1$ and $\ker M_{k}\not\perp\one$ by Lemma~\ref{lem:update}.
Thus, if $p_{k}=0$ and $k<h$, then $\npos(M_{k+1})\geq2$ by Lemma~\ref{lem:mainzero}, which implies that $\npos(M_{h})\geq2$.
If $p_{h}=0$ and $p_{1},\ldots,p_{h-1}>0$, then $\npos(M_{h})=1$ and $\nzero(M_{h})=\nzero(M_{h-1})+1$ by Lemma~\ref{lem:update}.
Thus, $\lambda_{2}(G)\leq1$ if and only if $p_{j}>0$ for $1\leq j<h$ and $p_{h}\geq0$ by~Lemma~\ref{lem:equ}.
By Lemmas~\ref{lem:equ} and~\ref{lem:update}, if $\lambda_{2}(G)\leq1$, then $\mi(G)=\nzero(M_{h})=\left\{\begin{smallmatrix*}[l]\nzero(M_{0})+1,&p_{h}=0;\\\nzero(M_{0}),&p_{h}\neq0.\end{smallmatrix*}\right.$
By Lemma~\ref{lem:base}, $\nzero(M_{0})=0$ if $\kappa+\tau=0$, and $\nzero(M_{0})=\kappa+\tau-1$ if $\kappa+\tau>0$.
This completes the proof.
\end{proof}

\begin{lemma}\label{lem:depth}
Let $G\cong D(\boldsymbol a,\boldsymbol b,\boldsymbol\mu;\kappa,\tau)$, where $\mu_{0}(a_{0}+2b_{0})>2$  if $\kappa+\tau=0$, and $\lambda_{2}(G)\leq1$.
If $\kappa+\tau=0$, then $h\leq4$;
if $\kappa+\tau>0$, then $h\leq3$.
\end{lemma}

\begin{proof}
Suppose that $\lambda_{2}(G)\leq1$.
Then $p_{j}>0$ for $1\leq j<h$ and $p_{h}\geq0$ by Lemma~\ref{prop:criterion}.
Let $r_{j}=c_{j}/p_{j}$ for $0\leq j<h$.
For $1\leq j\leq h$, let $s_{j}=r_{j-1}-\sigma_{j}/2$ and $\Delta_{j}=1+\mu_{j}s_{j}$.
Then $\Delta_{j}=p_{j}/(2p_{j-1})$ by~Eq.\eqref{eq:recurrence}.
For $1\leq j<h$, $\Delta_{j}>0$ and $r_{j}=s_{j}/(1+\mu_{j}s_{j})$.
Since $\mu_{j},\sigma_{j}\geq1$, $s_{j}>-1/\mu_{j}\geq-1$ and $s_{j}\leq r_{j-1}-1/2$ for $1\leq j<h$.
Let $f(x)=x/(1+x)$ for $x>-1$.
For $-1<x<y$, $f(y)-f(x)=(y-x)/((1+x)(1+y))>0$.
For $1\leq j<h$, $f(s_{j})-r_{j}=(\mu_{j}-1)s_{j}^2/((1+s_{j})(1+\mu_{j}s_{j}))\geq0$.
Then $r_{j}\leq f(s_{j})\leq f(r_{j-1}-1/2)$ for $1\leq j<h$.

Suppose that $\kappa+\tau=0$.
Since $\mu_{0}\sigma_{0}\geq3$ and $\mu_{0}\geq1$, $r_{0}=\sigma_{0}/(\mu_{0}\sigma_{0}-2)\leq3$.
Suppose that $h\geq5$.
Then $p_{1},p_{2},p_{3},p_{4}>0$ and $p_{5}\geq0$ by Lemma~\ref{prop:criterion}.
Since $f(x)$ is increasing, $r_{1}\leq5/7$, $r_{2}\leq3/17$, $r_{3}\leq-11/23$ and $r_{4}\leq-45$.
Then $s_{5}=r_{4}-\sigma_{5}/2\leq-91/2$ and $\Delta_{5}=1+\mu_{5}s_{5}\leq-89/2<0$.
Then $p_{5}=2p_{4}\Delta_{5}<0$, a contradiction.
Thus, $h\leq4$.

Suppose that $\kappa+\tau>0$.
Since $\mu_{0}\geq1$, $r_{0}=1/\mu_{0}\leq1$.
Suppose that $h\geq4$.
Then $p_{1},p_{2},p_{3}>0$ and $p_{4}\geq0$ by Lemma~\ref{prop:criterion}.
Since $f(x)$ is increasing, $r_{1}\leq1/3$, $r_{2}\leq-1/5$ and $r_{3}\leq-7/3$.
Then $s_{4}=r_{3}-\sigma_{4}/2\leq-17/6$ and $\Delta_{4}=1+\mu_{4}s_{4}\leq-11/6<0$.
Then $p_{4}=2p_{3}\Delta_{4}<0$, a contradiction.
Thus, $h\leq3$.
\end{proof}

\vspace{1.2cm}
\renewcommand{\refname}{References}
% Final Section 4 body. Root inserts the finite-graph plates at the marked point.
% The supplementary files named below must accompany the submitted manuscript.

\section{Proof of Theorem \ref{thm:main}}

\begin{lemma}\label{lem:classification}
Let $G=(C,Z)$ be a connected distance-hereditary split graph with $|V(G)|\geq3$.
Then the following hold:
\begin{enumerate}[label=\textup{(\arabic*)}]
\item $\lambda_{2}(G)\leq1$ if and only if $G\in \mathcal{D}_{1}\cup \mathcal{D}_{2}\cup \mathcal{D}_{3}$.
\item If $G\in \mathcal{D}_{1}$, then $\mi(G)=\max\{\gamma+\delta-1,0\}$;
if $G\in \mathcal{D}_{2}$, then $\mi(G)=\nobreak\left\{\begin{smallmatrix*}[l]1,&p_{h}=0;\\0,&p_{h}\neq0.\end{smallmatrix*}\right.$;
if $G\in \mathcal{D}_{3}$, then $\mi(G)=\nobreak\left\{\begin{smallmatrix*}[l]\kappa+\tau,&p_{h}=0;\\\kappa+\tau-1,&p_{h}\neq0.\end{smallmatrix*}\right.$
\end{enumerate}
\end{lemma}

\begin{proof}[Proof of Lemma~\ref{lem:classification}]
Suppose that $\lambda_{2}(G)\leq1$.
%If $G=K_{2}$, then $G=D(2,0,0,0)\in \mathcal{D}_{1}$, $\lambda_{2}(G)=-1$ and $\mi(G)=1$.
%Let $|V(G)| \geq3$.
By Corollary~\ref{cor:shape}, Lemma~\ref{prop:criterion} and Lemma~\ref{lem:depth}, $G\in \mathcal{D}_{1}\cup \mathcal{D}_{2}\cup \mathcal{D}_{3}$.
If $G\in \mathcal{D}_{1}\cup \mathcal{D}_{2}\cup \mathcal{D}_{3}$, then $G$ is a distance-hereditary split graph  by Lemma~\ref{lem:laminar}.
If $G\in \mathcal{D}_{2}\cup \mathcal{D}_{3}$, then $\lambda_{2}(G)\leq1$ and the stated values of $\mi(G)$ follow from Lemma~\ref{prop:criterion}.

Suppose that $G=D(\alpha,\beta,\gamma,\delta)\in \mathcal{D}_{1}$.
Let $k=\nobreak\left\{\begin{smallmatrix*}[l]\gamma,&\gamma+\delta>0;\\1,&\gamma+\delta=0.\end{smallmatrix*}\right.$
Let $H=D\bigl((\alpha,1),\allowbreak(\beta,0),\allowbreak(1);k,\delta\bigr)$.
Then $G\subset H\in \mathcal{D}_{3}$.
Thus, $\lambda_{2}(H)\leq1$ and $\mi(H)=k+\delta-1=\max\{\gamma+\delta-1,0\}$ by Lemma~\ref{prop:criterion}.
Then $\lambda_{2}(G)\leq1$ by Lemma~\ref{lem:interlacing}.
By Lemma~\ref{lem:courant-fischer}, $\lambda_{1}(G)\geq 2|E(G)|/|V(G)|>1$.
Thus, $\mi(G)\leq\mi(H)$  by Lemma~\ref{lem:interlacing}.
Thus, if $\gamma+\delta=0$, then $\mi(G)=\mi(H)=0$.
Let $\gamma+\delta>0$.
For $1\leq i\leq\gamma$, let $u_{i}\in C_{2}$ and let $u_{i}',u_{i}''$ be its two pendant neighbors.
For $1\leq j\leq\delta$, let $v_{j},w_{j}\in C_{3}$ and  $z_{j}$ be their common neighbor.
Let $t_{1},\ldots,t_{\gamma},s_{1},\ldots,s_{\delta}\in\mathbb R$ satisfy $\sum_{i=1}^{\gamma}t_{i}+2\sum_{j=1}^{\delta}s_{j}=0$.
Let $x\in\mathbb R^{|V(G)|}$ be the column vector with $x_{u_{i}}=x_{u_{i}'}=x_{u_{i}''}=t_{i}$ for $1\leq i\leq\gamma$, $x_{v_{j}}=x_{w_{j}}=s_{j}$ and $x_{z_{j}}=2s_{j}$ for $1\leq j\leq\delta$, and $x_{v}=0$ for every other vertex $v$ of $G$.
Then $A(G)x=x$.
Thus, $\mi(G)\geq\gamma+\delta-1$.
Therefore, $\mi(G)=\max\{\gamma+\delta-1,0\}$.
\end{proof}

The multiplicity of the eigenvalue $1$ is given by Lemma~\ref{lem:classification}.
We now give the proof of Theorem~\ref{thm:main}.
Let $n$ be a sufficiently large positive integer.

\begin{proof}[Proof of Theorem~\ref{thm:main}]
By Eq.\eqref{eq:recurrence}, the graphs $A_{i}\in \mathcal{D}_{2}$ ($1\leq i\leq355$), $\mathcal F_{j}(n)\in \mathcal{D}_{2}$  ($1\leq j\leq82$), and $\mathcal F_{j}(n)\in \mathcal{D}_{3}$ ($83\leq j\leq94$).
By Lemma~\ref{lem:classification}, it remains to prove that every graph in $\mathcal{D}_{1}\cup \mathcal{D}_{2}\cup \mathcal{D}_{3}$ is an induced subgraph of a graph in $\mathcal H$.
Every graph in $\mathcal{D}_{1}$ is an induced subgraph of $D\bigl((0,1),(0,0),(1);n,n\bigr)\in \mathcal{D}_{3}$.
Thus, it suffices to consider $G\cong D(\boldsymbol a,\boldsymbol b,\boldsymbol\mu;\kappa,\tau)\in \mathcal{D}_{2}\cup \mathcal{D}_{3}$.
Let $\boldsymbol\sigma=(\sigma_{0},\sigma_{1},\ldots,\sigma_{h})$, where $\sigma_{j}=a_{j}+2b_{j}$ ($0\leq j\leq h$).
Since $C_{j}$ ($0\leq j\leq h$) is nonempty, $\sigma_{0}\geq0$ and $\sigma_{i}\geq1$ ($1\leq i\leq h$).
Since $0\leq b_{j}\leq\lfloor\sigma_{j}/2\rfloor$ and $a_{j}=\sigma_{j}-2b_{j}$, the possible graphs are determined by the ranges of $\boldsymbol\sigma$, $\boldsymbol\mu$ and the choices of $\kappa,\tau$.

Let $r_{j}=c_{j}/p_{j}$ for $0\leq j<h$.
Then  $r_{j}=\frac{r_{j-1}-\sigma_{j}/2}{1+\mu_{j}(r_{j-1}-\sigma_{j}/2)}$ and $p_{h}=2p_{h-1}\left(1+r_{h-1}-\frac{\sigma_{h}}{2}\right)$ by~Eq.\eqref{eq:recurrence}.
Then $p_{j}>0$ if and only if  $1+\mu_{j}\left(r_{j-1}-\frac{\sigma_{j}}{2}\right)>0$.
Then $p_{h}\geq0$ if and only if   $r_{h-1}\geq\sigma_{h}/2-1$.
We distinguish the following two cases.

\noindent$\textbf{Case 1.}$ $G\in \mathcal{D}_{3}$.

Then $r_{0}=1/\mu_{0}$.
Without loss of generality, let $a_{0}=b_{0}=\sigma_{0}=0$ and $\kappa=\tau=n$, where $n$ is a sufficiently large positive integer.
Since $1\leq h\leq 3$, we distinguish the following three subcases.

\noindent$\textbf{Subcase 1.1.}$ $h=1$.

Then  $\mu_{0}(\sigma_{1}-2)\leq2$ for $r_{0}\geq\sigma_{1}/2-1$ and $r_{0}=1/\mu_{0}$.
Thus, if $\mu_{0}=1$, then $1\leq\sigma_{1}\leq4$; if $\mu_{0}=2$, then $1\leq\sigma_{1}\leq3$; if $\mu_{0}\geq3$, then $1\leq\sigma_{1}\leq2$.

Suppose that $\mu_{0}=1$.
Then $1\leq a_{1}+2b_{1}=\sigma_{1}\leq4$, which gives $0\leq b_{1}\leq2$.
If $b_{1}=0$, then $1\leq a_{1}\leq4$ and $G\subset D\bigl((0,4),(0,0),(1);n,n\bigr)\subset\mathcal F_{83}(n)$.
If $b_{1}=1$, then $0\leq a_{1}\leq2$ and $G\subset D\bigl((0,2),(0,1),(1);n,n\bigr)\subset\mathcal F_{84}(n)$.
If $b_{1}=2$, then $a_{1}=0$ and $G\cong D\bigl((0,0),(0,2),(1);n,n\bigr)\subset\mathcal F_{86}(n)$.

Suppose that $\mu_{0}=2$.
Then $1\leq a_{1}+2b_{1}=\sigma_{1}\leq3$, which gives $0\leq b_{1}\leq1$.
If $b_{1}=0$, then $1\leq a_{1}\leq3$ and $G\subset D\bigl((0,3),(0,0),(2);n,n\bigr)\subset\mathcal F_{87}(n)$.
If $b_{1}=1$, then $0\leq a_{1}\leq1$ and $G\subset D\bigl((0,1),(0,1),(2);n,n\bigr)\subset\mathcal F_{88}(n)$.

Suppose that $\mu_{0}\geq3$.
Let $n\geq\mu_{0}$.
Then $1\leq a_{1}+2b_{1}=\sigma_{1}\leq2$, which gives $0\leq b_{1}\leq1$.
If $b_{1}=0$, then $1\leq a_{1}\leq2$ and $G\subset D\bigl((0,2),(0,0),(n);n,n\bigr)\subset\mathcal F_{83}(n)$.
If $b_{1}=1$, then $a_{1}=0$ and $G\subset D\bigl((0,0),(0,1),(n);n,n\bigr)\subset\mathcal F_{84}(n)$.

\noindent$\textbf{Subcase 1.2.}$ $h=2$.

Then $1+\mu_{1}(1/\mu_{0}-\sigma_{1}/2)>0$, $r_{1}=\frac{r_{0}-\sigma_{1}/2}{1+\mu_{1}(r_{0}-\sigma_{1}/2)}$ and $r_{1}\geq\sigma_{2}/2-1$.
Thus, $\sigma_{1}/2\leq 1/\mu_{0}+1/(\mu_{1}+2)$.
Thus, if $\mu_{0}=1$, then $1\leq\sigma_{1}\leq2$; if $\mu_{0}\geq2$, then $\sigma_{1}=1$.
If $\mu_{0}=1$, then $0\leq r_{1}=\frac{1-\sigma_{1}/2}{1+\mu_{1}(1-\sigma_{1}/2)}\leq1/3$, which implies that $1\leq\sigma_{2}\leq2$.
If $\mu_{0}=2$, then $r_{1}=0$, which implies that $1\leq\sigma_{2}\leq2$.
If $\mu_{0}\geq3$, then $r_{1}<0$, which implies that $\sigma_{2}=1$.
If $\mu_{0}\geq3$, then $(\mu_{0}-2)\mu_{1}\leq4$ by $\sigma_{1}=1$ and $\sigma_{1}/2\leq 1/\mu_{0}+1/(\mu_{1}+2)$.

% Begin embedded drawing: figures/F83_F84_F86.tex
% Requires tikz and the decorations.pathreplacing library (already in the main file).
% Solid vertices form the clique C; its edges are omitted. An edge from an
% X-set to a C-block means complete adjacency between the two sets.
\begin{figure}[!htbp]
\centering
\begingroup
\newcommand{\Fpanel}[3]{%
\begin{tikzpicture}[x=0.9cm,y=0.9cm,line width=0.5pt,every node/.style={font=\small}]
\draw[rounded corners=5pt,gray] (-0.2,-0.45) rectangle (8.15,0.95);
\node[anchor=east] at (-0.32,0.2) {$C$};
\draw[rounded corners=5pt,gray] (-0.2,-0.85) rectangle (8.15,-2.15);
\node[anchor=east] at (-0.32,-1.45) {$Y$};
\node[draw,rounded corners=3pt,gray,minimum width=37mm,minimum height=11mm,inner sep=0pt] (C0) at (2.05,0.2) {};
\node[anchor=north east,xshift=-2pt,yshift=-2pt] at (C0.north east) {$C_0$};
\node[draw,rounded corners=3pt,gray,minimum width=13mm,minimum height=11mm,inner sep=0pt] (C1) at (5,0.2) {};
\node[anchor=north east,xshift=-2pt,yshift=-2pt] at (C1.north east) {$C_1$};
\node[draw,rounded corners=3pt,gray,minimum width=13mm,minimum height=11mm,inner sep=0pt] (C2) at (7.15,0.2) {};
\node[anchor=north east,xshift=-2pt,yshift=-2pt] at (C2.north east) {$C_2$};
\node[draw,ellipse,fill=white,minimum width=12mm,minimum height=7mm] (X0) at (1.15,1.65) {$X_0$};
\node[above] at (X0.north) {$|X_0|=1$};
\node[draw,ellipse,fill=white,minimum width=12mm,minimum height=7mm] (X1) at (4.65,1.65) {$X_1$};
\node[above] at (X1.north) {$|X_1|=n$};
\draw (X0.south)--(C0.north);
\draw (X1.south)--(C0.north);
\draw (X1.south)--(C1.north);
% n disjoint pairs in C_0, each with a common degree-two neighbour.
\draw (0.18,0)--(0.37,-1.2)--(0.56,0);
\draw (1.25,0)--(1.44,-1.2)--(1.63,0);
\foreach \x in {0.18,0.56,1.25,1.63} \fill (\x,0) circle (2pt);
\foreach \x in {0.37,1.44} \draw[fill=white] (\x,-1.2) circle (2.4pt);
\node at (0.9,0) {$\cdots$}; \node at (0.9,-1.2) {$\cdots$};
\draw[decorate,decoration={brace,mirror,amplitude=3pt}] (0.37,-1.48)--(1.44,-1.48) node[midway,below=4pt] {$n$};
% n further C_0 vertices, each with two pendant neighbours.
\foreach \x in {2.45,3.65} {
  \draw (\x,0)--(\x-0.18,-1.2) (\x,0)--(\x+0.18,-1.2);
  \fill (\x,0) circle (2pt);
  \draw[fill=white] (\x-0.18,-1.2) circle (2.4pt);
  \draw[fill=white] (\x+0.18,-1.2) circle (2.4pt);
}
\node at (3.05,0) {$\cdots$}; \node at (3.05,-1.2) {$\cdots$};
\draw[decorate,decoration={brace,mirror,amplitude=3pt}] (2.45,-1.48)--(3.65,-1.48) node[midway,below=4pt] {$n$};
% A block with argument 2 contains two plain vertices; argument 1 has
% one vertex with one pendant neighbour.
\ifnum#1=2 \fill (4.75,0) circle (2pt) (5.25,0) circle (2pt);
\else \draw (5,0)--(5,-1.2); \fill (5,0) circle (2pt); \draw[fill=white] (5,-1.2) circle (2.4pt); \fi
\ifnum#2=2 \fill (6.9,0) circle (2pt) (7.4,0) circle (2pt);
\else \draw (7.15,0)--(7.15,-1.2); \fill (7.15,0) circle (2pt); \draw[fill=white] (7.15,-1.2) circle (2.4pt); \fi
\node at (4,-2.55) {$\mathcal F_{#3}(n)$};
\end{tikzpicture}%
}
\resizebox{\textwidth}{!}{%
\Fpanel{2}{2}{83}\hspace{3mm}%
\Fpanel{2}{1}{84}%
}
\endgroup
\caption{The graphs $\mathcal F_{83}(n)$ and $\mathcal F_{84}(n)$.}
\label{fig:F83-F84-F86}
\end{figure}
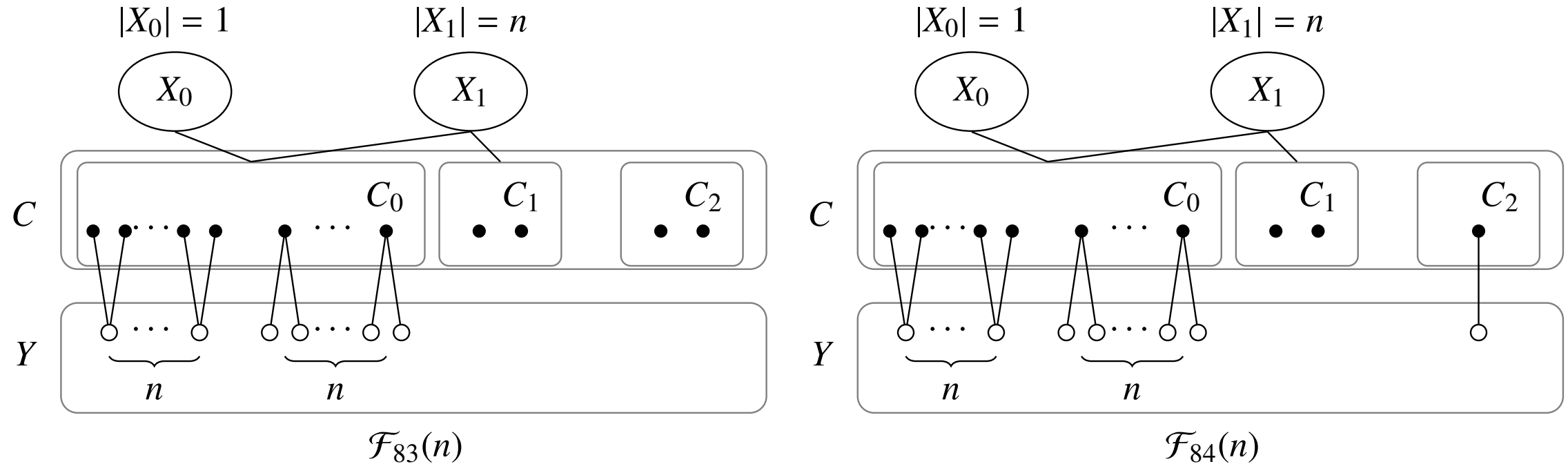
% End embedded drawing: figures/F83_F84_F86.tex

Suppose that $\mu_{0}=1$.
Let $n\geq\mu_{1}$.
Then $1\leq a_{j}+2b_{j}=\sigma_{j}\leq2$ for $j=1,2$.
Thus, $0\leq b_{1},b_{2}\leq1$.
If~$(b_{1},b_{2})=(0,0)$, then $G\subset D\bigl((0,2,2),(0,0,0),(1,n);n,n\bigr)=\mathcal F_{83}(n)$ (see Figure~\ref{fig:F83-F84-F86}).
If~$(b_{1},b_{2})=(0,1)$, then $G\subset D\bigl((0,2,0),(0,0,1),(1,n);n,n\bigr)=\mathcal F_{84}(n)$ (see Figure~\ref{fig:F83-F84-F86}).
If~$(b_{1},b_{2})=(1,0)$, then $G\subset D\bigl((0,0,2),(0,1,0),(1,n);n,n\bigr)=\mathcal F_{85}(n)$.
If~$(b_{1},b_{2})=(1,1)$, then $G\subset D\bigl((0,0,0),(0,1,1),(1,n);n,n\bigr)=\mathcal F_{86}(n)$.

Suppose that $\mu_{0}=2$.
Let $n\geq\mu_{1}$.
Then $a_{1}=1$, $b_{1}=0$ and $1\leq a_{2}+2b_{2}=\sigma_{2}\leq2$.
If $b_{2}=0$, then $1\leq a_{2}\leq2$ and $G\subset D\bigl((0,1,2),(0,0,0),(2,n);n,n\bigr)=\mathcal F_{87}(n)$.
If $b_{2}=1$, then $a_{2}=0$ and $G\subset D\bigl((0,1,0),(0,0,1),(2,n);n,n\bigr)=\mathcal F_{88}(n)$.

Suppose that $\mu_{0}\geq3$.
Then $a_{1}=a_{2}=1$, $b_{1}=b_{2}=0$ and $(\mu_{0}-2)\mu_{1}\leq4$.
If $\mu_{0}=3$, then $1\leq\mu_{1}\leq4$ and $G\subset D\bigl((0,1,1),(0,0,0),(3,4);n,n\bigr)=\mathcal F_{89}(n)$.
If $\mu_{0}=4$, then $1\leq\mu_{1}\leq2$ and $G\subset D\bigl((0,1,1),(0,0,0),(4,2);n,n\bigr)=\mathcal F_{90}(n)$.
If $5\leq\mu_{0}\leq6$, then $\mu_{1}=1$ and $G\subset D\bigl((0,1,1),(0,0,0),(6,1);n,n\bigr)=\mathcal F_{91}(n)$.

\noindent$\textbf{Subcase 1.3.}$ $h=3$.

Then $1+\mu_{j}(r_{j-1}-\sigma_{j}/2)>0$ for $j=1,2$ and $r_{2}\geq\sigma_{3}/2-1\geq-1/2$.
Since $r_{2}=\frac{r_{1}-\sigma_{2}/2}{1+\mu_{2}(r_{1}-\sigma_{2}/2)}\geq-1/2$, $r_{1}\geq\sigma_{2}/2-1/(\mu_{2}+2)>0$.
Since $r_{1}=\frac{1/\mu_{0}-\sigma_{1}/2}{1+\mu_{1}(1/\mu_{0}-\sigma_{1}/2)}$,  $1/\mu_{0}-\sigma_{1}/2>0$ , which implies that $\mu_{0}=1$ and $\sigma_{1}=1$.
Then $r_{1}=1/(\mu_{1}+2)\leq1/3$.
Thus, $\frac13 \geq r_{1}\geq \sigma_{2}/2-1/(\mu_{2}+2)\geq \frac{\sigma_{2}}{2}-\frac13$, which implies that $\sigma_{2}=1$.
Then $r_{2}=\frac{r_{1}-1/2}{1+\mu_{2}(r_{1}-1/2)}<0$.
Since $r_{2}\geq\sigma_{3}/2-1$, $\sigma_{3}=1$.
Since $r_{1}=1/(\mu_{1}+2)\geq\sigma_{2}/2-1/(\mu_{2}+2)$, $\mu_{1}\mu_{2}\leq4$.
Thus, $a_{1}=a_{2}=a_{3}=1$ and $b_{1}=b_{2}=b_{3}=0$.
If $\mu_{1}=1$, then $1\leq\mu_{2}\leq4$ and $G\subset D\bigl((0,1,1,1),(0,0,0,0),(1,1,4);n,n\bigr)=\mathcal F_{92}(n)$.
If $\mu_{1}=2$, then $1\leq\mu_{2}\leq2$ and $G\subset D\bigl((0,1,1,1),(0,0,0,0),(1,2,2);n,n\bigr)=\mathcal F_{93}(n)$.
If $3\leq\mu_{1}\leq4$, then $\mu_{2}=1$ and $G\subset D\bigl((0,1,1,1),(0,0,0,0),(1,4,1);n,n\bigr)=\mathcal F_{94}(n)$.
Thus, every graph in $\mathcal{D}_{3}$ is covered.

\noindent$\textbf{Case 2.}$ $G\in \mathcal{D}_{2}$.

Then $\mu_{0}\sigma_{0}\geq3$ and $1\leq h\leq4$.
Thus, $0<r_{0}=\sigma_{0}/(\mu_{0}\sigma_{0}-2)\leq\mu_{0}\sigma_{0}/(\mu_{0}\sigma_{0}-2)\leq3$.
If $h\leq3$ and the parameters $\sigma_{1},\ldots,\sigma_{h}$ and $\mu_{0},\ldots,\mu_{h-1}$ satisfy the conditions for $\mathcal{D}_{3}$ with $r_{0}=1/\mu_{0}$, then $G\subset D(\boldsymbol a,\boldsymbol b,\boldsymbol\mu;1,0)\in \mathcal{D}_{3}$.
For $h\leq3$, suppose that these parameters do not satisfy the conditions for $\mathcal{D}_{3}$ with $r_{0}=1/\mu_{0}$.

For $1\leq j<h$, since
$1+\mu_{j}(r_{j-1}-\sigma_{j}/2)>0$, we have $\frac{1}{\mu_{j}}-r_{j}=\frac{1}{\mu_{j}\bigl[1+\mu_{j}(r_{j-1}-\sigma_{j}/2)\bigr]}>0$, which implies that  $r_{j}<1/\mu_{j}$.
If $r_{j}\geq-1/2$, then $r_{j-1}\geq\sigma_{j}/2-1/(\mu_{j}+2)$.

\noindent$\textbf{Subcase 2.1.}$ $h=1$.

Then $r_{0}\geq\sigma_{1}/2-1$, and hence $\sigma_{1}\leq2+2r_{0}$.
Since these parameters do not satisfy the conditions for $\mathcal{D}_{3}$ with $r_{0}=1/\mu_{0}$, $1/\mu_{0}<\sigma_{1}/2-1$, which implies that $\sigma_{1}>2+2/\mu_{0}$.
If $\mu_{0}=1$, then $\sigma_{1}\geq5$ and $\sigma_{0}\geq3$ for $\mu_{0}\sigma_{0}\geq3$.
Thus, if $\mu_{0}=1$, then $r_{0}=\sigma_{0}/(\sigma_{0}-2)\geq3/2$, which gives  $3\leq\sigma_{0}\leq6$.
Thus, if $\mu_{0}=1$, then $5\leq\sigma_{1}\leq\lfloor2+2\sigma_{0}/(\sigma_{0}-2)\rfloor$.
If $\mu_{0}=2$, then $\sigma_{0}\geq2$ and $\sigma_{1}\geq4$ for $\mu_{0}\sigma_{0}\geq3$ and $\sigma_{1}>2+2/\mu_{0}$.
Thus, if $\mu_{0}=2$, then $r_{0}\geq\sigma_{1}/2-1\geq1$.
Thus, if $\mu_{0}=2$, then  $r_{0}=\sigma_{0}/(2\sigma_{0}-2)\geq1$ gives $\sigma_{0}\leq2$.
Thus, if $\mu_{0}=2$, then $\sigma_{0}=2$, $r_{0}=1$, and $\sigma_{1}=4$.
If $\mu_{0}\geq3$, then $\sigma_{1}\geq3$, and hence $r_{0}\geq1/2$.
Thus, if $\mu_{0}\geq3$, then  $r_{0}=\sigma_{0}/(\mu_{0}\sigma_{0}-2)\geq1/2$, which implies that $(\mu_{0}-2)\sigma_{0}\leq2$.
If $\mu_{0}=3$, then $1\leq\sigma_{0}\leq2$.
If $\mu_{0}=3$ and $\sigma_{0}=1$, then $r_{0}=1$, and hence $3\leq\sigma_{1}\leq4$; if $\mu_{0}=3$ and  $\sigma_{0}=2$, then $r_{0}=1/2$, and hence $\sigma_{1}=3$.
If $\mu_{0}=4$, then $\sigma_{0}=1$, $r_{0}=1/2$ and $\sigma_{1}=3$.
If $\mu_{0}\geq5$, then $(\mu_{0}-2)\sigma_{0}\geq3$, a contradiction.

Suppose that $\mu_{0}=1$.
If $\sigma_{0}=3,4,5$ or $6$, then $\sigma_{1}\leq8,6,5$ or $5$, respectively, for $\sigma_{1}\leq\lfloor 2+2\sigma_{0}/(\sigma_{0}-2)\rfloor$.
Note that $0\leq b_{j}\leq\lfloor\sigma_{j}/2\rfloor$.
If $\sigma_{0}=3$, then $G\subset D\bigl((3-2b_{0},8-2b_{1}),(b_{0},b_{1}),(1);0,0\bigr)\subset\mathcal F_{i}(n)$ for some $1\leq i\leq16$.
If $\sigma_{0}=4$, then $G\subset D\bigl((4-2b_{0},6-2b_{1}),(b_{0},b_{1}),(1);0,0\bigr)\subset\mathcal F_{i}(n)$ for some $17\leq i\leq34$.
If $5\leq\sigma_{0}\leq6$, then $G\subset D\bigl((6-2b_{0},5-2b_{1}),(b_{0},b_{1}),(1);0,0\bigr)\subset\mathcal F_{i}(n)$ for some $35\leq i\leq50$.

Suppose that $\mu_{0}=2$.
Then $\sigma_{0}=2$ and $\sigma_{1}=4$.
Then $0\leq b_{0}\leq1$ and $0\leq b_{1}\leq2$.
Thus, $G\cong D\bigl((2-2b_{0},4-2b_{1}),(b_{0},b_{1}),(2);0,0\bigr)\subset\mathcal F_{i}(n)$ for some $83\leq i\leq86$.

Suppose that $\mu_{0}=3$.
If  $\sigma_{0}=1$ or $2$, then $3\leq\sigma_{1}\leq4$ or $\sigma_{1}=3$ , respectively.
Thus, if $\sigma_{0}=1$, then $G\subset D\bigl((1,4-2b_{1}),(0,b_{1}),(3);0,0\bigr)\subset\mathcal F_{i}(n)$ for some $83\leq i\leq86$.
If $\sigma_{0}=2$, then $G\cong D\bigl((2-2b_{0},3-2b_{1}),(b_{0},b_{1}),(3);0,0\bigr)\subset\mathcal F_{i}(n)$ for some $87\leq i\leq88$.

Suppose that $\mu_{0}=4$.
Then $\sigma_{0}=1$ and $\sigma_{1}=3$.
Thus, we have $G\cong D\bigl((1,3-2b_{1}),(0,b_{1}),(4);0,0\bigr)\subset\mathcal F_{i}(n)$ for some $87\leq i\leq88$.

\noindent$\textbf{Subcase 2.2.}$ $h=2$.

Then $1+\mu_{1}(r_{0}-\sigma_{1}/2)>0$ and $r_{1}\geq\sigma_{2}/2-1$.
Since $r_{1}<1/\mu_{1}\leq1$, $\sigma_{2}<4$.
Thus, $1\leq\sigma_{2}\leq3$.

\begin{fact}\label{f221}
If $\sigma_{2}=1$, then $r_{0}\geq\sigma_{1}/2-1/(\mu_{1}+2)$ and $\sigma_{1}/2-1/(\mu_{1}+2)>1/\mu_{0}$.
If $\sigma_{2}=2$, then $r_{0}\geq\sigma_{1}/2$ and $\sigma_{1}\mu_{0}>2$.
If $\sigma_{2}=3$, then $\mu_{1}=1$ and $r_{0}\geq\sigma_{1}/2+1$.
\end{fact}
\begin{proof}
If $\sigma_{2}=1$, then $r_{1}\geq-1/2$.
Note that if $r_{j}\geq-1/2$, then $r_{j-1}\geq\sigma_{j}/2-1/(\mu_{j}+2)$.
Thus, if $\sigma_{2}=1$, then $r_{0}\geq\sigma_{1}/2-1/(\mu_{1}+2)$.
If $\sigma_{2}=2$, then $r_{1}\geq0$.
Thus, if $\sigma_{2}=2$, then  $r_{0}\geq\sigma_{1}/2$ for $r_1=\frac{r_0-\sigma_1/2}{1+\mu_1(r_0-\sigma_1/2)}$.
If $\sigma_{2}=3$, then $r_{1}\geq1/2$.
Thus, if $\sigma_{2}=3$, then $\mu_{1}=1$ for $r_{1}<1/\mu_{1}$.
Thus, if $\sigma_{2}=3$, then $r_{0}\geq\sigma_{1}/2+1$ for $r_1=\frac{r_0-\sigma_1/2}{1+\mu_1(r_0-\sigma_1/2)}$.
If $\sigma_{2}=1$, then $\sigma_{1}/2-1/(\mu_{1}+2)>1/\mu_{0}$ for the parameters do not satisfy the conditions for $\mathcal{D}_{3}$ with $r_{0}=1/\mu_{0}$.
Similarly, if $\sigma_{2}=2$, then $\sigma_{1}\mu_{0}>2$.
If $\sigma_{2}=3$, then $r_{1}\geq\sigma_{2}/2-1=1/2$.
Since $r_{1}<1/\mu_{1}$, we have $1/2\leq r_{1}<1/\mu_{1}$, which implies that $\mu_{1}=1$.
\end{proof}

Then $r_{0}\geq\sigma_{1}/2-1/3\geq1/6$ by Fact~\ref{f221}.
Then $\mu_{0}-2/\sigma_{0}=1/r_{0}\leq6$, and hence $1\leq\mu_{0}\leq8$.
Then $1\leq\sigma_{1}\leq6$ for $r_{0}\leq3$.

\begin{fact}\label{f22bounds}
If $\mu_{0}=1,2,3,4,5,6,7$ or $8$, then $r_{0}\geq7/6,2/3,5/14,3/10,1/4,1/4,1/6$ or $1/6$, and $\sigma_{0}\leq14,4,10,3,2,1,2$ or $1$, respectively.
\end{fact}
\begin{proof}
Suppose that $1\leq\mu_{0}\leq2$.
If $\mu_{0}=1$ and $\sigma_{2}=1$, then $\sigma_{1}/2-1/(\mu_{1}+2)>1$ gives $\sigma_{1}\geq3$.
Thus, if $\mu_{0}=1$ and $\sigma_{2}=1$, then $r_{0}\geq\sigma_{1}/2-1/(\mu_{1}+2)\geq3/2-1/3=7/6$.
If $\mu_{0}=1$ and $\sigma_{2}=2$, then $\sigma_{1}\geq3$ for $\sigma_{1}\mu_{0}>2$.
Thus, if $\mu_{0}=1$ and $\sigma_{2}=2$, then $r_{0}\geq\sigma_{1}/2\geq3/2$.
If $\mu_{0}=1$ and $\sigma_{2}= 3$, then $r_{0}\geq \sigma_{1}/2+1 \geq 3/2$.
If $\mu_{0}=2$ and $\sigma_{2}=1$, then $\sigma_{1}/2-1/(\mu_{1}+2)>1/2$ gives $\sigma_{1}\geq2$.
Thus, if $\mu_{0}=2$ and $\sigma_{2}=1$, then $r_{0}\geq\sigma_{1}/2-1/(\mu_{1}+2)\geq2/3$.
If $\mu_{0}=2$ and $\sigma_{2}=2$, then $\sigma_{1}\geq2$ for $\sigma_{1}\mu_{0}>2$.
Thus, $r_{0}\geq\sigma_{1}/2\geq1$.
If $\mu_{0}=2$ and $\sigma_{2}=3$, then $r_{0}\geq\sigma_{1}/2+1>1$.

Suppose that $3\leq\mu_{0}\leq8$.
If $\sigma_{1}\geq2$, then $r_{0}\geq2/3$ by Fact~\ref{f221}.
Similarly, if $\sigma_{1}=1$ and $\sigma_{2}=2$ or $3$, then $r_{0}\geq1/2$.
It remains to consider $\sigma_{1}=\sigma_{2}=1$.
Then $1/2-1/(\mu_{1}+2)>1/\mu_{0}$ by Fact~\ref{f221}.
Thus, if $\sigma_{1}=\sigma_{2}=1$ and $\mu_{0}=3,4,5,6,7$ or $8$, then $\mu_{1}\geq 5,3,2,2,1$ or $1$, respectively.
Thus, if $\sigma_{1}=\sigma_{2}=1$ and $\mu_{0}=3,4,5,6,7$ or $8$, then $r_{0}\geq \sigma_{1}/2-1/(\mu_{1}+2)\geq5/14,3/10,1/4,1/4,1/6,1/6$, respectively.

Since $1/r_{0}=\mu_{0}-2/\sigma_{0}$, $\sigma_{0}\leq14,4,10,3,2,1,2$ or $1$, respectively.
\end{proof}

Since $1\leq\sigma_{2}\leq3$, we distinguish the following three subsubcases.

\noindent$\textbf{Subsubcase 2.2.1.}$ $\sigma_{2}=3$.

\begin{fact}\label{f222}
Suppose that $\sigma_{2}=3$. Then $\mu_{0}=\mu_{1}=1$ and $3\leq\sigma_{0}\leq6$.
If $\sigma_{0}=3,4,5$ or $6$, then $\sigma_{1}\leq4,2,1$ or $1$, respectively.
\end{fact}
\begin{proof}
Suppose that $\sigma_{2}=3$.
Then $\mu_{1}=1$ and $r_{0}\geq\sigma_{1}/2+1\geq3/2$ by Fact~\ref{f221}.
Since $r_{0}=\sigma_{0}/(\mu_{0}\sigma_{0}-2)$, $\mu_{0}\leq2/\sigma_{0}+2/3\leq8/3$.
Thus, $\mu_{0}\leq2$.
If $\mu_{0}=2$, then $\sigma_{0}\geq2$ for $\mu_{0}\sigma_{0}\geq3$, and hence $\mu_{0}\leq5/3$, a contradiction.
Thus, $\mu_{0}=1$.

Since $r_{0}\geq3/2$, $3\leq\sigma_{0}\leq6$.
Since $r_{0}\geq\sigma_{1}/2+1$, $\sigma_{1}\leq2\sigma_{0}/(\sigma_{0}-2)-2$.
Thus, if $\sigma_{0}=3,4,5$ or $6$, then $\sigma_{1}\leq4,2,1$ or $1$, respectively.
\end{proof}

If $\sigma_{0}=3$, then $1\leq\sigma_{1}\leq4$, $0\leq b_{0}\leq1$, $0\leq b_{1}\leq2$ and $0\leq b_{2}\leq1$ by Fact~\ref{f222}.
Thus, if $\sigma_{0}=3$, then  $G\subset D\bigl((3-2b_{0},4-2b_{1},3-2b_{2}), (b_{0},b_{1},b_{2}),(1,1);0,0\bigr)\subset\mathcal F_{j}(n)$ for some $51\leq j\leq62$.
If $\sigma_{0}=4$, then $1\leq\sigma_{1}\leq2$, $0\leq b_{0}\leq2$, $0\leq b_{1}\leq1$ and $0\leq b_{2}\leq1$ by Fact~\ref{f222}.
Thus, if $\sigma_{0}=4$, then  $G\subset D\bigl((4-2b_{0},2-2b_{1},3-2b_{2}),
(b_{0},b_{1},b_{2}),(1,1);0,0\bigr)
\subset\mathcal F_{j}(n)$ for some $63\leq j\leq74$.
If $\sigma_{0}=5$ or $6$, then $\sigma_{1}=1$,  $0\leq b_{0}\leq\lfloor\sigma_{0}/2\rfloor$, $b_{1}=0$ and $0\leq b_{2}\leq1$ by Fact~\ref{f222}.
Thus, if $\sigma_{0}=5$ or $6$, then $G\subset D\bigl((6-2b_{0},1,3-2b_{2}),
(b_{0},0,b_{2}),(1,1);0,0\bigr)
\subset\mathcal F_{j}(n)$ for some $75\leq j\leq82$.

\noindent$\textbf{Subsubcase 2.2.2.}$ $\sigma_{2}=2$.

\begin{fact}\label{f223}
Suppose that $\sigma_{2}=2$. Then $1\leq\mu_{0}\leq4$ and $\mu_{1}\geq1$.
If $\mu_{0}=1$, then $3\leq\sigma_{0}\leq6$ and $3\leq\sigma_{1}\leq\lfloor2\sigma_{0}/(\sigma_{0}-2)\rfloor$.
If $\mu_{0}=2$, then $(\sigma_{0},\sigma_{1})=(2,2)$.
If $\mu_{0}=3$, then $\sigma_{0}=1$ and $1\leq\sigma_{1}\leq2$, or $(\sigma_{0},\sigma_{1})=(2,1)$.
If $\mu_{0}=4$, then $(\sigma_{0},\sigma_{1})=(1,1)$.
\end{fact}
\begin{proof}
Suppose that $\sigma_{2}=2$.
Then $r_{0}\geq\sigma_{1}/2>1/\mu_{0}$ by Fact~\ref{f221}.

If $\mu_{0}=1$, then $\sigma_{1}\geq3$ and $r_{0}\geq3/2$.
Thus, if $\mu_{0}=1$, then  $3\leq\sigma_{0}\leq6$ for $\mu_{0}\sigma_{0}\geq3$ and $r_0=\frac{\sigma_0}{\mu_0\sigma_0-2}$.
If $\mu_{0}=1$ and $\sigma_{0}=3,4,5$ or $6$, by  Fact~\ref{f221}, then $\sigma_{1}\leq6,4,3$ or $3$, respectively.

If $\mu_{0}\geq2$, then $r_{0}\geq1/2$.
Thus, if $\mu_{0}\geq2$, then  $(\mu_{0}-2)\sigma_{0}\leq2$ for $r_0=\frac{\sigma_0}{\mu_0\sigma_0-2}$.
If $\mu_{0}=2$, then $\sigma_{1}\geq2$ and $r_{0}\geq1$.
Thus, if $\mu_{0}=2$, then $\sigma_{0}=2$ for $\mu_{0}\sigma_{0}\geq3$ and $r_0=\frac{\sigma_0}{\mu_0\sigma_0-2}$.
Thus, if $\mu_{0}=2$, then $r_0=1$ and $\sigma_{1}=2$.
If $\mu_{0}=3$, then $1\leq\sigma_{0}\leq2$  for $r_0=\frac{\sigma_0}{\mu_0\sigma_0-2}$.
If $\mu_{0}=3$ and $\sigma_{0}=1$, then $r_{0}=1$ and $\sigma_{1}\leq2$;
if $\mu_{0}=3$ and $\sigma_{0}=2$, then $r_{0}=1/2$ and $\sigma_{1}=1$.
If $\mu_{0}=4$, then $\sigma_{0}=1$ and $\sigma_{1}=1$ for $r_0=\frac{\sigma_0}{\mu_0\sigma_0-2}$ and $r_{0}\geq\sigma_{1}/2$.
If $\mu_{0}\geq5$, then
$(\mu_{0}-2)\sigma_{0}>2$, a contradiction.
By the definition of $D(\boldsymbol a,\boldsymbol b,\boldsymbol\mu;\kappa,\tau)$, $\sigma_{1}\geq1$ and $\mu_{1}\geq1$.
\end{proof}

If $\mu_{0}=1$ and $\sigma_{0}=3$, then $3\leq\sigma_{1}\leq6$, $0\leq b_{0}\leq1$, $0\leq b_{1}\leq\lfloor\sigma_{1}/2\rfloor$ and $0\leq b_{2}\leq1$ by Fact~\ref{f223}.
Thus, if $\mu_{0}=1$ and $\sigma_{0}=3$, then $G\subset D\bigl((3-2b_{0},6-2b_{1},2-2b_{2}), (b_{0},b_{1},b_{2}),(1,n);0,0\bigr)=\mathcal F_{j}(n)$ for some $1\leq j\leq16$.
If $\mu_{0}=1$ and $\sigma_{0}=4$, then $3\leq\sigma_{1}\leq4$, $0\leq b_{0}\leq2$, $0\leq b_{1}\leq\lfloor\sigma_{1}/2\rfloor$ and $0\leq b_{2}\leq1$ by Fact~\ref{f223}.
Thus, if $\mu_{0}=1$ and $\sigma_{0}=4$, then $G\subset D\bigl((4-2b_{0},4-2b_{1},2-2b_{2}),
(b_{0},b_{1},b_{2}),(1,n);0,0\bigr)
=\mathcal F_{j}(n)$ for some $17\leq j\leq34$.
If $\mu_{0}=1$ and $\sigma_{0}=5$ or $6$, then $\sigma_{1}=3$, $0\leq b_{0}\leq\lfloor\sigma_{0}/2\rfloor\leq3$, $0\leq b_{1}\leq1$ and $0\leq b_{2}\leq1$ by Fact~\ref{f223}.
Thus, if $\mu_{0}=1$ and $\sigma_{0}=5$ or $6$, then $G\subset D\bigl((6-2b_{0},3-2b_{1},2-2b_{2}),
(b_{0},b_{1},b_{2}),(1,n);0,0\bigr)
=\mathcal F_{j}(n)$ for some $35\leq j\leq50$.

If $\mu_{0}=2$, then $\sigma_{0}=\sigma_{1}=2$ and $0\leq b_{0},b_{1},b_{2}\leq1$ by Fact~\ref{f223}.
Thus, if $\mu_{0}=2$, then $G\subset D\bigl((2-2b_{0},2-2b_{1},2-2b_{2}),
(b_{0},b_{1},b_{2}),(2,n);0,0\bigr)
\subset\mathcal F_{j}(n)$ for some $83\leq j\leq86$.
If $\mu_{0}=3$ and $\sigma_{0}=1$, then $1\leq\sigma_{1}\leq2$, $b_{0}=0$, $0\leq b_{1}\leq\lfloor\sigma_{1}/2\rfloor$ and $0\leq b_{2}\leq1$ by Fact~\ref{f223}.
Thus, if $\mu_{0}=3$ and $\sigma_{0}=1$, then  $G\subset D\bigl((1,2-2b_{1},2-2b_{2}),(0,b_{1},b_{2}),(3,n);0,0\bigr)\subset\mathcal F_{j}(n)$ for some $83\leq j\leq86$.
If $\mu_{0}=3$ and $\sigma_{0}=2$, then $\sigma_{1}=1$, $b_{1}=0$ and $0\leq b_{0},b_{2}\leq1$ by Fact~\ref{f223}.
Thus, if $\mu_{0}=3$ and $\sigma_{0}=2$, then $G\subset D\bigl((2-2b_{0},1,2-2b_{2}), (b_{0},0,b_{2}),(3,n);0,0\bigr)\subset\mathcal F_{j}(n)$ for some $87\leq j\leq88$.
If $\mu_{0}=4$, then $\sigma_{0}=\sigma_{1}=1$, $b_{0}=b_{1}=0$ and $0\leq b_{2}\leq1$ by Fact~\ref{f223}.
Thus, if $\mu_{0}=4$, then $G\subset D\bigl((1,1,2-2b_{2}), (0,0,b_{2}),(4,n);0,0\bigr)\subset\mathcal F_{j}(n)$ for some $87\leq j\leq88$.
There is no   $\mu_{0}\geq 5$ by Fact~\ref{f223}.

\noindent$\textbf{Subsubcase 2.2.3.}$ $\sigma_{2}=1$.

Then  $a_{2}=1$ and $b_{2}=0$.
Let $G'$ be the graph obtained from $G$ by adding a vertex to $C_{2}$ that is adjacent to every vertex of $C$ and has no neighbors in $Z$.
Then $\sigma_{2}=2$ for $G'$ and $G\subset G'$.
If $r_{0}\geq\sigma_{1}/2$, then $r_{1}\geq0$ by Fact~\ref{f221} and  $r_1=\frac{r_0-\sigma_1/2}{1+\mu_1(r_0-\sigma_1/2)}$.
Thus, if $r_{0}\geq\sigma_{1}/2$, then $G'$ belongs to $\mathcal{D}_{2}$ and satisfies the conditions in Subsubcase~2.2.2.
Thus, if $r_{0}\geq\sigma_{1}/2$, then $G\subset G'\subset H$ for some graph $H$ considered in Subsubcase~2.2.2.
Let $r_{0}<\sigma_{1}/2$.

Then  $1\leq\mu_{1}\leq\lfloor1/(\sigma_{1}/2-r_{0})-2\rfloor$, $\sigma_{1}/2>1/\mu_{0}$ and $\mu_{1}>1/(\sigma_{1}/2-1/\mu_{0})-2$ by Fact~\ref{f221}.

Then $1\leq\mu_{0}\leq8$ and $1\leq\sigma_{1}\leq6$, and the corresponding upper bounds for $\sigma_{0}$ are given by Fact~\ref{f22bounds}.
Then $2r_{0}<\sigma_{1}\leq 2r_{0}+2/3$ by Fact~\ref{f221}.

Suppose that $\mu_{0}=1$.
Then $r_{0}=\sigma_{0}/(\sigma_{0}-2)$ and $\sigma_{1}\geq3$ by $\sigma_{1}/2>1/\mu_{0}$.
Then $\sigma_{0}\geq4$ by $r_{0}<\sigma_{1}/2\leq3$.
If $\sigma_{0}=4$, then $r_{0}=2$ and $\sigma_{1}\geq5$ by $r_{0}<\sigma_{1}/2$.
Thus, if $\sigma_{0}=4$, then $r_{0}\geq\sigma_{1}/2-1/3\geq13/6>2$ by Fact~\ref{f221}, a contradiction.
Thus, $\sigma_{0}\geq5$ gives $r_{0}\leq5/3$.
Thus, $\sigma_{1}=3$ or $4$ by $r_{0}<\sigma_{1}/2\leq r_{0}+1/3$.
If $\sigma_{1}=4$, then $r_{0}=5/3$ by $r_{0}\geq\sigma_{1}/2-1/3$.
Thus, if $\sigma_{1}=4$, then $\sigma_{0}=5$ and $\mu_{1}=1$ by $r_{0}=\sigma_{0}/(\sigma_{0}-2)$ and $\mu_{1}\leq\lfloor1/(\sigma_{1}/2-r_{0})-2\rfloor$.
If $\sigma_{1}=3$, then $7\leq\sigma_{0}\leq14$ by $r_{0}<3/2$ and Fact~\ref{f22bounds}.
Thus, if $\sigma_{1}=3$, then $\mu_{1}\leq\lfloor8/(\sigma_{0}-6)\rfloor$ by $r_{0}=\sigma_{0}/(\sigma_{0}-2)$ and $\mu_{1}\leq\lfloor1/(\sigma_{1}/2-r_{0})-2\rfloor$.
Thus, if $\sigma_{0}=5$, $\sigma_{1}=4$ and $\mu_{1}=1$, then $G\subset D\bigl((5-2b_{0},4-2b_{1},1),\allowbreak(b_{0},b_{1},0),\allowbreak(1,1);0,0\bigr)=A_{i}$ for some $18\leq i\leq26$ (see Figure~\ref{fig:A18-A26}).
If $\sigma_{0}=7$, $\sigma_{1}=3$ and $1\leq\mu_{1}\leq8$, then $G\subset D\bigl((7-2b_{0},3-2b_{1},1),\allowbreak(b_{0},b_{1},0),\allowbreak(1,8);0,0\bigr)=A_{i}$ for some $31\leq i\leq38$.
If $\sigma_{0}=8$, $\sigma_{1}=3$ and $1\leq\mu_{1}\leq4$, then $G\subset D\bigl((8-2b_{0},3-2b_{1},1),\allowbreak(b_{0},b_{1},0),\allowbreak(1,4);0,0\bigr)=A_{i}$ for some $39\leq i\leq48$.
If $9\leq\sigma_{0}\leq10$, $\sigma_{1}=3$ and $1\leq\mu_{1}\leq2$, then $G\subset D\bigl((10-2b_{0},3-2b_{1},1),\allowbreak(b_{0},b_{1},0),\allowbreak(1,2);0,0\bigr)=A_{i}$ for some $55\leq i\leq66$.
If $11\leq\sigma_{0}\leq14$, $\sigma_{1}=3$ and $\mu_{1}=1$, then $G\subset D\bigl((14-2b_{0},3-2b_{1},1),\allowbreak(b_{0},b_{1},0),\allowbreak(1,1);0,0\bigr)=A_{i}$ for some $67\leq i\leq82$.

% Begin embedded drawing: figures/A18_A26.tex
% Requires tikz and graphicx. The filled vertices form one clique;
% each line from an X-vertex to a C_j box denotes complete adjacency.
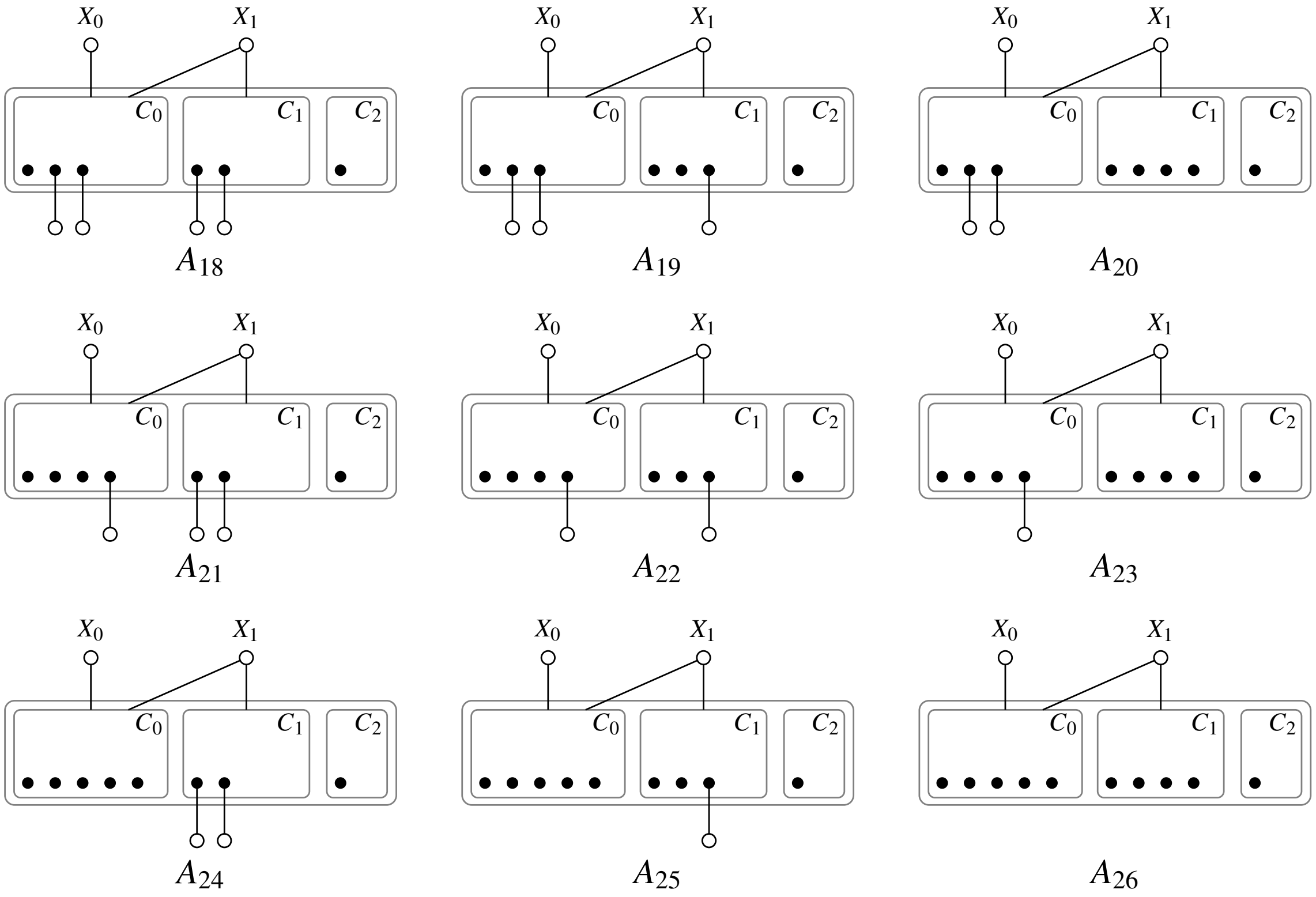
\begin{figure}[!htbp]
\centering
\resizebox{\textwidth}{!}{%
\begin{tikzpicture}[x=1cm,y=1cm,line width=0.5pt,
  every node/.style={font=\small}]
% #1: x-offset of the class; #2: a_j; #3: b_j.
\newcommand{\Aclass}[3]{%
  \ifnum#2>0
    \foreach \i in {1,...,#2}{%
      \pgfmathsetmacro{\xx}{#1+0.15+0.30*(\i-1)}
      \fill (\xx,0) circle (1.8pt);}%
  \fi
  \ifnum#3>0
    \foreach \i in {1,...,#3}{%
      \pgfmathsetmacro{\xx}{#1+0.15+0.30*(#2+\i-1)}
      \draw (\xx,0)--(\xx,-0.63);
      \fill (\xx,0) circle (1.8pt);
      \draw[fill=white] (\xx,-0.63) circle (2.1pt);}%
  \fi}
% Arguments: label, a_0, b_0, a_1, b_1.  Here a_2=1, b_2=0,
% and |X_0|=|X_1|=1 in all nine graphs.
\newcommand{\Aexample}[5]{%
  \draw[rounded corners=3pt,gray] (-0.1,-0.24) rectangle (4.18,0.90);
  \draw[rounded corners=2pt,gray] (0,-0.16) rectangle (1.68,0.80);
  \draw[rounded corners=2pt,gray] (1.85,-0.16) rectangle (3.23,0.80);
  \draw[rounded corners=2pt,gray] (3.42,-0.16) rectangle (4.08,0.80);
  \node[anchor=north east,inner sep=1pt,font=\scriptsize] at (1.68,0.80) {$C_0$};
  \node[anchor=north east,inner sep=1pt,font=\scriptsize] at (3.23,0.80) {$C_1$};
  \node[anchor=north east,inner sep=1pt,font=\scriptsize] at (4.08,0.80) {$C_2$};
  \draw (0.84,0.80)--(0.84,1.37);
  \draw (2.54,1.37)--(2.54,0.80);
  \draw (2.54,1.37)--(1.25,0.80);
  \draw[fill=white] (0.84,1.37) circle (2.1pt);
  \draw[fill=white] (2.54,1.37) circle (2.1pt);
  \node[above=1pt,font=\scriptsize] at (0.84,1.37) {$X_0$};
  \node[above=1pt,font=\scriptsize] at (2.54,1.37) {$X_1$};
  \Aclass{0}{#2}{#3}\Aclass{1.85}{#4}{#5}\Aclass{3.42}{1}{0}
  \node at (2.04,-1.0) {$A_{#1}$};}
\begin{scope}[shift={(0,0)}]      \Aexample{18}{1}{2}{0}{2} \end{scope}
\begin{scope}[shift={(5.0,0)}]   \Aexample{19}{1}{2}{2}{1} \end{scope}
\begin{scope}[shift={(10.0,0)}]    \Aexample{20}{1}{2}{4}{0} \end{scope}
\begin{scope}[shift={(0,-3.35)}]    \Aexample{21}{3}{1}{0}{2} \end{scope}
\begin{scope}[shift={(5.0,-3.35)}] \Aexample{22}{3}{1}{2}{1} \end{scope}
\begin{scope}[shift={(10.0,-3.35)}]  \Aexample{23}{3}{1}{4}{0} \end{scope}
\begin{scope}[shift={(0,-6.70)}]    \Aexample{24}{5}{0}{0}{2} \end{scope}
\begin{scope}[shift={(5.0,-6.70)}] \Aexample{25}{5}{0}{2}{1} \end{scope}
\begin{scope}[shift={(10.0,-6.70)}]  \Aexample{26}{5}{0}{4}{0} \end{scope}
\end{tikzpicture}%
}
\caption{The graphs $A_{18},\ldots,A_{26}$.}
\label{fig:A18-A26}
\end{figure}
% End embedded drawing: figures/A18_A26.tex

Suppose that $\mu_{0}=2$.
Then $2\leq\sigma_{0}\leq4$ and $\sigma_{1}\geq2$ by Fact~\ref{f22bounds}, $\mu_{0}\sigma_{0}\geq3$  and $\sigma_{1}/2>1/\mu_{0}$.
If $\sigma_{0}=2$, then $r_{0}=1$ and $\sigma_{1}\geq3$ by $r_{0}<\sigma_{1}/2$.
Thus, if $\sigma_{0}=2$, then $r_{0}\geq\sigma_{1}/2-1/3\geq7/6>1$ by Fact~\ref{f221}, a contradiction.
Thus, $3\leq\sigma_{0}\leq4$ and $r_{0}\leq3/4$ by $r_{0}=\sigma_{0}/(2\sigma_{0}-2)$.
Thus, $\sigma_{1}=2$ by $r_{0}\leq3/4$, $r_{0}<\sigma_{1}/2\leq r_{0}+1/3$ and $\sigma_{1}\geq2$.
Thus, $\mu_{1}\leq\lfloor2/(\sigma_{0}-2)\rfloor$ by $r_{0}=\sigma_{0}/(2\sigma_{0}-2)$ and  $\mu_{1}\leq\lfloor1/(\sigma_{1}/2-r_{0})-2\rfloor$.
If $\sigma_{0}=3$, $\sigma_{1}=2$ and $1\leq\mu_{1}\leq2$, then $G\subset D\bigl((3-2b_{0},2-2b_{1},1),\allowbreak(b_{0},b_{1},0),\allowbreak(2,2);0,0\bigr)=A_{i}$ for some $5\leq i\leq8$.
If $\sigma_{0}=4$, $\sigma_{1}=2$ and $\mu_{1}=1$, then $G\subset D\bigl((4-2b_{0},2-2b_{1},1),\allowbreak(b_{0},b_{1},0),\allowbreak(2,1);0,0\bigr)=A_{i}$ for some $12\leq i\leq17$.

Suppose that $\mu_{0}=3$.
Then $1\leq\sigma_{0}\leq10$ by Fact~\ref{f22bounds}.
If $\sigma_{0}=1$, then $r_{0}=1$ and $\sigma_{1}\geq3$ by $r_{0}<\sigma_{1}/2$.
Thus, if $\sigma_{0}=1$, then $r_{0}\geq\sigma_{1}/2-1/3\geq7/6>1$ by Fact~\ref{f221}, a contradiction.
If $\sigma_{0}=2$, then $r_{0}=1/2$ and $\sigma_{1}\geq2$ by $r_{0}<\sigma_{1}/2$.
Thus, if $\sigma_{0}=2$, then $r_{0}\geq\sigma_{1}/2-1/3\geq2/3>1/2$ by Fact~\ref{f221}, a contradiction.
Thus, $3\leq\sigma_{0}\leq10$, $r_{0}\leq3/7$ and $\sigma_{1}=1$ by $r_{0}=\sigma_{0}/(3\sigma_{0}-2)$, $\sigma_{1}/2\leq r_{0}+1/3$ and $\sigma_{1}\geq1$.
Thus, $\mu_{1}\geq5$ by $\sigma_{1}/2-1/(\mu_{1}+2)>1/\mu_{0}$.
Thus, $\mu_{1}\leq\lfloor4+8/(\sigma_{0}-2)\rfloor$ by $r_{0}=\sigma_{0}/(3\sigma_{0}-2)$ and $\mu_{1}\leq\lfloor1/(\sigma_{1}/2-r_{0})-2\rfloor$.
If $\sigma_{0}=3$, $\sigma_{1}=1$ and $5\leq\mu_{1}\leq12$, then $G\subset D\bigl((3-2b_{0},1,1),\allowbreak(b_{0},0,0),\allowbreak(3,12);0,0\bigr)=A_{i}$ for some $1\leq i\leq2$.
If $\sigma_{0}=4$, $\sigma_{1}=1$ and $5\leq\mu_{1}\leq8$, then $G\subset D\bigl((4-2b_{0},1,1),\allowbreak(b_{0},0,0),\allowbreak(3,8);0,0\bigr)=A_{i}$ for some $9\leq i\leq11$.
If $5\leq\sigma_{0}\leq6$, $\sigma_{1}=1$ and $5\leq\mu_{1}\leq6$, then $G\subset D\bigl((6-2b_{0},1,1),\allowbreak(b_{0},0,0),\allowbreak(3,6);0,0\bigr)=A_{i}$ for some $27\leq i\leq30$.
If $7\leq\sigma_{0}\leq10$, $\sigma_{1}=1$ and $\mu_{1}=5$, then $G\subset D\bigl((10-2b_{0},1,1),\allowbreak(b_{0},0,0),\allowbreak(3,5);0,0\bigr)=A_{i}$ for some $49\leq i\leq54$.

Suppose that $\mu_{0}=4$.
Then $1\leq\sigma_{0}\leq3$ by Fact~\ref{f22bounds}.
If $\sigma_{0}=1$, then $r_{0}=1/2$ and $\sigma_{1}\geq2$ by $r_{0}<\sigma_{1}/2$.
Thus, if $\sigma_{0}=1$, then $r_{0}\geq\sigma_{1}/2-1/3\geq2/3>1/2$ by Fact~\ref{f221}, a contradiction.
Thus, $2\leq\sigma_{0}\leq3$, $r_{0}\leq1/3$ and $\sigma_{1}=1$ by $r_{0}=\sigma_{0}/(4\sigma_{0}-2)$, $\sigma_{1}/2\leq r_{0}+1/3$ and $\sigma_{1}\geq1$.
Thus, $\mu_{1}\geq3$ by $\sigma_{1}/2-1/(\mu_{1}+2)>1/\mu_{0}$.
Thus, $\mu_{1}\leq\lfloor2+2/(\sigma_{0}-1)\rfloor$ by $r_{0}=\sigma_{0}/(4\sigma_{0}-2)$ and $\mu_{1}\leq\lfloor1/(\sigma_{1}/2-r_{0})-2\rfloor$.
If $\sigma_{0}=2$, $\sigma_{1}=1$ and $3\leq\mu_{1}\leq4$, then $G\subset D\bigl((2-2b_{0},1,1),\allowbreak(b_{0},0,0),\allowbreak(4,4);0,0\bigr)\subset\mathcal F_{89}(n)$.
If $\sigma_{0}=3$, $\sigma_{1}=1$ and $\mu_{1}=3$, then $G\subset D\bigl((3-2b_{0},1,1),\allowbreak(b_{0},0,0),\allowbreak(4,3);0,0\bigr)=A_{i}$ for some $3\leq i\leq4$.

Suppose that $\mu_{0}=5$.
Then $1\leq\sigma_{0}\leq2$ and $r_{0}\leq1/3$ by Fact~\ref{f22bounds} and $r_{0}=\sigma_{0}/(5\sigma_{0}-2)$.
Thus, $\sigma_{1}=1$ by $\sigma_{1}/2\leq r_{0}+1/3$ and $\sigma_{1}\geq1$.
Thus, $\mu_{1}\geq2$ by $\sigma_{1}/2-1/(\mu_{1}+2)>1/\mu_{0}$.
If $\sigma_{0}=1$ or $2$, then $r_{0}=1/3$ or $1/4$, respectively.
Thus, if $\sigma_{0}=1$ or $2$, then $\mu_{1}\leq4$ or $2$, respectively, by $\mu_{1}\leq\lfloor1/(\sigma_{1}/2-r_{0})-2\rfloor$.
If $\sigma_{0}=1$, $\sigma_{1}=1$ and $2\leq\mu_{1}\leq4$, then $G\subset D\bigl((1,1,1),\allowbreak(0,0,0),\allowbreak(5,4);0,0\bigr)\subset\mathcal F_{89}(n)$.
If $\sigma_{0}=2$, $\sigma_{1}=1$ and $\mu_{1}=2$, then $G\subset D\bigl((2-2b_{0},1,1),\allowbreak(b_{0},0,0),\allowbreak(5,2);0,0\bigr)\subset\mathcal F_{90}(n)$.

Suppose that $\mu_{0}=6$.
Then $\sigma_{0}=1$ and $r_{0}=1/4$ by Fact~\ref{f22bounds} and $r_{0}=\sigma_{0}/(6\sigma_{0}-2)$.
Thus, $\sigma_{1}=1$ by $\sigma_{1}/2\leq r_{0}+1/3$ and $\sigma_{1}\geq1$.
Thus, $\mu_{1}\geq2$ by $\sigma_{1}/2-1/(\mu_{1}+2)>1/\mu_{0}$.
Thus, $\mu_{1}=2$ by $\mu_{1}\geq2$ and $\mu_{1}\leq\lfloor1/(\sigma_{1}/2-r_{0})-2\rfloor$.
If $\sigma_{0}=1$, $\sigma_{1}=1$ and $\mu_{1}=2$, then $G\subset D\bigl((1,1,1),\allowbreak(0,0,0),\allowbreak(6,2);0,0\bigr)\subset\mathcal F_{90}(n)$.

Suppose that $\mu_{0}=7$.
Then $1\leq\sigma_{0}\leq2$ by Fact~\ref{f22bounds}.
If $\sigma_{0}=1$ or $2$, then $r_{0}=1/5$ or $1/6$, respectively.
Thus, $r_{0}\leq1/5$ and $\sigma_{1}=1$ by $\sigma_{1}/2\leq r_{0}+1/3$ and $\sigma_{1}\geq1$.
Thus, $\mu_{1}=1$ by $r_{0}\leq1/5$ and $\mu_{1}\leq\lfloor1/(\sigma_{1}/2-r_{0})-2\rfloor$.
If $1\leq\sigma_{0}\leq2$, $\sigma_{1}=1$ and $\mu_{1}=1$, then $G\subset D\bigl((2-2b_{0},1,1),\allowbreak(b_{0},0,0),\allowbreak(7,1);0,0\bigr)\subset\mathcal F_{91}(n)$.

Suppose that $\mu_{0}=8$.
Then $\sigma_{0}=1$ and $r_{0}=1/6$ by Fact~\ref{f22bounds} and $r_{0}=\sigma_{0}/(8\sigma_{0}-2)$.
Thus, $\sigma_{1}=1$ by $\sigma_{1}/2\leq r_{0}+1/3$ and $\sigma_{1}\geq1$.
Thus, $\mu_{1}=1$ by $\mu_{1}\geq1$ and $\mu_{1}\leq\lfloor1/(\sigma_{1}/2-r_{0})-2\rfloor$.
If $\sigma_{0}=1$, $\sigma_{1}=1$ and $\mu_{1}=1$, then $G\subset D\bigl((1,1,1),\allowbreak(0,0,0),\allowbreak(8,1);0,0\bigr)\subset\mathcal F_{91}(n)$.

\noindent$\textbf{Subcase 2.3.}$ $h=3$.

\begin{fact}\label{f23bounds}
Then $1\leq\sigma_{2},\sigma_{3}\leq2$, $(\sigma_{2},\sigma_{3})\neq(2,2)$, $1\leq\sigma_{1}\leq5$, $1\leq\mu_{1}\leq5$ and $1\leq\mu_{0}\leq3$.
If $(\sigma_{2},\sigma_{3})=(1,2)$, then $\mu_{1}=1$ and $r_{0}\geq\sigma_{1}/2+1$.
If $(\sigma_{2},\sigma_{3})=(2,1)$, then $\mu_{1}=\mu_{2}=1$ and $r_{0}\geq\sigma_{1}/2+2$.
If $(\sigma_{2},\sigma_{3})=(1,1)$, then $r_{1}\geq\mu_{2}/[2(\mu_{2}+2)]$.
If $\mu_{0}=1,2$ or $3$, then $\sigma_{0}\leq38,3$ or $1$, respectively.
\end{fact}
\begin{proof}
Then $r_{2}\geq\sigma_{3}/2-1$ and $r_{j}<1/\mu_{j}$ for $j=1,2$.
If $r_{0}-\sigma_{1}/2\leq0$, then $r_{1}=\frac{r_{0}-\sigma_{1}/2}{1+\mu_{1}(r_{0}-\sigma_{1}/2)}\leq0$.
If $r_{0}-\sigma_{1}/2>0$, then  $r_{1}\leq5/7$ for $r_{0}\leq3$.
If $r_{1}-\sigma_{2}/2\leq0$, then $r_{2}=\frac{r_{1}-\sigma_{2}/2}{1+\mu_{2}(r_{1}-\sigma_{2}/2)}\leq0$.
If $r_{1}-\sigma_{2}/2>0$, then $r_{2}\leq3/17$ for $r_{1}\leq5/7$.
Since $r_{2}\geq\sigma_{3}/2-1$, $1\leq\sigma_{3}\leq2$.
Since $r_{2}\geq-1/2$,  $r_{1}\geq\sigma_{2}/2-1/(\mu_{2}+2)\geq1/6$.
Since $r_{1}\leq5/7$, $\sigma_{2}/2-1/3\leq5/7$, which implies that $1\leq\sigma_{2}\leq2$.
Since  $1/6\leq r_{1}<1/\mu_{1}$, $\mu_{1}\leq5$.
Since $r_{1}=(r_{0}-\sigma_{1}/2)/(1+\mu_{1}(r_{0}-\sigma_{1}/2))$, $r_{0}=\sigma_{1}/2+r_{1}/(1-\mu_{1}r_{1})\geq\sigma_{1}/2+1/5\geq7/10$.
Then $\sigma_{1}/2+1/5\leq r_{0}\leq3$, and hence $1\leq\sigma_{1}\leq5$.
Also, $\mu_{0}=2/\sigma_{0}+1/r_{0}\leq2+10/7<4$, and hence $1\leq\mu_{0}\leq3$.

Suppose that $(\sigma_{2},\sigma_{3})=(1,1)$.
Then $r_{1}\geq \sigma_{2}/2-1/(\mu_{2}+2)=\mu_{2}/[2(\mu_{2}+2)]$.
Since $r_{1}<1/\mu_{1}$, $(\mu_{1}-2)\mu_{2}<4$.
Since $r_{0}=\sigma_{1}/2+r_{1}/(1-\mu_{1}r_{1})$, $r_{0}\geq\sigma_{1}/2+\mu_{2}/[4+(2-\mu_{1})\mu_{2}]$.

Suppose that $(\sigma_{2},\sigma_{3})=(1,2)$.
Then $r_{2}\geq0$, which implies that $r_{1}\geq1/2$.
Since $r_{1}<1/\mu_{1}$, $\mu_{1}=1$.
Thus, $r_{0}\geq\sigma_{1}/2+(1/2)/(1-1/2)=\sigma_{1}/2+1$.

Suppose that $(\sigma_{2},\sigma_{3})=(2,1)$.
Then $r_{1}\geq\sigma_{2}/2-1/(\mu_{2}+2)\geq2/3$.
Since $r_{1}<1/\mu_{1}$, $\mu_{1}=1$.
Since $r_{1}\leq5/7$, $\sigma_{2}/2-1/(\mu_{2}+2)\leq5/7$, and hence $\mu_{2}=1$.
Then $r_{0}\geq\sigma_{1}/2+(2/3)/(1-2/3)=\sigma_{1}/2+2$.

If $(\sigma_{2},\sigma_{3})=(2,2)$, then $r_{2}\geq0$ gives $r_{1}\geq1$, a contradiction to $r_{1}\leq5/7$.

If $\mu_{0}=2$, then $\sigma_{0}\leq3$ for $r_{0}\geq7/10$ and $\mu_{0}-2/\sigma_{0}=1/r_{0}$.
If $\mu_{0}=3$, then $\sigma_{0}\leq1$  for $r_{0}\geq7/10$ and $\mu_{0}-2/\sigma_{0}=1/r_{0}$.
Suppose that $\mu_{0}=1$.
If $\sigma_{2}=2$ or $\sigma_{3}=2$, then $r_{0}\geq\sigma_{1}/2+1\geq3/2$.
If $\sigma_{1}\geq2$, then $r_{0}\geq\sigma_{1}/2+1/5\geq6/5$.
Let $\sigma_{1}=\sigma_{2}=\sigma_{3}=1$.

The parameters  satisfy the conditions for $\mathcal{D}_{3}$ with $r_{0}=1/\mu_{0}=1$ if and only if $1/(\mu_{1}+2)\geq\mu_{2}/[2(\mu_{2}+2)]$.
Thus, $\mu_{1}\mu_{2}>4$.
Note that $(\mu_{1}-2)\mu_{2}<4$ and $1\leq\mu_{1}\leq5$.
If $\mu_{1}=1$, then $\mu_{2}\geq5$ and $r_0\geq\frac{\sigma_1}{2}+\frac{\mu_2}{4+(2-\mu_1)\mu_2}\geq1/2+5/9=19/18$.
If $\mu_{1}=2$, then $\mu_{2}\geq3$ and $r_{0}\geq1/2+3/4=5/4$.
If $\mu_{1}=3$, then $\mu_{2}\geq2$ and $r_{0}\geq1/2+2/2=3/2$.
If $\mu_{1}=4$, then $\mu_{2}=1$, which contradicts $\mu_{1}\mu_{2}>4$.
If $\mu_{1}=5$, then $\mu_{2}=1$ and $r_{0}\geq1/2+1=3/2$.
Thus, $r_{0}\geq19/18$ in each case.
Since $r_{0}=\sigma_{0}/(\sigma_{0}-2)$, $\sigma_{0}/(\sigma_{0}-2)\geq19/18$, which implies that $\sigma_{0}\leq38$.
\end{proof}

By Fact~\ref{f23bounds}, we distinguish the following three subsubcases.

\noindent$\textbf{Subsubcase 2.3.1.}$ $(\sigma_{2},\sigma_{3})=(1,2)$.

Then $\mu_{1}=1$ and $r_{0}\geq3/2$ by Fact~\ref{f23bounds}.
Then $\mu_{0}=2/\sigma_{0}+1/r_{0}\leq2/\sigma_{0}+2/3$.
If $\sigma_{0}=1$, then $\mu_{0}\geq3$ by $\mu_{0}\sigma_{0}>2$, which contradicts $\mu_{0}\leq8/3$.
Thus, $\sigma_{0}\geq2$ and $\mu_{0}\leq5/3$, which implies that $\mu_{0}=1$.
Then $3\leq\sigma_{0}\leq6$ by $r_{0}=\sigma_{0}/(\sigma_{0}-2)\geq3/2$ and $\mu_{0}\sigma_{0}>2$.
Then $\sigma_{1}\leq\lfloor2\sigma_{0}/(\sigma_{0}-2)-2\rfloor$ by $r_{0}\geq\sigma_{1}/2+1$.
Thus, if $\sigma_{0}=3,4,5$ or $6$, then $\sigma_{1}\leq4,2,1$ or $1$, respectively.
Let $n\geq\mu_{2}$.
If $\sigma_{0}=3$, then $G\subset D\bigl((3-2b_{0},4-2b_{1},1,2-2b_{3}),(b_{0},b_{1},0,b_{3}),(1,1,n);0,0\bigr)=\mathcal F_{j}(n)$ for some $51\leq j\leq62$.
If $\sigma_{0}=4$, then $G\subset D\bigl((4-2b_{0},2-2b_{1},1,2-2b_{3}),(b_{0},b_{1},0,b_{3}),(1,1,n);0,0\bigr)=\mathcal F_{j}(n)$ for some $63\leq j\leq74$.
If $\sigma_{0}=5$ or $6$, then $G\subset D\bigl((6-2b_{0},1,1,2-2b_{3}),(b_{0},0,0,b_{3}),(1,1,n);0,0\bigr)=\mathcal F_{j}(n)$ for some $75\leq j\leq82$.

\noindent$\textbf{Subsubcase 2.3.2.}$ $(\sigma_{2},\sigma_{3})=(2,1)$.

Then $\mu_{1}=\mu_{2}=1$ and $r_{0}\geq\sigma_{1}/2+2\geq5/2$ by Fact~\ref{f23bounds}.
Then $\mu_{0}=2/\sigma_{0}+1/r_{0}\leq2/\sigma_{0}+2/5$.
If $\sigma_{0}=1$, then $\mu_{0}\geq3$ by $\mu_{0}\sigma_{0}>2$, which contradicts $\mu_{0}\leq12/5$.
Thus, $\sigma_{0}\geq2$ and $\mu_{0}\leq7/5$, which implies that $\mu_{0}=1$.
Then $3\leq\sigma_{0}\leq10/3$ by $\mu_{0}\sigma_{0}>2$ and $r_{0}=\sigma_{0}/(\sigma_{0}-2)\geq5/2$.
Thus, $\sigma_{0}=3$ and $r_{0}=3$.
Thus, $1\leq\sigma_{1}\leq2$ by $r_{0}\geq\sigma_{1}/2+2$.
Thus, $G\subset D\bigl((3-2b_{0},2-2b_{1},2-2b_{2},1),(b_{0},b_{1},b_{2},0),(1,1,1);0,0\bigr)=A_{i}$ for some $91\leq i\leq98$.

\noindent$\textbf{Subsubcase 2.3.3.}$ $(\sigma_{2},\sigma_{3})=(1,1)$.

Then $a_{2}=a_{3}=1$ and $b_{2}=b_{3}=0$ by $\sigma_{j}=a_{j}+2b_{j}$ for $j=2,3$.
Let $G'$ be the graph obtained from $G$ by adding a vertex to $C_{3}$ that is adjacent to every vertex of $C$ and has no neighbors in $Z$.
Then $\sigma_{3}=2$ in $G'$ and $G\subset G'$.
If $r_{1}\geq1/2$, then $G'\in \mathcal{D}_{2}$ by $r_{2}\geq0$.
Thus, if $r_{1}\geq1/2$, then $G\subset G'\subset H$ for some graph $H$ considered in Subsubcase~2.3.1.
Let $r_{1}<1/2$.
Note that $r_{0}=\sigma_{0}/(\mu_{0}\sigma_{0}-2)$ and $r_{1}=(r_{0}-\sigma_{1}/2)/(1+\mu_{1}(r_{0}-\sigma_{1}/2))$.
By Fact~\ref{f23bounds}, $1\leq\mu_{2}\leq\lfloor2/(1-2r_{1})-2\rfloor$.

If $\mu_{0}=1$ and $\sigma_{1}=1$, then the parameters satisfy the conditions for $\mathcal{D}_{3}$ with $r_{0}=1/\mu_{0}=1$ if and only if $1/(\mu_{1}+2)\geq\mu_{2}/[2(\mu_{2}+2)]$.
Thus, if $\mu_{0}=1$ and $\sigma_{1}=1$, then $\mu_{1}\mu_{2}>4$.
By Fact~\ref{f23bounds} and $r_{0}=\sigma_{1}/2+r_{1}/(1-\mu_{1}r_{1})$, $r_{0}\geq\sigma_{1}/2+1/(6-\mu_{1})$.
Thus, $1\leq\sigma_{1}\leq\lfloor2r_{0}-2/(6-\mu_{1})\rfloor$.

Suppose that $\mu_{0}=1$.
Then $r_{0}=\sigma_{0}/(\sigma_{0}-2)$, $3\leq\sigma_{0}\leq38$, $1\leq\sigma_{1}\leq5$ and $1\leq\mu_{1}\leq5$ by Fact~\ref{f23bounds} and $\mu_{0}\sigma_{0}>2$.
By Fact~\ref{f23bounds} and $r_{1}<1/2$, $1/6\leq r_{1}<1/2$ and $1\leq\mu_{2}\leq\lfloor2/(1-2r_{1})-2\rfloor$.
If $\sigma_{1}=1$, then $b_{1}=0$ and $r_{1}=(\sigma_{0}+2)/[2(\sigma_{0}-2)+\mu_{1}(\sigma_{0}+2)]$.
Thus, if $\sigma_{1}=1$, then $\mu_{1}(\sigma_{0}+2)>8$ and $1\leq\mu_{2}\leq\lfloor4(\sigma_{0}+2)/[\mu_{1}(\sigma_{0}+2)-8]\rfloor$ by $r_{1}<1/2$ and $\mu_{2}\leq\lfloor2/(1-2r_{1})-2\rfloor$.
If $\sigma_{1}=1$, then $\mu_{1}\mu_{2}>4$, for the parameters do not satisfy the conditions for $\mathcal{D}_{3}$ with $r_{0}=1$.
If $\sigma_{1}=1$, then  $\mu_{1}(\sigma_{0}+2)-8=(\mu_{1}-2)(\sigma_{0}+2)+2(\sigma_{0}-2)$.
Thus, if $\sigma_{1}=1$, then $(\mu_{1}-2)\mu_{2}<4$ by $\sigma_{0}>2$ and $\mu_{2}[\mu_{1}(\sigma_{0}+2)-8]\leq4(\sigma_{0}+2)$.
If $\sigma_{1}=1$ and $\mu_{1}=3$, then $2\leq\mu_{2}\leq\lfloor4(\sigma_{0}+2)/(3\sigma_{0}-2)\rfloor\leq\lfloor20/7\rfloor=2$ by $\mu_{1}\mu_{2}>4$ and $\sigma_{0}\geq3$.
If $\sigma_{1}=1$ and $\mu_{1}=4$, then $\mu_{2}=1$ by $(\mu_{1}-2)\mu_{2}<4$, a contradiction to $\mu_{1}\mu_{2}>4$.
If $\sigma_{1}=1$ and $\mu_{1}=5$, then $\mu_{2}=1$ by $(\mu_{1}-2)\mu_{2}<4$.
Thus, if $\sigma_{1}=1$, then $\mu_{1}=1,2,3$ or $5$.
If $\sigma_{1}=1$ and $\mu_{1}=1$, then $\sigma_{0}\geq7$ and $5\leq\mu_{2}\leq\lfloor4+32/(\sigma_{0}-6)\rfloor$ by $\mu_{1}(\sigma_{0}+2)>8$, $\mu_{1}\mu_{2}>4$ and $\mu_{2}\leq\lfloor4(\sigma_{0}+2)/[\mu_{1}(\sigma_{0}+2)-8]\rfloor$.
If $\sigma_{1}=1$ and $\mu_{1}=1$, then $\mu_{2}\leq36,20,14,12,10$ or $9$ for $\sigma_{0}=7,8,9,10,11$ or $12$, respectively, by $\mu_{2}\leq\lfloor4+32/(\sigma_{0}-6)\rfloor$.
If $\sigma_{1}=1$ and $\mu_{1}=1$, then $\mu_{2}\leq8,7,6$ or $5$ for $13\leq\sigma_{0}\leq14$, $15\leq\sigma_{0}\leq16$, $17\leq\sigma_{0}\leq22$ or $23\leq\sigma_{0}\leq38$, respectively, by $\mu_{2}\leq\lfloor4+32/(\sigma_{0}-6)\rfloor$.
Thus, if $\sigma_{1}=1$, $\mu_{1}=1$ and $\sigma_{0}=7$, then $G\subset D\bigl((7-2b_{0},1,1,1),\allowbreak(b_{0},0,0,0),\allowbreak(1,1,36);0,0\bigr)=A_{i}$ for some $220\leq i\leq223$.
Thus, if $\sigma_{1}=1$, $\mu_{1}=1$ and $\sigma_{0}=8$, then $G\subset D\bigl((8-2b_{0},1,1,1),\allowbreak(b_{0},0,0,0),\allowbreak(1,1,20);0,0\bigr)=A_{i}$ for some $224\leq i\leq228$.
Thus, if $\sigma_{1}=1$, $\mu_{1}=1$ and $\sigma_{0}=9$, then $G\subset D\bigl((9-2b_{0},1,1,1),\allowbreak(b_{0},0,0,0),\allowbreak(1,1,14);0,0\bigr)=A_{i}$ for some $249\leq i\leq253$.
Thus, if $\sigma_{1}=1$, $\mu_{1}=1$ and $\sigma_{0}=10$, then $G\subset D\bigl((10-2b_{0},1,1,1),\allowbreak(b_{0},0,0,0),\allowbreak(1,1,12);0,0\bigr)=A_{i}$ for some $254\leq i\leq259$.
Thus, if $\sigma_{1}=1$, $\mu_{1}=1$ and $\sigma_{0}=11$, then $G\subset D\bigl((11-2b_{0},1,1,1),\allowbreak(b_{0},0,0,0),\allowbreak(1,1,10);0,0\bigr)=A_{i}$ for some $278\leq i\leq283$.
Thus, if $\sigma_{1}=1$, $\mu_{1}=1$ and $\sigma_{0}=12$, then $G\subset D\bigl((12-2b_{0},1,1,1),\allowbreak(b_{0},0,0,0),\allowbreak(1,1,9);0,0\bigr)=A_{i}$ for some $284\leq i\leq290$.
Thus, if $\sigma_{1}=1$, $\mu_{1}=1$ and $13\leq\sigma_{0}\leq14$, then $G\subset D\bigl((14-2b_{0},1,1,1),\allowbreak(b_{0},0,0,0),\allowbreak(1,1,8);0,0\bigr)=A_{i}$ for some $305\leq i\leq312$.
Thus, if $\sigma_{1}=1$, $\mu_{1}=1$ and $15\leq\sigma_{0}\leq16$, then $G\subset D\bigl((16-2b_{0},1,1,1),\allowbreak(b_{0},0,0,0),\allowbreak(1,1,7);0,0\bigr)=A_{i}$ for some $313\leq i\leq321$.
Thus, if $\sigma_{1}=1$, $\mu_{1}=1$ and $17\leq\sigma_{0}\leq22$, then $G\subset D\bigl((22-2b_{0},1,1,1),\allowbreak(b_{0},0,0,0),\allowbreak(1,1,6);0,0\bigr)=A_{i}$ for some $322\leq i\leq333$.
Thus, if $\sigma_{1}=1$, $\mu_{1}=1$ and $23\leq\sigma_{0}\leq38$, then $G\subset D\bigl((38-2b_{0},1,1,1),\allowbreak(b_{0},0,0,0),\allowbreak(1,1,5);0,0\bigr)=A_{i}$ for some $334\leq i\leq353$.
If $\sigma_{1}=1$ and $\mu_{1}=2$, then $3\leq\mu_{2}\leq\lfloor2+8/(\sigma_{0}-2)\rfloor$ by $\mu_{1}\mu_{2}>4$ and $\mu_{2}\leq\lfloor4(\sigma_{0}+2)/[2(\sigma_{0}+2)-8]\rfloor$.
Thus, if $\sigma_{1}=1$ and $\mu_{1}=2$, then $3\leq\sigma_{0}\leq10$ by $3\leq2+8/(\sigma_{0}-2)$.
If $\sigma_{1}=1$ and $\mu_{1}=2$, then $\mu_{2}\leq10,6,4$ or $3$ for $\sigma_{0}=3$, $\sigma_{0}=4$, $5\leq\sigma_{0}\leq6$ or $7\leq\sigma_{0}\leq10$, respectively, by $\mu_{2}\leq\lfloor2+8/(\sigma_{0}-2)\rfloor$.
Thus, if $\sigma_{1}=1$, $\mu_{1}=2$ and $\sigma_{0}=3$, then $G\subset D\bigl((3-2b_{0},1,1,1),\allowbreak(b_{0},0,0,0),\allowbreak(1,2,10);0,0\bigr)=A_{i}$ for some $83\leq i\leq84$.
Thus, if $\sigma_{1}=1$, $\mu_{1}=2$ and $\sigma_{0}=4$, then $G\subset D\bigl((4-2b_{0},1,1,1),\allowbreak(b_{0},0,0,0),\allowbreak(1,2,6);0,0\bigr)=A_{i}$ for some $139\leq i\leq141$.
Thus, if $\sigma_{1}=1$, $\mu_{1}=2$ and $5\leq\sigma_{0}\leq6$, then $G\subset D\bigl((6-2b_{0},1,1,1),\allowbreak(b_{0},0,0,0),\allowbreak(1,2,4);0,0\bigr)=A_{i}$ for some $184\leq i\leq187$.
Thus, if $\sigma_{1}=1$, $\mu_{1}=2$ and $7\leq\sigma_{0}\leq10$, then $G\subset D\bigl((10-2b_{0},1,1,1),\allowbreak(b_{0},0,0,0),\allowbreak(1,2,3);0,0\bigr)=A_{i}$ for some $260\leq i\leq265$.
If $\sigma_{1}=1$ and $\mu_{1}=3$, then $\mu_{2}=2$ and $2[3(\sigma_{0}+2)-8]\leq4(\sigma_{0}+2)$, which implies that $3\leq\sigma_{0}\leq6$.
If $\sigma_{1}=1$ and $\mu_{1}=5$, then $\mu_{2}=1$ and $5(\sigma_{0}+2)-8\leq4(\sigma_{0}+2)$, which implies that $3\leq\sigma_{0}\leq6$.
Thus, if $\sigma_{1}=1$ and $\mu_{1}=3$, then $G\subset D\bigl((6-2b_{0},1,1,1),\allowbreak(b_{0},0,0,0),\allowbreak(1,3,2);0,0\bigr)=A_{i}$ for some $188\leq i\leq191$.
Thus, if $\sigma_{1}=1$ and $\mu_{1}=5$, then $G\subset D\bigl((6-2b_{0},1,1,1),\allowbreak(b_{0},0,0,0),\allowbreak(1,5,1);0,0\bigr)=A_{i}$ for some $192\leq i\leq195$.
If $\sigma_{1}=2$, then $r_{1}=2/(\sigma_{0}-2+2\mu_{1})$.
Thus, if $\sigma_{1}=2$, then $\sigma_{0}+2\mu_{1}>6$ and $1\leq\mu_{2}\leq\lfloor8/(\sigma_{0}+2\mu_{1}-6)\rfloor$ by $r_{1}<1/2$ and $\mu_{2}\leq\lfloor2/(1-2r_{1})-2\rfloor$.
Thus, if $\sigma_{1}=2$, then $7\leq\sigma_{0}+2\mu_{1}\leq14$.
If $\sigma_{1}=2$ and $\mu_{1}=1$, then $5\leq\sigma_{0}\leq12$ and $1\leq\mu_{2}\leq\lfloor8/(\sigma_{0}-4)\rfloor$.
If $\sigma_{1}=2$ and $\mu_{1}=1$, then $\mu_{2}\leq8,4,2$ or $1$ for $\sigma_{0}=5$, $\sigma_{0}=6$, $7\leq\sigma_{0}\leq8$ or $9\leq\sigma_{0}\leq12$, respectively, by $\mu_{2}\leq\lfloor8/(\sigma_{0}-4)\rfloor$.
Thus, if $\sigma_{1}=2$, $\mu_{1}=1$ and $\sigma_{0}=5$, then $G\subset D\bigl((5-2b_{0},2-2b_{1},1,1),\allowbreak(b_{0},b_{1},0,0),\allowbreak(1,1,8);0,0\bigr)=A_{i}$ for some $178\leq i\leq183$.
Thus, if $\sigma_{1}=2$, $\mu_{1}=1$ and $\sigma_{0}=6$, then $G\subset D\bigl((6-2b_{0},2-2b_{1},1,1),\allowbreak(b_{0},b_{1},0,0),\allowbreak(1,1,4);0,0\bigr)=A_{i}$ for some $196\leq i\leq203$.
Thus, if $\sigma_{1}=2$, $\mu_{1}=1$ and $7\leq\sigma_{0}\leq8$, then $G\subset D\bigl((8-2b_{0},2-2b_{1},1,1),\allowbreak(b_{0},b_{1},0,0),\allowbreak(1,1,2);0,0\bigr)=A_{i}$ for some $229\leq i\leq238$.
Thus, if $\sigma_{1}=2$, $\mu_{1}=1$ and $9\leq\sigma_{0}\leq12$, then $G\subset D\bigl((12-2b_{0},2-2b_{1},1,1),\allowbreak(b_{0},b_{1},0,0),\allowbreak(1,1,1);0,0\bigr)=A_{i}$ for some $291\leq i\leq304$.
If $\sigma_{1}=2$ and $\mu_{1}=2$, then $3\leq\sigma_{0}\leq10$ and $1\leq\mu_{2}\leq\lfloor8/(\sigma_{0}-2)\rfloor$.
If $\sigma_{1}=2$ and $\mu_{1}=2$, then $\mu_{2}\leq8,4,2$ or $1$ for $\sigma_{0}=3$, $\sigma_{0}=4$, $5\leq\sigma_{0}\leq6$ or $7\leq\sigma_{0}\leq10$, respectively, by $\mu_{2}\leq\lfloor8/(\sigma_{0}-2)\rfloor$.
Thus, if $\sigma_{1}=2$, $\mu_{1}=2$ and $\sigma_{0}=3$, then $G\subset D\bigl((3-2b_{0},2-2b_{1},1,1),\allowbreak(b_{0},b_{1},0,0),\allowbreak(1,2,8);0,0\bigr)=A_{i}$ for some $87\leq i\leq90$.
Thus, if $\sigma_{1}=2$, $\mu_{1}=2$ and $\sigma_{0}=4$, then $G\subset D\bigl((4-2b_{0},2-2b_{1},1,1),\allowbreak(b_{0},b_{1},0,0),\allowbreak(1,2,4);0,0\bigr)=A_{i}$ for some $142\leq i\leq147$.
Thus, if $\sigma_{1}=2$, $\mu_{1}=2$ and $5\leq\sigma_{0}\leq6$, then $G\subset D\bigl((6-2b_{0},2-2b_{1},1,1),\allowbreak(b_{0},b_{1},0,0),\allowbreak(1,2,2);0,0\bigr)=A_{i}$ for some $204\leq i\leq211$.
Thus, if $\sigma_{1}=2$, $\mu_{1}=2$ and $7\leq\sigma_{0}\leq10$, then $G\subset D\bigl((10-2b_{0},2-2b_{1},1,1),\allowbreak(b_{0},b_{1},0,0),\allowbreak(1,2,1);0,0\bigr)=A_{i}$ for some $266\leq i\leq277$.
If $\sigma_{1}=2$ and $\mu_{1}=3$, then $3\leq\sigma_{0}\leq8$ and $1\leq\mu_{2}\leq\lfloor8/\sigma_{0}\rfloor$.
If $\sigma_{1}=2$, $\mu_{1}=3$ and $3\leq\sigma_{0}\leq4$, then $1\leq\mu_{2}\leq2$.
If $\sigma_{1}=2$, $\mu_{1}=3$ and $5\leq\sigma_{0}\leq8$, then $\mu_{2}=1$.
Thus, if $\sigma_{1}=2$, $\mu_{1}=3$ and $3\leq\sigma_{0}\leq4$, then $G\subset D\bigl((4-2b_{0},2-2b_{1},1,1),\allowbreak(b_{0},b_{1},0,0),\allowbreak(1,3,2);0,0\bigr)=A_{i}$ for some $148\leq i\leq153$.
Thus, if $\sigma_{1}=2$, $\mu_{1}=3$ and $5\leq\sigma_{0}\leq8$, then $G\subset D\bigl((8-2b_{0},2-2b_{1},1,1),\allowbreak(b_{0},b_{1},0,0),\allowbreak(1,3,1);0,0\bigr)=A_{i}$ for some $239\leq i\leq248$.
If $\sigma_{1}=2$ and $\mu_{1}=4$, then $3\leq\sigma_{0}\leq6$ and $\mu_{2}=1$ by $7\leq\sigma_{0}+2\mu_{1}\leq14$ and $1\leq\mu_{2}\leq\lfloor8/(\sigma_{0}+2)\rfloor$.
If $\sigma_{1}=2$ and $\mu_{1}=5$, then $3\leq\sigma_{0}\leq4$ and $\mu_{2}=1$ by $7\leq\sigma_{0}+2\mu_{1}\leq14$ and $1\leq\mu_{2}\leq\lfloor8/(\sigma_{0}+4)\rfloor$.
Thus, if $\sigma_{1}=2$ and $\mu_{1}=4$, then $G\subset D\bigl((6-2b_{0},2-2b_{1},1,1),\allowbreak(b_{0},b_{1},0,0),\allowbreak(1,4,1);0,0\bigr)=A_{i}$ for some $212\leq i\leq219$.
Thus, if $\sigma_{1}=2$ and $\mu_{1}=5$, then $G\subset D\bigl((4-2b_{0},2-2b_{1},1,1),\allowbreak(b_{0},b_{1},0,0),\allowbreak(1,5,1);0,0\bigr)=A_{i}$ for some $154\leq i\leq159$.
If $\sigma_{1}=3$, then $r_{0}\geq3/2+1/(6-\mu_{1})\geq17/10$, which implies that $3\leq\sigma_{0}\leq4$ by $r_{0}=\sigma_{0}/(\sigma_{0}-2)$.
If $\sigma_{1}=3$ and $\sigma_{0}=3$, then $r_{0}=3$ and $r_{1}=3/(2+3\mu_{1})$.
Thus, if $\sigma_{1}=3$ and $\sigma_{0}=3$, then $2\leq\mu_{1}\leq5$ and $1\leq\mu_{2}\leq\lfloor12/(3\mu_{1}-4)\rfloor$ by $r_{1}<1/2$ and $\mu_{2}\leq\lfloor2/(1-2r_{1})-2\rfloor$.
If $\sigma_{1}=3$ and $\sigma_{0}=3$, then $\mu_{2}\leq6,2$ or $1$ for $\mu_{1}=2$, $\mu_{1}=3$ or $4\leq\mu_{1}\leq5$, respectively, by $\mu_{2}\leq\lfloor12/(3\mu_{1}-4)\rfloor$.
Thus, if $\sigma_{1}=3$, $\sigma_{0}=3$ and $\mu_{1}=2$, then $G\subset D\bigl((3-2b_{0},3-2b_{1},1,1),\allowbreak(b_{0},b_{1},0,0),\allowbreak(1,2,6);0,0\bigr)=A_{i}$ for some $99\leq i\leq102$.
Thus, if $\sigma_{1}=3$, $\sigma_{0}=3$ and $\mu_{1}=3$, then $G\subset D\bigl((3-2b_{0},4-2b_{1},1,1),\allowbreak(b_{0},b_{1},0,0),\allowbreak(1,3,2);0,0\bigr)=A_{i}$ for some $109\leq i\leq114$.
Thus, if $\sigma_{1}=3$, $\sigma_{0}=3$ and $4\leq\mu_{1}\leq5$, then $G\subset D\bigl((3-2b_{0},4-2b_{1},1,1),\allowbreak(b_{0},b_{1},0,0),\allowbreak(1,5,1);0,0\bigr)=A_{i}$ for some $115\leq i\leq120$.
If $\sigma_{1}=3$ and $\sigma_{0}=4$, then $r_{0}=2$ and $r_{1}=1/(\mu_{1}+2)$.
If $\sigma_{1}=3$ and $\sigma_{0}=4$, then $1\leq\mu_{2}\leq\lfloor4/\mu_{1}\rfloor$ by $\mu_{2}\leq\lfloor2/(1-2r_{1})-2\rfloor$.
Thus, if $\sigma_{1}=3$ and $\sigma_{0}=4$, then $1\leq\mu_{1}\leq4$.
If $\sigma_{1}=3$ and $\sigma_{0}=4$, then $\mu_{2}\leq4,2$ or $1$ for $\mu_{1}=1$, $\mu_{1}=2$ or $3\leq\mu_{1}\leq4$, respectively, by $\mu_{2}\leq\lfloor4/\mu_{1}\rfloor$.
Thus, if $\sigma_{1}=3$, $\sigma_{0}=4$ and $\mu_{1}=1$, then $G\subset D\bigl((4-2b_{0},3-2b_{1},1,1),\allowbreak(b_{0},b_{1},0,0),\allowbreak(1,1,4);0,0\bigr)=A_{i}$ for some $160\leq i\leq165$.
Thus, if $\sigma_{1}=3$, $\sigma_{0}=4$ and $\mu_{1}=2$, then $G\subset D\bigl((4-2b_{0},3-2b_{1},1,1),\allowbreak(b_{0},b_{1},0,0),\allowbreak(1,2,2);0,0\bigr)=A_{i}$ for some $166\leq i\leq171$.
Thus, if $\sigma_{1}=3$, $\sigma_{0}=4$ and $3\leq\mu_{1}\leq4$, then $G\subset D\bigl((4-2b_{0},3-2b_{1},1,1),\allowbreak(b_{0},b_{1},0,0),\allowbreak(1,4,1);0,0\bigr)=A_{i}$ for some $172\leq i\leq177$.
If $\sigma_{1}=4$, then $r_{0}\geq2+1/(6-\mu_{1})\geq11/5$, which implies that $\sigma_{0}=3$ by $r_{0}=\sigma_{0}/(\sigma_{0}-2)$ and $\sigma_{0}\geq3$.
If $\sigma_{1}=4$, then $r_{0}=3$ and $r_{1}=1/(\mu_{1}+1)$.
Thus, if $\sigma_{1}=4$, then $2\leq\mu_{1}\leq5$ and $1\leq\mu_{2}\leq\lfloor4/(\mu_{1}-1)\rfloor$ by $r_{1}<1/2$ and $\mu_{2}\leq\lfloor2/(1-2r_{1})-2\rfloor$.
If $\sigma_{1}=4$, then $\mu_{2}\leq4,2$ or $1$ for $\mu_{1}=2$, $\mu_{1}=3$ or $4\leq\mu_{1}\leq5$, respectively, by $\mu_{2}\leq\lfloor4/(\mu_{1}-1)\rfloor$.
Thus, if $\sigma_{1}=4$ and $\mu_{1}=2$, then $G\subset D\bigl((3-2b_{0},4-2b_{1},1,1),\allowbreak(b_{0},b_{1},0,0),\allowbreak(1,2,4);0,0\bigr)=A_{i}$ for some $103\leq i\leq108$.
Thus, if $\sigma_{1}=4$ and $\mu_{1}=3$, then $G\subset D\bigl((3-2b_{0},4-2b_{1},1,1),\allowbreak(b_{0},b_{1},0,0),\allowbreak(1,3,2);0,0\bigr)=A_{i}$ for some $109\leq i\leq114$.
Thus, if $\sigma_{1}=4$ and $4\leq\mu_{1}\leq5$, then $G\subset D\bigl((3-2b_{0},4-2b_{1},1,1),\allowbreak(b_{0},b_{1},0,0),\allowbreak(1,5,1);0,0\bigr)=A_{i}$ for some $115\leq i\leq120$.
If $\sigma_{1}=5$, then $r_{0}\geq5/2+1/(6-\mu_{1})\geq27/10$, which implies that $\sigma_{0}=3$ by $r_{0}=\sigma_{0}/(\sigma_{0}-2)$ and $\sigma_{0}\geq3$.
If $\sigma_{1}=5$, then $r_{0}=3$ and $r_{1}=1/(\mu_{1}+2)$.
If $\sigma_{1}=5$, then $1\leq\mu_{2}\leq\lfloor4/\mu_{1}\rfloor$ by $\mu_{2}\leq\lfloor2/(1-2r_{1})-2\rfloor$.
Thus, if $\sigma_{1}=5$, then $1\leq\mu_{1}\leq4$.
If $\sigma_{1}=5$, then $\mu_{2}\leq4,2$ or $1$ for $\mu_{1}=1$, $\mu_{1}=2$ or $3\leq\mu_{1}\leq4$, respectively, by $\mu_{2}\leq\lfloor4/\mu_{1}\rfloor$.
Thus, if $\sigma_{1}=5$ and $\mu_{1}=1$, then $G\subset D\bigl((3-2b_{0},5-2b_{1},1,1),\allowbreak(b_{0},b_{1},0,0),\allowbreak(1,1,4);0,0\bigr)=A_{i}$ for some $121\leq i\leq126$.
Thus, if $\sigma_{1}=5$ and $\mu_{1}=2$, then $G\subset D\bigl((3-2b_{0},5-2b_{1},1,1),\allowbreak(b_{0},b_{1},0,0),\allowbreak(1,2,2);0,0\bigr)=A_{i}$ for some $127\leq i\leq132$.
Thus, if $\sigma_{1}=5$ and $3\leq\mu_{1}\leq4$, then $G\subset D\bigl((3-2b_{0},5-2b_{1},1,1),\allowbreak(b_{0},b_{1},0,0),\allowbreak(1,4,1);0,0\bigr)=A_{i}$ for some $133\leq i\leq138$.

Suppose that $\mu_{0}=2$.
Then $2\leq\sigma_{0}\leq3$ by Fact~\ref{f23bounds} and $\mu_{0}\sigma_{0}>2$. Then $\sigma_{1}=1$ by $r_{0}\leq1$ and $r_{1}>0$. If $\sigma_{0}=2$, then $\mu_{1}\mu_{2}\leq4$ by $r_{0}=1$ and $r_{1}\geq\mu_{2}/[2(\mu_{2}+2)]$. If $\sigma_{0}=3$, then $(\mu_{1}+2)\mu_{2}\leq4$ by $r_{0}=3/4$ and $r_{1}\geq\mu_{2}/[2(\mu_{2}+2)]$.
If $\sigma_{0}\leq2$, $\sigma_{1}\leq1$, $\mu_{1}\leq1$ and $\mu_{2}\leq4$, then $G\subset D\bigl((2-2b_{0},1,1,1),\allowbreak(b_{0},0,0,0),\allowbreak(2,1,4);0,0\bigr)\subset\mathcal F_{92}(n)$.
If $\sigma_{0}\leq2$, $\sigma_{1}\leq1$, $\mu_{1}\leq2$ and $\mu_{2}\leq2$, then $G\subset D\bigl((2-2b_{0},1,1,1),\allowbreak(b_{0},0,0,0),\allowbreak(2,2,2);0,0\bigr)\subset\mathcal F_{93}(n)$.
If $\sigma_{0}\leq2$, $\sigma_{1}\leq1$, $\mu_{1}\leq4$ and $\mu_{2}\leq1$, then $G\subset D\bigl((2-2b_{0},1,1,1),\allowbreak(b_{0},0,0,0),\allowbreak(2,4,1);0,0\bigr)\subset\mathcal F_{94}(n)$.
If $\sigma_{0}\leq3$, $\sigma_{1}\leq1$, $\mu_{1}\leq2$ and $\mu_{2}\leq1$, then $G\subset D\bigl((3-2b_{0},1,1,1),\allowbreak(b_{0},0,0,0),\allowbreak(2,2,1);0,0\bigr)=A_{i}$ for some $85\leq i\leq86$.

Suppose that $\mu_{0}=3$.
Then $\sigma_{0}=1$ and $r_{0}=1$ by Fact~\ref{f23bounds} and $r_{0}=\sigma_{0}/(3\sigma_{0}-2)$. Then $\sigma_{1}=1$ by $r_{1}>0$ and $r_{0}=1$. Thus, $\mu_{1}\mu_{2}\leq4$ by $r_{1}\geq\mu_{2}/[2(\mu_{2}+2)]$.
If $\sigma_{0}\leq1$, $\sigma_{1}\leq1$, $\mu_{1}\leq1$ and $\mu_{2}\leq4$, then $G\subset D\bigl((1-2b_{0},1,1,1),\allowbreak(b_{0},0,0,0),\allowbreak(3,1,4);0,0\bigr)\subset\mathcal F_{92}(n)$.
If $\sigma_{0}\leq1$, $\sigma_{1}\leq1$, $\mu_{1}\leq2$ and $\mu_{2}\leq2$, then $G\subset D\bigl((1-2b_{0},1,1,1),\allowbreak(b_{0},0,0,0),\allowbreak(3,2,2);0,0\bigr)\subset\mathcal F_{93}(n)$.
If $\sigma_{0}\leq1$, $\sigma_{1}\leq1$, $\mu_{1}\leq4$ and $\mu_{2}\leq1$, then $G\subset D\bigl((1-2b_{0},1,1,1),\allowbreak(b_{0},0,0,0),\allowbreak(3,4,1);0,0\bigr)\subset\mathcal F_{94}(n)$.

\noindent$\textbf{Subcase 2.4.}$ $h=4$.

\begin{fact}\label{f24bounds}
Then $\boldsymbol\sigma=(3,1,1,1,1)$ and $\boldsymbol\mu=(1,1,1,1)$.
\end{fact}
\begin{proof}
Then $r_{3}\geq\sigma_{4}/2-1\geq-1/2$.
By $r_{3}\geq-1/2$ and $r_{3}=(r_{2}-\sigma_{3}/2)/(1+\mu_{3}(r_{2}-\sigma_{3}/2))$, $r_{2}\geq\sigma_{3}/2-1/(\mu_{3}+2)\geq1/6$.
Since $r_{j}<1/\mu_{j}$ for $j=1,2$, the expression for $r_{j}$ gives $r_{j-1}=\sigma_{j}/2+r_{j}/(1-\mu_{j}r_{j})$.
Thus, $r_{1}\geq1/2+(1/6)/(1-1/6)=7/10$ and $r_{0}\geq1/2+(7/10)/(1-7/10)=17/6$.

Since $r_{0}=1/(\mu_{0}-2/\sigma_{0})$, $\mu_{0}\leq2/\sigma_{0}+6/17$.
If $\sigma_{0}=1$, then $\mu_{0}\leq2$; if $\sigma_{0}=2$, then $\mu_{0}\leq1$.
Both cases contradict $\mu_{0}\sigma_{0}>2$.
Thus, $\sigma_{0}\geq3$, which implies that $\mu_{0}=1$.
Then $\sigma_{0}/(\sigma_{0}-2)=r_{0}\geq17/6$, and hence $\sigma_{0}\leq34/11$.
Thus, $\sigma_{0}=3$ and $r_{0}=3$.

Since $r_{1}\geq7/10$ and $r_{1}<1/\mu_{1}$, $\mu_{1}=1$.
Then $3=r_{0}=\sigma_{1}/2+r_{1}/(1-r_{1})\geq\sigma_{1}/2+7/3$.
Thus, $\sigma_{1}=1$ and $r_{1}=5/7$.
Since $r_{2}\geq1/6$ and $r_{2}<1/\mu_{2}$, $\mu_{2}\leq5$.
Then $5/7=r_{1}\geq\sigma_{2}/2+1/(6-\mu_{2})$, which implies that $\sigma_{2}=1$ and $1/(6-\mu_{2})\leq3/14$.
Thus, $\mu_{2}=1$ and $r_{2}=3/17$.

Since $3/17=r_{2}\geq\sigma_{3}/2-1/(\mu_{3}+2)\geq\sigma_{3}/2-1/3$, $\sigma_{3}=1$.
Then $1/(\mu_{3}+2)\geq1/2-3/17=11/34$, which implies that $\mu_{3}=1$.
Thus, $r_{3}=-11/23$.
Since $-11/23=r_{3}\geq\sigma_{4}/2-1$, $\sigma_{4}=1$.
\end{proof}

Then $b_{1}=b_{2}=b_{3}=b_{4}=0$ and $0\leq b_{0}\leq1$ by Fact~\ref{f24bounds}.
If $b_{0}=0$, then $G\cong D\bigl((3,1,1,1,1),(0,0,0,0,0),(1,1,1,1);0,0\bigr)=A_{355}$.
If $b_{0}=1$, then $G\cong D\bigl((1,1,1,1,1),(1,0,0,0,0),(1,1,1,1);0,0\bigr)=A_{354}$.
\end{proof}

\begin{table}[p]
\centering
\begingroup\footnotesize
\setlength{\tabcolsep}{3pt}
\renewcommand{\arraystretch}{1.0}
\caption{Graphs $A_{1},\ldots,A_{355}$.}
\label{tab:s5-explicit-1}
\begin{tabular*}{\textwidth}{@{\extracolsep{\fill}}cccc@{\hspace{7mm}}cccc}
\toprule
$A_{i}$ & $\boldsymbol a$ & $\boldsymbol b$ & $\boldsymbol\mu$ & $A_{i}$ & $\boldsymbol a$ & $\boldsymbol b$ & $\boldsymbol\mu$\\
\midrule
$A_{1}$ & $(1,1,1)$ & $(1,0,0)$ & $(3,12)$ & $A_{46}$ & $(6,3,1)$ & $(1,0,0)$ & $(1,4)$\\
$A_{2}$ & $(3,1,1)$ & $(0,0,0)$ & $(3,12)$ & $A_{47}$ & $(8,1,1)$ & $(0,1,0)$ & $(1,4)$\\
$A_{3}$ & $(1,1,1)$ & $(1,0,0)$ & $(4,3)$ & $A_{48}$ & $(8,3,1)$ & $(0,0,0)$ & $(1,4)$\\
$A_{4}$ & $(3,1,1)$ & $(0,0,0)$ & $(4,3)$ & $A_{49}$ & $(0,1,1)$ & $(5,0,0)$ & $(3,5)$\\
$A_{5}$ & $(1,0,1)$ & $(1,1,0)$ & $(2,2)$ & $A_{50}$ & $(2,1,1)$ & $(4,0,0)$ & $(3,5)$\\
$A_{6}$ & $(1,2,1)$ & $(1,0,0)$ & $(2,2)$ & $A_{51}$ & $(4,1,1)$ & $(3,0,0)$ & $(3,5)$\\
$A_{7}$ & $(3,0,1)$ & $(0,1,0)$ & $(2,2)$ & $A_{52}$ & $(6,1,1)$ & $(2,0,0)$ & $(3,5)$\\
$A_{8}$ & $(3,2,1)$ & $(0,0,0)$ & $(2,2)$ & $A_{53}$ & $(8,1,1)$ & $(1,0,0)$ & $(3,5)$\\
$A_{9}$ & $(0,1,1)$ & $(2,0,0)$ & $(3,8)$ & $A_{54}$ & $(10,1,1)$ & $(0,0,0)$ & $(3,5)$\\
$A_{10}$ & $(2,1,1)$ & $(1,0,0)$ & $(3,8)$ & $A_{55}$ & $(0,1,1)$ & $(5,1,0)$ & $(1,2)$\\
$A_{11}$ & $(4,1,1)$ & $(0,0,0)$ & $(3,8)$ & $A_{56}$ & $(0,3,1)$ & $(5,0,0)$ & $(1,2)$\\
$A_{12}$ & $(0,0,1)$ & $(2,1,0)$ & $(2,1)$ & $A_{57}$ & $(2,1,1)$ & $(4,1,0)$ & $(1,2)$\\
$A_{13}$ & $(0,2,1)$ & $(2,0,0)$ & $(2,1)$ & $A_{58}$ & $(2,3,1)$ & $(4,0,0)$ & $(1,2)$\\
$A_{14}$ & $(2,0,1)$ & $(1,1,0)$ & $(2,1)$ & $A_{59}$ & $(4,1,1)$ & $(3,1,0)$ & $(1,2)$\\
$A_{15}$ & $(2,2,1)$ & $(1,0,0)$ & $(2,1)$ & $A_{60}$ & $(4,3,1)$ & $(3,0,0)$ & $(1,2)$\\
$A_{16}$ & $(4,0,1)$ & $(0,1,0)$ & $(2,1)$ & $A_{61}$ & $(6,1,1)$ & $(2,1,0)$ & $(1,2)$\\
$A_{17}$ & $(4,2,1)$ & $(0,0,0)$ & $(2,1)$ & $A_{62}$ & $(6,3,1)$ & $(2,0,0)$ & $(1,2)$\\
$A_{18}$ & $(1,0,1)$ & $(2,2,0)$ & $(1,1)$ & $A_{63}$ & $(8,1,1)$ & $(1,1,0)$ & $(1,2)$\\
$A_{19}$ & $(1,2,1)$ & $(2,1,0)$ & $(1,1)$ & $A_{64}$ & $(8,3,1)$ & $(1,0,0)$ & $(1,2)$\\
$A_{20}$ & $(1,4,1)$ & $(2,0,0)$ & $(1,1)$ & $A_{65}$ & $(10,1,1)$ & $(0,1,0)$ & $(1,2)$\\
$A_{21}$ & $(3,0,1)$ & $(1,2,0)$ & $(1,1)$ & $A_{66}$ & $(10,3,1)$ & $(0,0,0)$ & $(1,2)$\\
$A_{22}$ & $(3,2,1)$ & $(1,1,0)$ & $(1,1)$ & $A_{67}$ & $(0,1,1)$ & $(7,1,0)$ & $(1,1)$\\
$A_{23}$ & $(3,4,1)$ & $(1,0,0)$ & $(1,1)$ & $A_{68}$ & $(0,3,1)$ & $(7,0,0)$ & $(1,1)$\\
$A_{24}$ & $(5,0,1)$ & $(0,2,0)$ & $(1,1)$ & $A_{69}$ & $(2,1,1)$ & $(6,1,0)$ & $(1,1)$\\
$A_{25}$ & $(5,2,1)$ & $(0,1,0)$ & $(1,1)$ & $A_{70}$ & $(2,3,1)$ & $(6,0,0)$ & $(1,1)$\\
$A_{26}$ & $(5,4,1)$ & $(0,0,0)$ & $(1,1)$ & $A_{71}$ & $(4,1,1)$ & $(5,1,0)$ & $(1,1)$\\
$A_{27}$ & $(0,1,1)$ & $(3,0,0)$ & $(3,6)$ & $A_{72}$ & $(4,3,1)$ & $(5,0,0)$ & $(1,1)$\\
$A_{28}$ & $(2,1,1)$ & $(2,0,0)$ & $(3,6)$ & $A_{73}$ & $(6,1,1)$ & $(4,1,0)$ & $(1,1)$\\
$A_{29}$ & $(4,1,1)$ & $(1,0,0)$ & $(3,6)$ & $A_{74}$ & $(6,3,1)$ & $(4,0,0)$ & $(1,1)$\\
$A_{30}$ & $(6,1,1)$ & $(0,0,0)$ & $(3,6)$ & $A_{75}$ & $(8,1,1)$ & $(3,1,0)$ & $(1,1)$\\
$A_{31}$ & $(1,1,1)$ & $(3,1,0)$ & $(1,8)$ & $A_{76}$ & $(8,3,1)$ & $(3,0,0)$ & $(1,1)$\\
$A_{32}$ & $(1,3,1)$ & $(3,0,0)$ & $(1,8)$ & $A_{77}$ & $(10,1,1)$ & $(2,1,0)$ & $(1,1)$\\
$A_{33}$ & $(3,1,1)$ & $(2,1,0)$ & $(1,8)$ & $A_{78}$ & $(10,3,1)$ & $(2,0,0)$ & $(1,1)$\\
$A_{34}$ & $(3,3,1)$ & $(2,0,0)$ & $(1,8)$ & $A_{79}$ & $(12,1,1)$ & $(1,1,0)$ & $(1,1)$\\
$A_{35}$ & $(5,1,1)$ & $(1,1,0)$ & $(1,8)$ & $A_{80}$ & $(12,3,1)$ & $(1,0,0)$ & $(1,1)$\\
$A_{36}$ & $(5,3,1)$ & $(1,0,0)$ & $(1,8)$ & $A_{81}$ & $(14,1,1)$ & $(0,1,0)$ & $(1,1)$\\
$A_{37}$ & $(7,1,1)$ & $(0,1,0)$ & $(1,8)$ & $A_{82}$ & $(14,3,1)$ & $(0,0,0)$ & $(1,1)$\\
$A_{38}$ & $(7,3,1)$ & $(0,0,0)$ & $(1,8)$ & $A_{83}$ & $(1,1,1,1)$ & $(1,0,0,0)$ & $(1,2,10)$\\
$A_{39}$ & $(0,1,1)$ & $(4,1,0)$ & $(1,4)$ & $A_{84}$ & $(3,1,1,1)$ & $(0,0,0,0)$ & $(1,2,10)$\\
$A_{40}$ & $(0,3,1)$ & $(4,0,0)$ & $(1,4)$ & $A_{85}$ & $(1,1,1,1)$ & $(1,0,0,0)$ & $(2,2,1)$\\
$A_{41}$ & $(2,1,1)$ & $(3,1,0)$ & $(1,4)$ & $A_{86}$ & $(3,1,1,1)$ & $(0,0,0,0)$ & $(2,2,1)$\\
$A_{42}$ & $(2,3,1)$ & $(3,0,0)$ & $(1,4)$ & $A_{87}$ & $(1,0,1,1)$ & $(1,1,0,0)$ & $(1,2,8)$\\
$A_{43}$ & $(4,1,1)$ & $(2,1,0)$ & $(1,4)$ & $A_{88}$ & $(1,2,1,1)$ & $(1,0,0,0)$ & $(1,2,8)$\\
$A_{44}$ & $(4,3,1)$ & $(2,0,0)$ & $(1,4)$ & $A_{89}$ & $(3,0,1,1)$ & $(0,1,0,0)$ & $(1,2,8)$\\
$A_{45}$ & $(6,1,1)$ & $(1,1,0)$ & $(1,4)$ & $A_{90}$ & $(3,2,1,1)$ & $(0,0,0,0)$ & $(1,2,8)$\\
\bottomrule
\end{tabular*}
\endgroup
\end{table}
\begin{table}[p]
\centering
\begingroup\footnotesize
\setlength{\tabcolsep}{3pt}
\renewcommand{\arraystretch}{1.0}
\addtocounter{table}{-1}
\caption{(continued)}
\label{tab:s5-explicit-2}
\begin{tabular*}{\textwidth}{@{\extracolsep{\fill}}cccc@{\hspace{7mm}}cccc}
\toprule
$A_{i}$ & $\boldsymbol a$ & $\boldsymbol b$ & $\boldsymbol\mu$ & $A_{i}$ & $\boldsymbol a$ & $\boldsymbol b$ & $\boldsymbol\mu$\\
\midrule
$A_{91}$ & $(1,0,0,1)$ & $(1,1,1,0)$ & $(1,1,1)$ & $A_{136}$ & $(3,1,1,1)$ & $(0,2,0,0)$ & $(1,4,1)$\\
$A_{92}$ & $(1,0,2,1)$ & $(1,1,0,0)$ & $(1,1,1)$ & $A_{137}$ & $(3,3,1,1)$ & $(0,1,0,0)$ & $(1,4,1)$\\
$A_{93}$ & $(1,2,0,1)$ & $(1,0,1,0)$ & $(1,1,1)$ & $A_{138}$ & $(3,5,1,1)$ & $(0,0,0,0)$ & $(1,4,1)$\\
$A_{94}$ & $(1,2,2,1)$ & $(1,0,0,0)$ & $(1,1,1)$ & $A_{139}$ & $(0,1,1,1)$ & $(2,0,0,0)$ & $(1,2,6)$\\
$A_{95}$ & $(3,0,0,1)$ & $(0,1,1,0)$ & $(1,1,1)$ & $A_{140}$ & $(2,1,1,1)$ & $(1,0,0,0)$ & $(1,2,6)$\\
$A_{96}$ & $(3,0,2,1)$ & $(0,1,0,0)$ & $(1,1,1)$ & $A_{141}$ & $(4,1,1,1)$ & $(0,0,0,0)$ & $(1,2,6)$\\
$A_{97}$ & $(3,2,0,1)$ & $(0,0,1,0)$ & $(1,1,1)$ & $A_{142}$ & $(0,0,1,1)$ & $(2,1,0,0)$ & $(1,2,4)$\\
$A_{98}$ & $(3,2,2,1)$ & $(0,0,0,0)$ & $(1,1,1)$ & $A_{143}$ & $(0,2,1,1)$ & $(2,0,0,0)$ & $(1,2,4)$\\
$A_{99}$ & $(1,1,1,1)$ & $(1,1,0,0)$ & $(1,2,6)$ & $A_{144}$ & $(2,0,1,1)$ & $(1,1,0,0)$ & $(1,2,4)$\\
$A_{100}$ & $(1,3,1,1)$ & $(1,0,0,0)$ & $(1,2,6)$ & $A_{145}$ & $(2,2,1,1)$ & $(1,0,0,0)$ & $(1,2,4)$\\
$A_{101}$ & $(3,1,1,1)$ & $(0,1,0,0)$ & $(1,2,6)$ & $A_{146}$ & $(4,0,1,1)$ & $(0,1,0,0)$ & $(1,2,4)$\\
$A_{102}$ & $(3,3,1,1)$ & $(0,0,0,0)$ & $(1,2,6)$ & $A_{147}$ & $(4,2,1,1)$ & $(0,0,0,0)$ & $(1,2,4)$\\
$A_{103}$ & $(1,0,1,1)$ & $(1,2,0,0)$ & $(1,2,4)$ & $A_{148}$ & $(0,0,1,1)$ & $(2,1,0,0)$ & $(1,3,2)$\\
$A_{104}$ & $(1,2,1,1)$ & $(1,1,0,0)$ & $(1,2,4)$ & $A_{149}$ & $(0,2,1,1)$ & $(2,0,0,0)$ & $(1,3,2)$\\
$A_{105}$ & $(1,4,1,1)$ & $(1,0,0,0)$ & $(1,2,4)$ & $A_{150}$ & $(2,0,1,1)$ & $(1,1,0,0)$ & $(1,3,2)$\\
$A_{106}$ & $(3,0,1,1)$ & $(0,2,0,0)$ & $(1,2,4)$ & $A_{151}$ & $(2,2,1,1)$ & $(1,0,0,0)$ & $(1,3,2)$\\
$A_{107}$ & $(3,2,1,1)$ & $(0,1,0,0)$ & $(1,2,4)$ & $A_{152}$ & $(4,0,1,1)$ & $(0,1,0,0)$ & $(1,3,2)$\\
$A_{108}$ & $(3,4,1,1)$ & $(0,0,0,0)$ & $(1,2,4)$ & $A_{153}$ & $(4,2,1,1)$ & $(0,0,0,0)$ & $(1,3,2)$\\
$A_{109}$ & $(1,0,1,1)$ & $(1,2,0,0)$ & $(1,3,2)$ & $A_{154}$ & $(0,0,1,1)$ & $(2,1,0,0)$ & $(1,5,1)$\\
$A_{110}$ & $(1,2,1,1)$ & $(1,1,0,0)$ & $(1,3,2)$ & $A_{155}$ & $(0,2,1,1)$ & $(2,0,0,0)$ & $(1,5,1)$\\
$A_{111}$ & $(1,4,1,1)$ & $(1,0,0,0)$ & $(1,3,2)$ & $A_{156}$ & $(2,0,1,1)$ & $(1,1,0,0)$ & $(1,5,1)$\\
$A_{112}$ & $(3,0,1,1)$ & $(0,2,0,0)$ & $(1,3,2)$ & $A_{157}$ & $(2,2,1,1)$ & $(1,0,0,0)$ & $(1,5,1)$\\
$A_{113}$ & $(3,2,1,1)$ & $(0,1,0,0)$ & $(1,3,2)$ & $A_{158}$ & $(4,0,1,1)$ & $(0,1,0,0)$ & $(1,5,1)$\\
$A_{114}$ & $(3,4,1,1)$ & $(0,0,0,0)$ & $(1,3,2)$ & $A_{159}$ & $(4,2,1,1)$ & $(0,0,0,0)$ & $(1,5,1)$\\
$A_{115}$ & $(1,0,1,1)$ & $(1,2,0,0)$ & $(1,5,1)$ & $A_{160}$ & $(0,1,1,1)$ & $(2,1,0,0)$ & $(1,1,4)$\\
$A_{116}$ & $(1,2,1,1)$ & $(1,1,0,0)$ & $(1,5,1)$ & $A_{161}$ & $(0,3,1,1)$ & $(2,0,0,0)$ & $(1,1,4)$\\
$A_{117}$ & $(1,4,1,1)$ & $(1,0,0,0)$ & $(1,5,1)$ & $A_{162}$ & $(2,1,1,1)$ & $(1,1,0,0)$ & $(1,1,4)$\\
$A_{118}$ & $(3,0,1,1)$ & $(0,2,0,0)$ & $(1,5,1)$ & $A_{163}$ & $(2,3,1,1)$ & $(1,0,0,0)$ & $(1,1,4)$\\
$A_{119}$ & $(3,2,1,1)$ & $(0,1,0,0)$ & $(1,5,1)$ & $A_{164}$ & $(4,1,1,1)$ & $(0,1,0,0)$ & $(1,1,4)$\\
$A_{120}$ & $(3,4,1,1)$ & $(0,0,0,0)$ & $(1,5,1)$ & $A_{165}$ & $(4,3,1,1)$ & $(0,0,0,0)$ & $(1,1,4)$\\
$A_{121}$ & $(1,1,1,1)$ & $(1,2,0,0)$ & $(1,1,4)$ & $A_{166}$ & $(0,1,1,1)$ & $(2,1,0,0)$ & $(1,2,2)$\\
$A_{122}$ & $(1,3,1,1)$ & $(1,1,0,0)$ & $(1,1,4)$ & $A_{167}$ & $(0,3,1,1)$ & $(2,0,0,0)$ & $(1,2,2)$\\
$A_{123}$ & $(1,5,1,1)$ & $(1,0,0,0)$ & $(1,1,4)$ & $A_{168}$ & $(2,1,1,1)$ & $(1,1,0,0)$ & $(1,2,2)$\\
$A_{124}$ & $(3,1,1,1)$ & $(0,2,0,0)$ & $(1,1,4)$ & $A_{169}$ & $(2,3,1,1)$ & $(1,0,0,0)$ & $(1,2,2)$\\
$A_{125}$ & $(3,3,1,1)$ & $(0,1,0,0)$ & $(1,1,4)$ & $A_{170}$ & $(4,1,1,1)$ & $(0,1,0,0)$ & $(1,2,2)$\\
$A_{126}$ & $(3,5,1,1)$ & $(0,0,0,0)$ & $(1,1,4)$ & $A_{171}$ & $(4,3,1,1)$ & $(0,0,0,0)$ & $(1,2,2)$\\
$A_{127}$ & $(1,1,1,1)$ & $(1,2,0,0)$ & $(1,2,2)$ & $A_{172}$ & $(0,1,1,1)$ & $(2,1,0,0)$ & $(1,4,1)$\\
$A_{128}$ & $(1,3,1,1)$ & $(1,1,0,0)$ & $(1,2,2)$ & $A_{173}$ & $(0,3,1,1)$ & $(2,0,0,0)$ & $(1,4,1)$\\
$A_{129}$ & $(1,5,1,1)$ & $(1,0,0,0)$ & $(1,2,2)$ & $A_{174}$ & $(2,1,1,1)$ & $(1,1,0,0)$ & $(1,4,1)$\\
$A_{130}$ & $(3,1,1,1)$ & $(0,2,0,0)$ & $(1,2,2)$ & $A_{175}$ & $(2,3,1,1)$ & $(1,0,0,0)$ & $(1,4,1)$\\
$A_{131}$ & $(3,3,1,1)$ & $(0,1,0,0)$ & $(1,2,2)$ & $A_{176}$ & $(4,1,1,1)$ & $(0,1,0,0)$ & $(1,4,1)$\\
$A_{132}$ & $(3,5,1,1)$ & $(0,0,0,0)$ & $(1,2,2)$ & $A_{177}$ & $(4,3,1,1)$ & $(0,0,0,0)$ & $(1,4,1)$\\
$A_{133}$ & $(1,1,1,1)$ & $(1,2,0,0)$ & $(1,4,1)$ & $A_{178}$ & $(1,0,1,1)$ & $(2,1,0,0)$ & $(1,1,8)$\\
$A_{134}$ & $(1,3,1,1)$ & $(1,1,0,0)$ & $(1,4,1)$ & $A_{179}$ & $(1,2,1,1)$ & $(2,0,0,0)$ & $(1,1,8)$\\
$A_{135}$ & $(1,5,1,1)$ & $(1,0,0,0)$ & $(1,4,1)$ & $A_{180}$ & $(3,0,1,1)$ & $(1,1,0,0)$ & $(1,1,8)$\\
\bottomrule
\end{tabular*}
\endgroup
\end{table}
\begin{table}[p]
\centering
\begingroup\footnotesize
\setlength{\tabcolsep}{3pt}
\renewcommand{\arraystretch}{1.0}
\addtocounter{table}{-1}
\caption{(continued)}
\label{tab:s5-explicit-3}
\begin{tabular*}{\textwidth}{@{\extracolsep{\fill}}cccc@{\hspace{7mm}}cccc}
\toprule
$A_{i}$ & $\boldsymbol a$ & $\boldsymbol b$ & $\boldsymbol\mu$ & $A_{i}$ & $\boldsymbol a$ & $\boldsymbol b$ & $\boldsymbol\mu$\\
\midrule
$A_{181}$ & $(3,2,1,1)$ & $(1,0,0,0)$ & $(1,1,8)$ & $A_{226}$ & $(4,1,1,1)$ & $(2,0,0,0)$ & $(1,1,20)$\\
$A_{182}$ & $(5,0,1,1)$ & $(0,1,0,0)$ & $(1,1,8)$ & $A_{227}$ & $(6,1,1,1)$ & $(1,0,0,0)$ & $(1,1,20)$\\
$A_{183}$ & $(5,2,1,1)$ & $(0,0,0,0)$ & $(1,1,8)$ & $A_{228}$ & $(8,1,1,1)$ & $(0,0,0,0)$ & $(1,1,20)$\\
$A_{184}$ & $(0,1,1,1)$ & $(3,0,0,0)$ & $(1,2,4)$ & $A_{229}$ & $(0,0,1,1)$ & $(4,1,0,0)$ & $(1,1,2)$\\
$A_{185}$ & $(2,1,1,1)$ & $(2,0,0,0)$ & $(1,2,4)$ & $A_{230}$ & $(0,2,1,1)$ & $(4,0,0,0)$ & $(1,1,2)$\\
$A_{186}$ & $(4,1,1,1)$ & $(1,0,0,0)$ & $(1,2,4)$ & $A_{231}$ & $(2,0,1,1)$ & $(3,1,0,0)$ & $(1,1,2)$\\
$A_{187}$ & $(6,1,1,1)$ & $(0,0,0,0)$ & $(1,2,4)$ & $A_{232}$ & $(2,2,1,1)$ & $(3,0,0,0)$ & $(1,1,2)$\\
$A_{188}$ & $(0,1,1,1)$ & $(3,0,0,0)$ & $(1,3,2)$ & $A_{233}$ & $(4,0,1,1)$ & $(2,1,0,0)$ & $(1,1,2)$\\
$A_{189}$ & $(2,1,1,1)$ & $(2,0,0,0)$ & $(1,3,2)$ & $A_{234}$ & $(4,2,1,1)$ & $(2,0,0,0)$ & $(1,1,2)$\\
$A_{190}$ & $(4,1,1,1)$ & $(1,0,0,0)$ & $(1,3,2)$ & $A_{235}$ & $(6,0,1,1)$ & $(1,1,0,0)$ & $(1,1,2)$\\
$A_{191}$ & $(6,1,1,1)$ & $(0,0,0,0)$ & $(1,3,2)$ & $A_{236}$ & $(6,2,1,1)$ & $(1,0,0,0)$ & $(1,1,2)$\\
$A_{192}$ & $(0,1,1,1)$ & $(3,0,0,0)$ & $(1,5,1)$ & $A_{237}$ & $(8,0,1,1)$ & $(0,1,0,0)$ & $(1,1,2)$\\
$A_{193}$ & $(2,1,1,1)$ & $(2,0,0,0)$ & $(1,5,1)$ & $A_{238}$ & $(8,2,1,1)$ & $(0,0,0,0)$ & $(1,1,2)$\\
$A_{194}$ & $(4,1,1,1)$ & $(1,0,0,0)$ & $(1,5,1)$ & $A_{239}$ & $(0,0,1,1)$ & $(4,1,0,0)$ & $(1,3,1)$\\
$A_{195}$ & $(6,1,1,1)$ & $(0,0,0,0)$ & $(1,5,1)$ & $A_{240}$ & $(0,2,1,1)$ & $(4,0,0,0)$ & $(1,3,1)$\\
$A_{196}$ & $(0,0,1,1)$ & $(3,1,0,0)$ & $(1,1,4)$ & $A_{241}$ & $(2,0,1,1)$ & $(3,1,0,0)$ & $(1,3,1)$\\
$A_{197}$ & $(0,2,1,1)$ & $(3,0,0,0)$ & $(1,1,4)$ & $A_{242}$ & $(2,2,1,1)$ & $(3,0,0,0)$ & $(1,3,1)$\\
$A_{198}$ & $(2,0,1,1)$ & $(2,1,0,0)$ & $(1,1,4)$ & $A_{243}$ & $(4,0,1,1)$ & $(2,1,0,0)$ & $(1,3,1)$\\
$A_{199}$ & $(2,2,1,1)$ & $(2,0,0,0)$ & $(1,1,4)$ & $A_{244}$ & $(4,2,1,1)$ & $(2,0,0,0)$ & $(1,3,1)$\\
$A_{200}$ & $(4,0,1,1)$ & $(1,1,0,0)$ & $(1,1,4)$ & $A_{245}$ & $(6,0,1,1)$ & $(1,1,0,0)$ & $(1,3,1)$\\
$A_{201}$ & $(4,2,1,1)$ & $(1,0,0,0)$ & $(1,1,4)$ & $A_{246}$ & $(6,2,1,1)$ & $(1,0,0,0)$ & $(1,3,1)$\\
$A_{202}$ & $(6,0,1,1)$ & $(0,1,0,0)$ & $(1,1,4)$ & $A_{247}$ & $(8,0,1,1)$ & $(0,1,0,0)$ & $(1,3,1)$\\
$A_{203}$ & $(6,2,1,1)$ & $(0,0,0,0)$ & $(1,1,4)$ & $A_{248}$ & $(8,2,1,1)$ & $(0,0,0,0)$ & $(1,3,1)$\\
$A_{204}$ & $(0,0,1,1)$ & $(3,1,0,0)$ & $(1,2,2)$ & $A_{249}$ & $(1,1,1,1)$ & $(4,0,0,0)$ & $(1,1,14)$\\
$A_{205}$ & $(0,2,1,1)$ & $(3,0,0,0)$ & $(1,2,2)$ & $A_{250}$ & $(3,1,1,1)$ & $(3,0,0,0)$ & $(1,1,14)$\\
$A_{206}$ & $(2,0,1,1)$ & $(2,1,0,0)$ & $(1,2,2)$ & $A_{251}$ & $(5,1,1,1)$ & $(2,0,0,0)$ & $(1,1,14)$\\
$A_{207}$ & $(2,2,1,1)$ & $(2,0,0,0)$ & $(1,2,2)$ & $A_{252}$ & $(7,1,1,1)$ & $(1,0,0,0)$ & $(1,1,14)$\\
$A_{208}$ & $(4,0,1,1)$ & $(1,1,0,0)$ & $(1,2,2)$ & $A_{253}$ & $(9,1,1,1)$ & $(0,0,0,0)$ & $(1,1,14)$\\
$A_{209}$ & $(4,2,1,1)$ & $(1,0,0,0)$ & $(1,2,2)$ & $A_{254}$ & $(0,1,1,1)$ & $(5,0,0,0)$ & $(1,1,12)$\\
$A_{210}$ & $(6,0,1,1)$ & $(0,1,0,0)$ & $(1,2,2)$ & $A_{255}$ & $(2,1,1,1)$ & $(4,0,0,0)$ & $(1,1,12)$\\
$A_{211}$ & $(6,2,1,1)$ & $(0,0,0,0)$ & $(1,2,2)$ & $A_{256}$ & $(4,1,1,1)$ & $(3,0,0,0)$ & $(1,1,12)$\\
$A_{212}$ & $(0,0,1,1)$ & $(3,1,0,0)$ & $(1,4,1)$ & $A_{257}$ & $(6,1,1,1)$ & $(2,0,0,0)$ & $(1,1,12)$\\
$A_{213}$ & $(0,2,1,1)$ & $(3,0,0,0)$ & $(1,4,1)$ & $A_{258}$ & $(8,1,1,1)$ & $(1,0,0,0)$ & $(1,1,12)$\\
$A_{214}$ & $(2,0,1,1)$ & $(2,1,0,0)$ & $(1,4,1)$ & $A_{259}$ & $(10,1,1,1)$ & $(0,0,0,0)$ & $(1,1,12)$\\
$A_{215}$ & $(2,2,1,1)$ & $(2,0,0,0)$ & $(1,4,1)$ & $A_{260}$ & $(0,1,1,1)$ & $(5,0,0,0)$ & $(1,2,3)$\\
$A_{216}$ & $(4,0,1,1)$ & $(1,1,0,0)$ & $(1,4,1)$ & $A_{261}$ & $(2,1,1,1)$ & $(4,0,0,0)$ & $(1,2,3)$\\
$A_{217}$ & $(4,2,1,1)$ & $(1,0,0,0)$ & $(1,4,1)$ & $A_{262}$ & $(4,1,1,1)$ & $(3,0,0,0)$ & $(1,2,3)$\\
$A_{218}$ & $(6,0,1,1)$ & $(0,1,0,0)$ & $(1,4,1)$ & $A_{263}$ & $(6,1,1,1)$ & $(2,0,0,0)$ & $(1,2,3)$\\
$A_{219}$ & $(6,2,1,1)$ & $(0,0,0,0)$ & $(1,4,1)$ & $A_{264}$ & $(8,1,1,1)$ & $(1,0,0,0)$ & $(1,2,3)$\\
$A_{220}$ & $(1,1,1,1)$ & $(3,0,0,0)$ & $(1,1,36)$ & $A_{265}$ & $(10,1,1,1)$ & $(0,0,0,0)$ & $(1,2,3)$\\
$A_{221}$ & $(3,1,1,1)$ & $(2,0,0,0)$ & $(1,1,36)$ & $A_{266}$ & $(0,0,1,1)$ & $(5,1,0,0)$ & $(1,2,1)$\\
$A_{222}$ & $(5,1,1,1)$ & $(1,0,0,0)$ & $(1,1,36)$ & $A_{267}$ & $(0,2,1,1)$ & $(5,0,0,0)$ & $(1,2,1)$\\
$A_{223}$ & $(7,1,1,1)$ & $(0,0,0,0)$ & $(1,1,36)$ & $A_{268}$ & $(2,0,1,1)$ & $(4,1,0,0)$ & $(1,2,1)$\\
$A_{224}$ & $(0,1,1,1)$ & $(4,0,0,0)$ & $(1,1,20)$ & $A_{269}$ & $(2,2,1,1)$ & $(4,0,0,0)$ & $(1,2,1)$\\
$A_{225}$ & $(2,1,1,1)$ & $(3,0,0,0)$ & $(1,1,20)$ & $A_{270}$ & $(4,0,1,1)$ & $(3,1,0,0)$ & $(1,2,1)$\\
\bottomrule
\end{tabular*}
\endgroup
\end{table}
\begin{table}[p]
\centering
\begingroup\footnotesize
\setlength{\tabcolsep}{3pt}
\renewcommand{\arraystretch}{1.0}
\addtocounter{table}{-1}
\caption{(continued)}
\label{tab:s5-explicit-4}
\begin{tabular*}{\textwidth}{@{\extracolsep{\fill}}cccc@{\hspace{7mm}}cccc}
\toprule
$A_{i}$ & $\boldsymbol a$ & $\boldsymbol b$ & $\boldsymbol\mu$ & $A_{i}$ & $\boldsymbol a$ & $\boldsymbol b$ & $\boldsymbol\mu$\\
\midrule
$A_{271}$ & $(4,2,1,1)$ & $(3,0,0,0)$ & $(1,2,1)$ & $A_{316}$ & $(6,1,1,1)$ & $(5,0,0,0)$ & $(1,1,7)$\\
$A_{272}$ & $(6,0,1,1)$ & $(2,1,0,0)$ & $(1,2,1)$ & $A_{317}$ & $(8,1,1,1)$ & $(4,0,0,0)$ & $(1,1,7)$\\
$A_{273}$ & $(6,2,1,1)$ & $(2,0,0,0)$ & $(1,2,1)$ & $A_{318}$ & $(10,1,1,1)$ & $(3,0,0,0)$ & $(1,1,7)$\\
$A_{274}$ & $(8,0,1,1)$ & $(1,1,0,0)$ & $(1,2,1)$ & $A_{319}$ & $(12,1,1,1)$ & $(2,0,0,0)$ & $(1,1,7)$\\
$A_{275}$ & $(8,2,1,1)$ & $(1,0,0,0)$ & $(1,2,1)$ & $A_{320}$ & $(14,1,1,1)$ & $(1,0,0,0)$ & $(1,1,7)$\\
$A_{276}$ & $(10,0,1,1)$ & $(0,1,0,0)$ & $(1,2,1)$ & $A_{321}$ & $(16,1,1,1)$ & $(0,0,0,0)$ & $(1,1,7)$\\
$A_{277}$ & $(10,2,1,1)$ & $(0,0,0,0)$ & $(1,2,1)$ & $A_{322}$ & $(0,1,1,1)$ & $(11,0,0,0)$ & $(1,1,6)$\\
$A_{278}$ & $(1,1,1,1)$ & $(5,0,0,0)$ & $(1,1,10)$ & $A_{323}$ & $(2,1,1,1)$ & $(10,0,0,0)$ & $(1,1,6)$\\
$A_{279}$ & $(3,1,1,1)$ & $(4,0,0,0)$ & $(1,1,10)$ & $A_{324}$ & $(4,1,1,1)$ & $(9,0,0,0)$ & $(1,1,6)$\\
$A_{280}$ & $(5,1,1,1)$ & $(3,0,0,0)$ & $(1,1,10)$ & $A_{325}$ & $(6,1,1,1)$ & $(8,0,0,0)$ & $(1,1,6)$\\
$A_{281}$ & $(7,1,1,1)$ & $(2,0,0,0)$ & $(1,1,10)$ & $A_{326}$ & $(8,1,1,1)$ & $(7,0,0,0)$ & $(1,1,6)$\\
$A_{282}$ & $(9,1,1,1)$ & $(1,0,0,0)$ & $(1,1,10)$ & $A_{327}$ & $(10,1,1,1)$ & $(6,0,0,0)$ & $(1,1,6)$\\
$A_{283}$ & $(11,1,1,1)$ & $(0,0,0,0)$ & $(1,1,10)$ & $A_{328}$ & $(12,1,1,1)$ & $(5,0,0,0)$ & $(1,1,6)$\\
$A_{284}$ & $(0,1,1,1)$ & $(6,0,0,0)$ & $(1,1,9)$ & $A_{329}$ & $(14,1,1,1)$ & $(4,0,0,0)$ & $(1,1,6)$\\
$A_{285}$ & $(2,1,1,1)$ & $(5,0,0,0)$ & $(1,1,9)$ & $A_{330}$ & $(16,1,1,1)$ & $(3,0,0,0)$ & $(1,1,6)$\\
$A_{286}$ & $(4,1,1,1)$ & $(4,0,0,0)$ & $(1,1,9)$ & $A_{331}$ & $(18,1,1,1)$ & $(2,0,0,0)$ & $(1,1,6)$\\
$A_{287}$ & $(6,1,1,1)$ & $(3,0,0,0)$ & $(1,1,9)$ & $A_{332}$ & $(20,1,1,1)$ & $(1,0,0,0)$ & $(1,1,6)$\\
$A_{288}$ & $(8,1,1,1)$ & $(2,0,0,0)$ & $(1,1,9)$ & $A_{333}$ & $(22,1,1,1)$ & $(0,0,0,0)$ & $(1,1,6)$\\
$A_{289}$ & $(10,1,1,1)$ & $(1,0,0,0)$ & $(1,1,9)$ & $A_{334}$ & $(0,1,1,1)$ & $(19,0,0,0)$ & $(1,1,5)$\\
$A_{290}$ & $(12,1,1,1)$ & $(0,0,0,0)$ & $(1,1,9)$ & $A_{335}$ & $(2,1,1,1)$ & $(18,0,0,0)$ & $(1,1,5)$\\
$A_{291}$ & $(0,0,1,1)$ & $(6,1,0,0)$ & $(1,1,1)$ & $A_{336}$ & $(4,1,1,1)$ & $(17,0,0,0)$ & $(1,1,5)$\\
$A_{292}$ & $(0,2,1,1)$ & $(6,0,0,0)$ & $(1,1,1)$ & $A_{337}$ & $(6,1,1,1)$ & $(16,0,0,0)$ & $(1,1,5)$\\
$A_{293}$ & $(2,0,1,1)$ & $(5,1,0,0)$ & $(1,1,1)$ & $A_{338}$ & $(8,1,1,1)$ & $(15,0,0,0)$ & $(1,1,5)$\\
$A_{294}$ & $(2,2,1,1)$ & $(5,0,0,0)$ & $(1,1,1)$ & $A_{339}$ & $(10,1,1,1)$ & $(14,0,0,0)$ & $(1,1,5)$\\
$A_{295}$ & $(4,0,1,1)$ & $(4,1,0,0)$ & $(1,1,1)$ & $A_{340}$ & $(12,1,1,1)$ & $(13,0,0,0)$ & $(1,1,5)$\\
$A_{296}$ & $(4,2,1,1)$ & $(4,0,0,0)$ & $(1,1,1)$ & $A_{341}$ & $(14,1,1,1)$ & $(12,0,0,0)$ & $(1,1,5)$\\
$A_{297}$ & $(6,0,1,1)$ & $(3,1,0,0)$ & $(1,1,1)$ & $A_{342}$ & $(16,1,1,1)$ & $(11,0,0,0)$ & $(1,1,5)$\\
$A_{298}$ & $(6,2,1,1)$ & $(3,0,0,0)$ & $(1,1,1)$ & $A_{343}$ & $(18,1,1,1)$ & $(10,0,0,0)$ & $(1,1,5)$\\
$A_{299}$ & $(8,0,1,1)$ & $(2,1,0,0)$ & $(1,1,1)$ & $A_{344}$ & $(20,1,1,1)$ & $(9,0,0,0)$ & $(1,1,5)$\\
$A_{300}$ & $(8,2,1,1)$ & $(2,0,0,0)$ & $(1,1,1)$ & $A_{345}$ & $(22,1,1,1)$ & $(8,0,0,0)$ & $(1,1,5)$\\
$A_{301}$ & $(10,0,1,1)$ & $(1,1,0,0)$ & $(1,1,1)$ & $A_{346}$ & $(24,1,1,1)$ & $(7,0,0,0)$ & $(1,1,5)$\\
$A_{302}$ & $(10,2,1,1)$ & $(1,0,0,0)$ & $(1,1,1)$ & $A_{347}$ & $(26,1,1,1)$ & $(6,0,0,0)$ & $(1,1,5)$\\
$A_{303}$ & $(12,0,1,1)$ & $(0,1,0,0)$ & $(1,1,1)$ & $A_{348}$ & $(28,1,1,1)$ & $(5,0,0,0)$ & $(1,1,5)$\\
$A_{304}$ & $(12,2,1,1)$ & $(0,0,0,0)$ & $(1,1,1)$ & $A_{349}$ & $(30,1,1,1)$ & $(4,0,0,0)$ & $(1,1,5)$\\
$A_{305}$ & $(0,1,1,1)$ & $(7,0,0,0)$ & $(1,1,8)$ & $A_{350}$ & $(32,1,1,1)$ & $(3,0,0,0)$ & $(1,1,5)$\\
$A_{306}$ & $(2,1,1,1)$ & $(6,0,0,0)$ & $(1,1,8)$ & $A_{351}$ & $(34,1,1,1)$ & $(2,0,0,0)$ & $(1,1,5)$\\
$A_{307}$ & $(4,1,1,1)$ & $(5,0,0,0)$ & $(1,1,8)$ & $A_{352}$ & $(36,1,1,1)$ & $(1,0,0,0)$ & $(1,1,5)$\\
$A_{308}$ & $(6,1,1,1)$ & $(4,0,0,0)$ & $(1,1,8)$ & $A_{353}$ & $(38,1,1,1)$ & $(0,0,0,0)$ & $(1,1,5)$\\
$A_{309}$ & $(8,1,1,1)$ & $(3,0,0,0)$ & $(1,1,8)$ & $A_{354}$ & $(1,1,1,1,1)$ & $(1,0,0,0,0)$ & $(1,1,1,1)$\\
$A_{310}$ & $(10,1,1,1)$ & $(2,0,0,0)$ & $(1,1,8)$ & $A_{355}$ & $(3,1,1,1,1)$ & $(0,0,0,0,0)$ & $(1,1,1,1)$\\
$A_{311}$ & $(12,1,1,1)$ & $(1,0,0,0)$ & $(1,1,8)$ &  &  &  & \\
$A_{312}$ & $(14,1,1,1)$ & $(0,0,0,0)$ & $(1,1,8)$ &  &  &  & \\
$A_{313}$ & $(0,1,1,1)$ & $(8,0,0,0)$ & $(1,1,7)$ &  &  &  & \\
$A_{314}$ & $(2,1,1,1)$ & $(7,0,0,0)$ & $(1,1,7)$ &  &  &  & \\
$A_{315}$ & $(4,1,1,1)$ & $(6,0,0,0)$ & $(1,1,7)$ &  &  &  & \\
\bottomrule
\end{tabular*}
\endgroup
\end{table}
\begin{table}[p]
\centering
\begingroup\footnotesize
\setlength{\tabcolsep}{2pt}
\renewcommand{\arraystretch}{1.0}
\caption{Graphs $\mathcal F_{1}(n),\ldots,\mathcal F_{94}(n)$.}
\label{tab:s5-explicit-F}
\begin{tabular*}{\textwidth}{@{\extracolsep{\fill}}ccccc@{\hspace{5mm}}ccccc}
\toprule
$\mathcal F_{i}(n)$ & $\boldsymbol a$ & $\boldsymbol b$ & $\boldsymbol\mu$ & $\kappa=\tau$ & $\mathcal F_{i}(n)$ & $\boldsymbol a$ & $\boldsymbol b$ & $\boldsymbol\mu$ & $\kappa=\tau$\\
\midrule
$\mathcal F_{1}(n)$ & $(3,6,2)$ & $(0,0,0)$ & $(1,n)$ & $0$ & $\mathcal F_{48}(n)$ & $(0,3,0)$ & $(3,0,1)$ & $(1,n)$ & $0$\\
$\mathcal F_{2}(n)$ & $(3,6,0)$ & $(0,0,1)$ & $(1,n)$ & $0$ & $\mathcal F_{49}(n)$ & $(0,1,2)$ & $(3,1,0)$ & $(1,n)$ & $0$\\
$\mathcal F_{3}(n)$ & $(3,4,2)$ & $(0,1,0)$ & $(1,n)$ & $0$ & $\mathcal F_{50}(n)$ & $(0,1,0)$ & $(3,1,1)$ & $(1,n)$ & $0$\\
$\mathcal F_{4}(n)$ & $(3,4,0)$ & $(0,1,1)$ & $(1,n)$ & $0$ & $\mathcal F_{51}(n)$ & $(3,4,1,2)$ & $(0,0,0,0)$ & $(1,1,n)$ & $0$\\
$\mathcal F_{5}(n)$ & $(3,2,2)$ & $(0,2,0)$ & $(1,n)$ & $0$ & $\mathcal F_{52}(n)$ & $(3,4,1,0)$ & $(0,0,0,1)$ & $(1,1,n)$ & $0$\\
$\mathcal F_{6}(n)$ & $(3,2,0)$ & $(0,2,1)$ & $(1,n)$ & $0$ & $\mathcal F_{53}(n)$ & $(3,2,1,2)$ & $(0,1,0,0)$ & $(1,1,n)$ & $0$\\
$\mathcal F_{7}(n)$ & $(3,0,2)$ & $(0,3,0)$ & $(1,n)$ & $0$ & $\mathcal F_{54}(n)$ & $(3,2,1,0)$ & $(0,1,0,1)$ & $(1,1,n)$ & $0$\\
$\mathcal F_{8}(n)$ & $(3,0,0)$ & $(0,3,1)$ & $(1,n)$ & $0$ & $\mathcal F_{55}(n)$ & $(3,0,1,2)$ & $(0,2,0,0)$ & $(1,1,n)$ & $0$\\
$\mathcal F_{9}(n)$ & $(1,6,2)$ & $(1,0,0)$ & $(1,n)$ & $0$ & $\mathcal F_{56}(n)$ & $(3,0,1,0)$ & $(0,2,0,1)$ & $(1,1,n)$ & $0$\\
$\mathcal F_{10}(n)$ & $(1,6,0)$ & $(1,0,1)$ & $(1,n)$ & $0$ & $\mathcal F_{57}(n)$ & $(1,4,1,2)$ & $(1,0,0,0)$ & $(1,1,n)$ & $0$\\
$\mathcal F_{11}(n)$ & $(1,4,2)$ & $(1,1,0)$ & $(1,n)$ & $0$ & $\mathcal F_{58}(n)$ & $(1,4,1,0)$ & $(1,0,0,1)$ & $(1,1,n)$ & $0$\\
$\mathcal F_{12}(n)$ & $(1,4,0)$ & $(1,1,1)$ & $(1,n)$ & $0$ & $\mathcal F_{59}(n)$ & $(1,2,1,2)$ & $(1,1,0,0)$ & $(1,1,n)$ & $0$\\
$\mathcal F_{13}(n)$ & $(1,2,2)$ & $(1,2,0)$ & $(1,n)$ & $0$ & $\mathcal F_{60}(n)$ & $(1,2,1,0)$ & $(1,1,0,1)$ & $(1,1,n)$ & $0$\\
$\mathcal F_{14}(n)$ & $(1,2,0)$ & $(1,2,1)$ & $(1,n)$ & $0$ & $\mathcal F_{61}(n)$ & $(1,0,1,2)$ & $(1,2,0,0)$ & $(1,1,n)$ & $0$\\
$\mathcal F_{15}(n)$ & $(1,0,2)$ & $(1,3,0)$ & $(1,n)$ & $0$ & $\mathcal F_{62}(n)$ & $(1,0,1,0)$ & $(1,2,0,1)$ & $(1,1,n)$ & $0$\\
$\mathcal F_{16}(n)$ & $(1,0,0)$ & $(1,3,1)$ & $(1,n)$ & $0$ & $\mathcal F_{63}(n)$ & $(4,2,1,2)$ & $(0,0,0,0)$ & $(1,1,n)$ & $0$\\
$\mathcal F_{17}(n)$ & $(4,4,2)$ & $(0,0,0)$ & $(1,n)$ & $0$ & $\mathcal F_{64}(n)$ & $(4,2,1,0)$ & $(0,0,0,1)$ & $(1,1,n)$ & $0$\\
$\mathcal F_{18}(n)$ & $(4,4,0)$ & $(0,0,1)$ & $(1,n)$ & $0$ & $\mathcal F_{65}(n)$ & $(4,0,1,2)$ & $(0,1,0,0)$ & $(1,1,n)$ & $0$\\
$\mathcal F_{19}(n)$ & $(4,2,2)$ & $(0,1,0)$ & $(1,n)$ & $0$ & $\mathcal F_{66}(n)$ & $(4,0,1,0)$ & $(0,1,0,1)$ & $(1,1,n)$ & $0$\\
$\mathcal F_{20}(n)$ & $(4,2,0)$ & $(0,1,1)$ & $(1,n)$ & $0$ & $\mathcal F_{67}(n)$ & $(2,2,1,2)$ & $(1,0,0,0)$ & $(1,1,n)$ & $0$\\
$\mathcal F_{21}(n)$ & $(4,0,2)$ & $(0,2,0)$ & $(1,n)$ & $0$ & $\mathcal F_{68}(n)$ & $(2,2,1,0)$ & $(1,0,0,1)$ & $(1,1,n)$ & $0$\\
$\mathcal F_{22}(n)$ & $(4,0,0)$ & $(0,2,1)$ & $(1,n)$ & $0$ & $\mathcal F_{69}(n)$ & $(2,0,1,2)$ & $(1,1,0,0)$ & $(1,1,n)$ & $0$\\
$\mathcal F_{23}(n)$ & $(2,4,2)$ & $(1,0,0)$ & $(1,n)$ & $0$ & $\mathcal F_{70}(n)$ & $(2,0,1,0)$ & $(1,1,0,1)$ & $(1,1,n)$ & $0$\\
$\mathcal F_{24}(n)$ & $(2,4,0)$ & $(1,0,1)$ & $(1,n)$ & $0$ & $\mathcal F_{71}(n)$ & $(0,2,1,2)$ & $(2,0,0,0)$ & $(1,1,n)$ & $0$\\
$\mathcal F_{25}(n)$ & $(2,2,2)$ & $(1,1,0)$ & $(1,n)$ & $0$ & $\mathcal F_{72}(n)$ & $(0,2,1,0)$ & $(2,0,0,1)$ & $(1,1,n)$ & $0$\\
$\mathcal F_{26}(n)$ & $(2,2,0)$ & $(1,1,1)$ & $(1,n)$ & $0$ & $\mathcal F_{73}(n)$ & $(0,0,1,2)$ & $(2,1,0,0)$ & $(1,1,n)$ & $0$\\
$\mathcal F_{27}(n)$ & $(2,0,2)$ & $(1,2,0)$ & $(1,n)$ & $0$ & $\mathcal F_{74}(n)$ & $(0,0,1,0)$ & $(2,1,0,1)$ & $(1,1,n)$ & $0$\\
$\mathcal F_{28}(n)$ & $(2,0,0)$ & $(1,2,1)$ & $(1,n)$ & $0$ & $\mathcal F_{75}(n)$ & $(6,1,1,2)$ & $(0,0,0,0)$ & $(1,1,n)$ & $0$\\
$\mathcal F_{29}(n)$ & $(0,4,2)$ & $(2,0,0)$ & $(1,n)$ & $0$ & $\mathcal F_{76}(n)$ & $(6,1,1,0)$ & $(0,0,0,1)$ & $(1,1,n)$ & $0$\\
$\mathcal F_{30}(n)$ & $(0,4,0)$ & $(2,0,1)$ & $(1,n)$ & $0$ & $\mathcal F_{77}(n)$ & $(4,1,1,2)$ & $(1,0,0,0)$ & $(1,1,n)$ & $0$\\
$\mathcal F_{31}(n)$ & $(0,2,2)$ & $(2,1,0)$ & $(1,n)$ & $0$ & $\mathcal F_{78}(n)$ & $(4,1,1,0)$ & $(1,0,0,1)$ & $(1,1,n)$ & $0$\\
$\mathcal F_{32}(n)$ & $(0,2,0)$ & $(2,1,1)$ & $(1,n)$ & $0$ & $\mathcal F_{79}(n)$ & $(2,1,1,2)$ & $(2,0,0,0)$ & $(1,1,n)$ & $0$\\
$\mathcal F_{33}(n)$ & $(0,0,2)$ & $(2,2,0)$ & $(1,n)$ & $0$ & $\mathcal F_{80}(n)$ & $(2,1,1,0)$ & $(2,0,0,1)$ & $(1,1,n)$ & $0$\\
$\mathcal F_{34}(n)$ & $(0,0,0)$ & $(2,2,1)$ & $(1,n)$ & $0$ & $\mathcal F_{81}(n)$ & $(0,1,1,2)$ & $(3,0,0,0)$ & $(1,1,n)$ & $0$\\
$\mathcal F_{35}(n)$ & $(6,3,2)$ & $(0,0,0)$ & $(1,n)$ & $0$ & $\mathcal F_{82}(n)$ & $(0,1,1,0)$ & $(3,0,0,1)$ & $(1,1,n)$ & $0$\\
$\mathcal F_{36}(n)$ & $(6,3,0)$ & $(0,0,1)$ & $(1,n)$ & $0$ & $\mathcal F_{83}(n)$ & $(0,2,2)$ & $(0,0,0)$ & $(1,n)$ & $n$\\
$\mathcal F_{37}(n)$ & $(6,1,2)$ & $(0,1,0)$ & $(1,n)$ & $0$ & $\mathcal F_{84}(n)$ & $(0,2,0)$ & $(0,0,1)$ & $(1,n)$ & $n$\\
$\mathcal F_{38}(n)$ & $(6,1,0)$ & $(0,1,1)$ & $(1,n)$ & $0$ & $\mathcal F_{85}(n)$ & $(0,0,2)$ & $(0,1,0)$ & $(1,n)$ & $n$\\
$\mathcal F_{39}(n)$ & $(4,3,2)$ & $(1,0,0)$ & $(1,n)$ & $0$ & $\mathcal F_{86}(n)$ & $(0,0,0)$ & $(0,1,1)$ & $(1,n)$ & $n$\\
$\mathcal F_{40}(n)$ & $(4,3,0)$ & $(1,0,1)$ & $(1,n)$ & $0$ & $\mathcal F_{87}(n)$ & $(0,1,2)$ & $(0,0,0)$ & $(2,n)$ & $n$\\
$\mathcal F_{41}(n)$ & $(4,1,2)$ & $(1,1,0)$ & $(1,n)$ & $0$ & $\mathcal F_{88}(n)$ & $(0,1,0)$ & $(0,0,1)$ & $(2,n)$ & $n$\\
$\mathcal F_{42}(n)$ & $(4,1,0)$ & $(1,1,1)$ & $(1,n)$ & $0$ & $\mathcal F_{89}(n)$ & $(0,1,1)$ & $(0,0,0)$ & $(3,4)$ & $n$\\
$\mathcal F_{43}(n)$ & $(2,3,2)$ & $(2,0,0)$ & $(1,n)$ & $0$ & $\mathcal F_{90}(n)$ & $(0,1,1)$ & $(0,0,0)$ & $(4,2)$ & $n$\\
$\mathcal F_{44}(n)$ & $(2,3,0)$ & $(2,0,1)$ & $(1,n)$ & $0$ & $\mathcal F_{91}(n)$ & $(0,1,1)$ & $(0,0,0)$ & $(6,1)$ & $n$\\
$\mathcal F_{45}(n)$ & $(2,1,2)$ & $(2,1,0)$ & $(1,n)$ & $0$ & $\mathcal F_{92}(n)$ & $(0,1,1,1)$ & $(0,0,0,0)$ & $(1,1,4)$ & $n$\\
$\mathcal F_{46}(n)$ & $(2,1,0)$ & $(2,1,1)$ & $(1,n)$ & $0$ & $\mathcal F_{93}(n)$ & $(0,1,1,1)$ & $(0,0,0,0)$ & $(1,2,2)$ & $n$\\
$\mathcal F_{47}(n)$ & $(0,3,2)$ & $(3,0,0)$ & $(1,n)$ & $0$ & $\mathcal F_{94}(n)$ & $(0,1,1,1)$ & $(0,0,0,0)$ & $(1,4,1)$ & $n$\\
\bottomrule
\end{tabular*}
\endgroup
\end{table}

\end{document}